\documentclass[11pt, reqno]{amsart}
\usepackage{indentfirst}
\usepackage{amsmath,amssymb,latexsym,esint,cite,mathrsfs}
\usepackage{verbatim,wasysym}
\usepackage[left=2.4cm,right=2.4cm,top=2.6cm,bottom=2.6cm]{geometry}
\usepackage{tikz,enumitem,graphicx, subfig, microtype, color}
\usepackage{epic,eepic}
\usepackage[colorlinks=true,urlcolor=blue, citecolor=red,linkcolor=blue,
linktocpage,pdfpagelabels, bookmarksnumbered,bookmarksopen]{hyperref}
\usepackage[hyperpageref]{backref}
\usepackage[english]{babel}
\numberwithin{equation}{section}
\newtheorem{thm}{Theorem}[section]
\newtheorem{lem}[thm]{Lemma}

\newtheorem{Prop}[thm]{Proposition}
\newtheorem{Def}[thm]{Definition}
\newtheorem{Rem}[thm]{Remark}

\begin{document}
\title[Maxwell equation]{On a critical time-harmonic Maxwell equation in nonlocal media}
\author[M.\ Yang]{Minbo Yang}
\author[W.\ Ye]{Weiwei Ye}
\author[S.\ Zhang]{Shuijin Zhang}
\address{Minbo Yang, Shuijin Zhang \newline\indent Department of Mathematics, Zhejiang Normal University, \newline\indent
	Jinhua, Zhejiang, 321004, People's Republic of China}
\email{M. Yang: mbyang@zjnu.edu.cn;  S. Zhang: shuijinzhang@zjnu.edu.cn }
	\address{Weiwei Ye \newline\indent Department of Mathematics, Zhejiang Normal University, \newline\indent
	Jinhua, Zhejiang, 321004, People's Republic of China, \newline\indent
	Department of Mathematics, Fuyang Normal University, \newline\indent
	Fuyang, Anhui, 236037, People's Republic of China}
	\email{W. Ye: yeweiweime@163.com}
.
  
\begin{abstract}
In this paper, we study the existence of solutions for a critical time-harmonic Maxwell equation in nonlocal media
\begin{equation*}
\left\{
\begin{aligned}
&\nabla\times(\nabla\times u)+\lambda u=\left(I_{\alpha}\ast|u|^{2^{\ast}_{\alpha}}\right)|u|^{2^{\ast}_{\alpha}-2}u~~&&\mathrm{in}~\Omega,\\
&\nu\times u=0~~&&\mathrm{on}~\partial\Omega,
\end{aligned}
\right.
\end{equation*}
where $\Omega\subset\mathbb{R}^{3}$ is a bounded domain, either convex or with $\mathcal{C}^{1,1}$ boundary, $\nu$ is the exterior normal, $\lambda<0$ is a real parameter, $2^{\ast}_{\alpha}=3+\alpha$ with $0<\alpha<3$ is the upper critical exponent due to the Hardy-Littlewood-Sobolev inequality. By introducing some suitable Coulomb spaces involving curl operator $W^{\alpha,2^{\ast}_{\alpha}}_{0}(\mathrm{curl};\Omega)$, we are able to obtain the ground state solutions of the curl-curl equation via the method of constraining Nehari-Pankov manifold. Correspondingly, some sharp constants of the Sobolev-like inequalities with curl operator are obtained by a nonlocal version of the concentration-compactness principle.
\end{abstract}
\maketitle
%
%\medskip
\setlength{\parindent}{2em}
\section{\bf  Introduction and Main Results}
\subsection{Introduction}

Let $\Omega\subset\mathbb{R}^{3}$ be a bounded domain, we are concerned with the curl-curl equation
\begin{equation}\label{bounded domain case}
\left\{
\begin{aligned}
&\nabla\times(\nabla\times u)+\lambda u=f(x,u)~~&&\mathrm{in}~\Omega,\\
&\nu\times u=0~~&&\mathrm{on}~\partial\Omega,
\end{aligned}
\right.
\end{equation}
where $\lambda<0$ is a real parameter, $\nu:\partial\Omega\longrightarrow\mathbb{R}^{3}$ is the exterior normal. Equation (\ref{bounded domain case}) can be derived from the first order Maxwell equation \cite{Monk2003}
\begin{equation}\label{First order Maxwell in differential form}
\left \{
\begin{aligned}
&\nabla\times \mathcal{H}=\mathcal{J}+\partial_{t}\mathcal{D}, ~&&(\mathrm{Ampere's~circuital~law})\\
& div(\mathcal{D})=\varrho, ~ &&(\mathrm{Gauss's~ law})\\
&\partial_{t}\mathcal{B}+\nabla\times \mathcal{E}=0,  ~ && (\mathrm{Faraday's~  law~ of ~induction})\\
&div(\mathcal{B})=0, ~ &&(\mathrm{Gauss's ~law~ for~ magnetism})
\end{aligned}
\right.
\end{equation}
where $\mathcal{E},\mathcal{H},\mathcal{D},\mathcal{B}$ are corresponded to the electric field, magnetic induction, electric displacement and magnetic filed, respectively. $\mathcal{J}$ is the electric current intensity, and $\varrho$ is the electric charge density. Generally, these physical quantities satisfy the following constitutive equations (see \cite[Section 1.1.3]{Dorfler1986}):
\begin{equation}\label{constitutive equations}
\mathcal{\mathcal{J}}=\sigma \mathcal{E},~\mathcal{D}=\varepsilon \mathcal{E}+\mathcal{P}_{NL},~\mathcal{H}=\frac{1}{\mu}\mathcal{B}-\mathcal{M},
\end{equation}
where $\mathcal{P}_{NL}, \mathcal{M}$ denote the polarization field and magnetization filed respectively, $\varepsilon,\mu,\sigma$ are the electric permittivity, magnetic permeability and the electric conductivity.
Taking the special case with the absence of charges, currents and magnetization, namely, $\mathcal{J}=\mathcal{M}=0$, $\rho=0$, equation (\ref{First order Maxwell in differential form}) becomes the second curl-curl equation
\begin{equation}\label{second}
\nabla\times(\frac{1}{\mu}\nabla\times\mathcal{E})+\varepsilon\partial_{t}^{2}\mathcal{E}=-\partial_{t}^{2}\mathcal{P}_{NL}.
\end{equation}
As the electric field and polarization field are time harmonic with the ansatz $\mathcal{E}(x,t)=E(x)e^{i\omega t}$, $\mathcal{P}_{NL}(x,t)=P(x)e^{i\omega t}$, equation (\ref{second}) turns into the time-harmonic Maxwell equation
\begin{equation*}
\nabla\times(\frac{1}{\mu}\nabla\times E)+\varepsilon \omega^{2}E=\omega^{2} P.
\end{equation*}
In some Kerr-like medias, the polarization field function $\mathcal{P}_{NL}$  is usually chosen to be
 $\mathcal{P}_{NL}=\alpha(x)|\mathcal{E}|^{p-2}\mathcal{E}$ with $2\leq p\leq6$ for the purpose of simplifying  the model. Then by setting
\begin{equation*}
f(x,E)=\partial_{E}F(x,E)=\mu\omega^{2}\alpha(x)|E|^{p-2}E,
\end{equation*}
one can deduce the main equation (\ref{bounded domain case})
\begin{equation*}
\nabla\times(\nabla\times E)+\lambda E=f(x,E),
\end{equation*}
where $\lambda=-\mu\omega^{2}\varepsilon$. The boundary condition holds when $\Omega$ is surrounded by a perfect conductor, see \cite{Dorfler1986}.

Apparently, equation (\ref{bounded domain case}) has a variational structure and the solutions are the critical points of the functional
 \begin{equation}\label{functional 1}
J_{\lambda}(u)=\int_{\Omega}|\nabla\times u|^{2}dx+\frac{\lambda}{2}\int_{\Omega}|u|^{2}dx-\int_{\Omega}F(x,u)dx,
\end{equation}
which is well defined on the natural space
\begin{equation*}
X=W^{p}_{0}(\mathrm{curl};\Omega)=\overline{\mathcal{C}_{0}^{\infty}(\Omega,\mathbb{R}^{3})}^{||\cdot||_{W^{p}(\mathrm{curl};\Omega)}},
\end{equation*}
where
\begin{equation*}
W^{p}(\mathrm{curl};\Omega)=\{u\in L^{p}(\Omega,\mathbb{R}^{3}):\nabla\times u\in L^{2}(\Omega,\mathbb{R}^{3})\}
\end{equation*}
is a Banach space, see \cite{Bartsch2015-1}. By introducing the Helmholtz decomposition
\begin{equation*}
W^{p}_{0}(\mathrm{curl};\Omega)=X_{\Omega}\oplus X^c_{\Omega},
\end{equation*}
where
\begin{equation*}
\begin{aligned}
X_{\Omega}:&=\{v\in W_{0}^{p}(\mathrm{curl};\Omega):\int_{\Omega}\langle v,\varphi\rangle dx=0~\mathrm{for~every}~\varphi\in C^{\infty}_{0}(\Omega,\mathbb{R}^{3})~\mathrm{with}~\nabla\times\varphi=0\}\\
&=\{v\in W_{0}^{p}(\mathrm{curl};\Omega):\mathrm{div}~v=0~\mathrm{in~the~sense~of~distributions}\},
\end{aligned}
\end{equation*}
and
\begin{equation*}\label{W Omega}
X^c_{\Omega}:=\{w\in W^{p}_{0}(\mathrm{curl};\Omega):\int_{\Omega}\langle w,\nabla\times\varphi\rangle dx=0~\mathrm{for~all}~\varphi\in C^{\infty}_{0}(\Omega,\mathbb{R}^{3})\},
\end{equation*}
functional (\ref{functional 1}) can be rewritten as
\begin{equation*}\label{functional2}
\begin{aligned}
J_{\lambda}(u)=J_{\lambda}(v+w)&=\frac{1}{2}\int_{\Omega}|\nabla\times v|^{2}dx+\frac{\lambda}{2}\int_{
\Omega}|v+w|^{2}dx-\int_{\Omega}F(x,v+w)dx\\
&=\frac{1}{2}\int_{\Omega}|\nabla v|^{2}dx+\frac{\lambda}{2}\int_{
\Omega}|v+w|^{2}dx-\int_{\Omega}F(x,v+w)dx,
\end{aligned}
\end{equation*}
where $\nabla\times(\nabla\times v)=\nabla(\nabla\cdot v)-\nabla\cdot(\nabla v)=-\Delta v$ for $\mathrm{div}v=0$. Since the operator $\nabla\times(\nabla\times\cdot)$ has an infinite dimension kernel, i.e. $\nabla\times(\nabla\varphi)=0$ for $\varphi\in \mathcal{C}^{\infty}_{0}(\Omega)$, one can easily check that $J_{\lambda}$ has the strongly indefinite nature. Particularly, $J_{\lambda}$ has a linking geometry as $\lambda\leq0$, see \cite[pp.3]{Bartsch2015-1}.

To overcome the difficulty of the strong indefiniteness, by assuming that the additional condition $\nabla\cdot u=0$, then the curl-curl equation can be reduced to the classical elliptic equation
\begin{equation}\label{elliptic}
\left\{
\begin{aligned}
&-\Delta u+\lambda u=f(x,u)~~&&\mathrm{in}~\Omega,\\
&u=0~~&&\mathrm{on}~\partial\Omega.
\end{aligned}
\right.
\end{equation}
%where $\nabla\times(\nabla\times u)=\nabla(\nabla\cdot u)-\nabla\cdot(\nabla u)=\nabla(\nabla\cdot u)-\Delta u=-\Delta u$.
This elliptic equation has been widely studied in different dimensional spaces and topology regions, see the pioneering work of Brezis and Nirenberg \cite{Brezis1983} and the more references \cite{Capozzi1985,Cerami1986,Gao2018}. Essentially, the non-divergence condition is a Coulomb gauge condition, which holds in the gauge invariant field. This requires that the polarization field $\mathcal{P}_{NL}$ does not happen or linearly depends on $\mathcal{E}$, otherwise, it destroys the gauge invariance of the curl-curl equation. If the polarization field $\mathcal{P}_{NL}=0$, then the curl-curl equation becomes a linear time harmonic Maxwell equation, which has been extensively considered in \cite{Dorfler1986,Monk2003,Picard2001}. Physically there indeed exist some special cylindrically symmetric transverse electric and transverse magnetic which satisfy the non-divergence condition, and they have been studied by Stuart and Zhou in \cite{Stuart2001,Stuart2005}.

For the general case with $\nabla\cdot u\neq0$, the study of the curl-curl equation becomes much more challenge. The first attempt goes back to the pioneering work of Benci \cite{Benci2004}. Under some nonlinear assumptions on $W(t)$, the authors investigated the Born-Infeld static magnetic model
\begin{equation}\label{Born}
\nabla\times(\nabla\times A)=W'(|A|^{2})A,~~\mathrm{in}~\mathbb{R}^{3},
\end{equation}
where $A=\nabla\times B$ is a magnetic potential. In a suitable subspace, Azzollini et al. \cite{Azzollini2006} obtained the cylindrically symmetric solutions of (\ref{Born}) by the  Palais principle of symmetric criticality. By using the Hodge decomposition, the cylindrically symmetric solutions with a second form have also been constructed by D'Aprile and Siciliano in \cite{D'Aprile2011}. In fact, in some bounded domains with cylindrically symmetric, the similar solutions were obtained in \cite{Bartsch2015-1,Bartsch2017,Mederski2018}. Barstch et al. \cite{Bartsch2016} also analysed the spectrum of the curl-curl operator with cylindrically symmetric period potential $V(x)=V(r,x_{3})$, and considered the following time-harmonic Maxwell equation
\begin{equation}\label{Barstch}
\nabla\times(\nabla\times E)+V(x)E=\Gamma(x)|E|^{p-2}E,~~\mathrm{in}~\mathbb{R}^{3},
\end{equation}
where $\Gamma(x)$ is a period function with respect $x_{3}$. By the method of constraining symmetric sub-manifold, the cylindrically symmetric ground state solutions of (\ref{Barstch}) were obtained, one may see \cite{Zeng2017} for other extended results.

If the problem was set in some non-symmetric bounded domains or some cases with non-symmetric potential, the methods mentioned above do not work well. Moreover, due to the lack of weak-weak$^{\ast}$ continuity of $J^{'}_{\lambda}(u)$ , the abstract linking theorems established in \cite{Benci1979,Mederski2016,Guo2016} do not work any longer, and so we fail to look for the suitable $(PS)$ sequences. Even if we can obtain the bounded $(PS)$ sequence, we still do not know whether the weak limit is a critical point of the functional. Inspired by the work of Szulkin and Weth in \cite{Szulkin2009}, Barstch and  Mederski\cite{Bartsch2015-1} constructed a Nehari-Pankov manifold, which is homeomorphism with the upper unit ball of the subspace $X_{\Omega}$, and in where, the $(PS)$ sequence is obtained by the Ekeland variational principle. On the other hand, by the compact embedding
\begin{equation}\label{embed}
X_{\Omega}\hookrightarrow L^{p}(\Omega,\mathbb{R}^{3}),~~2\leq p<6,
\end{equation}
they succeeded in verifying the $(PS)_{c}^{\tau}$ condition, see Definition \ref{PSC}, which implies the weak-weak$^{\ast}$ continuity of $J^{'}_{\lambda}(u)$. We would also like to mention that the convexity assumption of the nonlinearity $f(x,u)$ plays a key role in finding the bounded $(PS)$ sequence, see also \cite{Bartsch2017} for the weakened version. For other related results, we may turn to \cite{Qin2016-1} for the asymptotically linear case and \cite{Mederski2019} for the case with supercritical growth at $0$ and subcritical growth at infinity.

For the critical case $p=6$, the embedding (\ref{embed}) above is not compact any more, then it is rather difficult to verify the $(\mathrm{PS})^{\tau}_{c}$ condition. Mederski \cite{Mederski2018} proposed a compactly perturbed method and proved that the $(\mathrm{PS})$ sequence contains a weakly convergence subsequence with a nontrivial limit point. Later, Mederski and Szulkin\cite{Mederski2021} established a general concentration-compactness lemma in $\mathbb{R}^{N}$ and obtained the sharp constant in the curl inequality. As an application, the authors dealt with the Brezis-Nirenberg type problem by an extend skill. In the entire space $\mathbb{R}^{3}$, the embedding above is also not compact, then a new critical point theory related to a new topological manifold has been established by Mederski et al. in \cite{Mederski2019}, there the compactness was recovered and the existence of multiply solutions was obtained. In a direct way, Mederski \cite{Mederski2016} established a global compactness lemma that accounts for the lack of weak-weak$^{\ast}$ continuity. Moreover, a Pohozaev identity has been established, which gives a criterion for the nonexistence of classical solution. In an earlier work, Bartsch \cite{Bartsch2016} showed that no interesting solution can be leaded under the fully radial symmetry assumption on the potential $V(x)$.

However, for some Kerr-type nonlinear mediums, the material law (\ref{constitutive equations}) between the electric field $\mathcal{E}$ and the displacement field $\mathcal{D}$ becomes more delicate, see \cite[(1.8)]{Bartsch2016},
\begin{equation*}	\mathcal{D}=\epsilon_{0}(n(x)^{2}\mathcal{E}+\mathcal{P}_{NL}(x,\mathcal{E}))~~\mathrm{with}~~\mathcal{P}_{NL}(x,\mathcal{E})=\chi^{(3)}(\mathcal{E}\cdot\mathcal{E})\mathcal{E},
\end{equation*}
where $n^{2}(x)=1+\chi^{(1)}(x)$ is the square of the refractive index and $\chi^{(1)}$, $\chi^{(3)}$ denote the linear and cubic susceptibilities of the medium respectively. Particularly, in some nonlocal optic materials, the refractive index $n(x)$ is quite dependent on the electric field $\mathcal{E}$ in a small neighbourhood, and the refractive index change $\triangle n$ can be represented in general form as
\begin{equation*}
	\triangle n(E)=s\int_{-\infty}^{+\infty}K(x-y)E(x)dx,
\end{equation*}
see \cite[(1)]{Krolikowski2001}. This phenomenological model is of great significance in the research  of laser beams and solitary waves in nonlocal nematic liquid crystals, see \cite{Krolikowski2001,Reimbert2006} and the reference therein. However, these articles are based on the nonlinear Schr\"{o}dinger equation, which is an asymptotic approximation of Maxwell's equations. To investigate more information about the electromagnetic waves in the nonlocal optic mediums, one need to deal with the full three-dimensional Maxwell problem.
Recently, Mandel \cite{Mandel2021} investigated the curl-curl equation with nonlocal nonlinearity
\begin{equation*}
\nabla\times(\nabla\times E)+E=(K(x)\ast|E|^{p})|E|^{p-2}E~~\mathrm{in}~\mathbb{R}^{3},
\end{equation*}
and the author proved that nonlocal media may admit ground states even though the corresponding local models do not admit. In there, the parameter $\lambda=-\mu\omega^{2}\varepsilon=1$ with $\varepsilon<0$ is corresponded to the new artificially produced metamaterials with negative reflexive, see \cite{Padilla2006}, and the kernel $K(x)=e^{-f(x)}$ is a exponent type responding function which expresses the nonlocal polarization of the nonlocal optical media, see \cite{Krolikowski2001} and \cite{Nikolov2003} for more cases with oscillatory kernel functions. What's more, nonlocality appears naturally in optical systems with a thermal \cite{L} and it is known to influence the propagation of electromagnetic
waves in plasmas \cite{BC}.  Nonlocality also has
attracted considerable interest as a means of eliminating
collapse and stabilizing multidimensional solitary waves \cite{Ba} and it plays an important
role in the theory of Bose-Einstein condensation \cite{D} where it accounts for the finite-range many-body interactions.

\subsection{Main Results}

In the present paper, we are interested in the curl-curl equation with critical convolution part, namely, we consider the curl-curl equation with Riesz potential interaction part
\begin{equation*}
\nabla\times(\nabla\times E)+\lambda E=(I_{\alpha}(x)\ast|E|^{p})|E|^{p-2}E~~\mathrm{in}~\mathbb{R}^{3},
\end{equation*}
where $I_{\alpha}:\mathbb{R}^{3}\longrightarrow R$ is the Riesz potential of order $\alpha\in(0,3)$
defined for $x\in\mathbb{R}^{3}\setminus\{0\}$ as
\begin{equation*}
I_{\alpha}(x)=\frac{A_{\alpha}}{|x|^{3-\alpha}},~~A^{\alpha}=\frac{\Gamma(\frac{3-\alpha}{2})}{\Gamma(\frac{\alpha}{2})\pi^{\frac{N}{2}}2^{\alpha}}.
\end{equation*}
The choice of normalisation constant $A^{\alpha}$ ensures that the kernel $I_{\alpha}$ enjoys the semigroup
property
\begin{equation*}
I_{\alpha+\beta}=I_{\alpha}\ast I_{\beta}~~\mathrm{for~each}~~ \alpha,\beta\in(0,3)~~ \mathrm{such~that}~~ \alpha+\beta<3,
\end{equation*}
see for example \cite[pp. 73-74]{du1970}. Indeed, the classical elliptic equation with Riesz potential has been widely studied, and it also has a rich physical background and mathematical research value, see \cite{Du2019,Moroz2017,Gao2017,Gao2018,Mercuri2016} and the reference therein.

We are going to consider the following Brezis-Nirenberg type problem for the curl-curl equation
\begin{equation}\label{nonlocal case}
\left\{
\begin{aligned}
&\nabla\times(\nabla\times u)+\lambda u=\left(I_{\alpha}\ast|u|^{2^{\ast}_{\alpha}}\right)|u|^{2^{\ast}_{\alpha}-2}u~~&&\mathrm{in}~\Omega,\\
&\nu\times u=0~~&&\mathrm{on}~\partial\Omega,
\end{aligned}
\right.
\end{equation}
where $\Omega\subset\mathbb{R}^{3}$ is a bounded domain, either convex or with $\mathcal{C}^{1,1}$ boundary, $\nu$ is the exterior normal, $\lambda<0$ is a real parameter, $2^{\ast}_{\alpha}=3+\alpha$ with $0<\alpha<3$ is the upper critical exponent in the sense of the following Hardy-Littlewood-Sobolev (HLS for short) inequality, see \cite{Lieb2001}.
\begin{Prop}\label{HLS}
Let $t,r\in(1,\infty)$, $\alpha\in(0,N)$ with $\frac{1}{t}+\frac{N-\alpha}{N}+\frac{1}{r}=2$. For $h\in L^{r}(\mathbb{R}^{N},\mathbb{R}^{N}), g\in L^{t}(\mathbb{R}^{N},\mathbb{R}^{N})$, there exists a sharp constant $C(r,t,N,\alpha)$ independent $g$ and $h$ such that
\begin{equation}\label{HLS1}
\int_{\mathbb{R}^{N}}(I_{\alpha}\ast |h|)|g|dx\leq C(r,t,N,\alpha)||h||_{L^{r}(\mathbb{R}^{N},\mathbb{R}^{N})}||g||_{L^{t}(\mathbb{R}^{N},\mathbb{R}^{N})}.
\end{equation}
If $t=r=\frac{2N}{N+\alpha}$, then there is a equality in (\ref{HLS1}) if and only if $h(x)=Cg(x)$ and
\begin{equation}\label{form}
h(x)=A(\gamma^{2}+|x-a|^{2})^{-\frac{N+\alpha}{2}}
\end{equation}
for some $A\in \mathbb{C}$, $0\neq \gamma\in \mathbb{R}$ and $a\in\mathbb{R}^{N}$.
\end{Prop}
In view of the HLS inequality, the functional corresponds to the nonlocal curl-curl equation
\begin{equation}\label{functional}
J_{\lambda}(u)=\frac{1}{2}\int_{\Omega}|\nabla\times u|^{2}dx+\frac{\lambda}{2}\int_{
\Omega}|u|^{2}dx-\frac{1}{2\cdot2^{\ast}_{\alpha}}\int_{\Omega}|I_{\alpha/2}\ast|u|^{2^{\ast}_{\alpha}}|^{2}dx
\end{equation}
is well defined on the natural space $W^{2^{\ast}}_{0}(\mathrm{curl};\Omega)$. However, due to the appearance of the convolution part, this space is not good enough for us to prove the coercive property of the functional. Therefore, it is necessary to introduce the Coulomb space
\begin{equation*}
Q^{\alpha,2^{\ast}_{\alpha}}(\Omega,\mathbb{R}^{3})=\{u:\int_{\Omega}|I_{\alpha/2}\ast|u|^{2^{\ast}_{\alpha}}|^{2}dx<\infty\},
\end{equation*}
see Definition \ref{coulomb space} below. Then, we may define the Coulomb space involve curl operator as
\begin{equation*}
W^{\alpha,2^{\ast}_{\alpha}}(\mathrm{curl};\Omega)=\{u\in Q^{\alpha,2^{\ast}_{\alpha}}(\Omega,\mathbb{R}^{3}):\nabla\times u\in L^{2}(\Omega,\mathbb{R}^{3})\}
\end{equation*}
is a Banach space (see Lemma \ref{complete1}) if provided with the norm
\begin{equation*}
||u||_{W^{\alpha,2^{\ast}_{\alpha}}(\mathrm{curl};\Omega)}:=(||u||_{Q^{\alpha,2^{\ast}_{\alpha}}}^{2}+|\nabla\times u|^{2}_{2})^{1/2}.
\end{equation*}
We also need the following space
\begin{equation*}
	W^{\alpha,2^{\ast}_{\alpha}}_{0}(\mathrm{curl};\Omega)=\overline{\mathcal{C}_{0}^{\infty}(\Omega,\mathbb{R}^{3})}^{||\cdot||_{W^{\alpha,2^{\ast}_{\alpha}}(\mathrm{curl};\Omega)}}.
\end{equation*}

In this way, we can easily check that the functional (\ref{functional}) is well defined on $W^{\alpha,2^{\ast}_{\alpha}}_{0}(\mathrm{curl};\Omega)$, see Lemma \ref{well define}. In order to obtain the Brezis-Lieb lemma in the dual space, we extend the linear functionals of the Coulomb space to a mix-norm space, see Proposition \ref{dual}. Correspondingly, to establish the Helmholtz decomposition on the work space $W^{\alpha,2^{\ast}_{\alpha}}_{0}(\mathrm{curl};\Omega)$ and $W^{\alpha,2^{\ast}_{\alpha}}_{0}(\mathrm{curl};\mathbb{R}^{3})$, see (\ref{bounded helmholtz}) and Lemma \ref{Helmholtz}, we introduce the following subspace
\begin{equation}\label{V W}
\begin{aligned}
&&\mathcal{V}_{\Omega}:&=\{v\in W^{\alpha,2^{\ast}_{\alpha}}_{0}(\mathrm{curl};\Omega):\int_{\Omega}\langle v,\varphi\rangle dx=0~\mathrm{for~every}~\varphi\in C^{\infty}_{0}(\Omega,\mathbb{R}^{3})~\mathrm{with}~\nabla\times\varphi=0\},\\
&&\mathcal{W}_{\Omega}:&=\{w\in W^{\alpha,2^{\ast}_{\alpha}}_{0}(\mathrm{curl};\Omega):\int_{\Omega}\langle w,\nabla\times\varphi\rangle dx=0~\mathrm{for~all}~\varphi\in C^{\infty}_{0}(\Omega,\mathbb{R}^{3})\}\\
&&~\quad~&=\{w\in W^{\alpha,2^{\ast}_{\alpha}}_{0}(\mathrm{curl};\Omega):\nabla\times w =0~\mathrm{in~the~sense~of~distributions}\}.
\end{aligned}
\end{equation}
Here and below $\langle\cdot,\cdot\rangle$ denote the inner product. Without misunderstanding, we shall write $\mathcal{V}_{\mathbb{R}^{3}}, \mathcal{W}_{\mathbb{R}^{3}}$ if $\Omega=\mathbb{R}^{3}$. Basing on this new decomposition, we can adapt the classical concentration-compactness lemma to suit the new situation. Owing to the concentration-compactness lemma, we obtain the  weak-weak$^{\ast}$ continuity of $J'_{\lambda}(u)$ on the Nehari-Pankov manifold (see (\ref{N lambda}))
 \begin{equation*}\label{N lambda 1}
\mathcal{N}_{\lambda}:=\{u\in W^{\alpha,2^{\ast}_{\alpha}}_{0}(\mathrm{curl};\Omega)\setminus(\widetilde{\mathcal{V}}_{\Omega}\oplus\mathcal{W}_{\Omega}):J'_{\lambda}(u)|_{\mathcal{R}u\oplus\widetilde{\mathcal{V}}_{\Omega}\oplus\mathcal{W}_{\Omega}}=0\},
\end{equation*}
where $\widetilde{\mathcal{V}}_{\Omega}$ is a subspace of $\mathcal{V}_{\Omega}$ on which the quadratic part of $J_{\lambda}$ (see \ref{quad}) is negative semi-definite. Meanwhile, the concentrated compactness lemma implies the $L^{2}(\Omega,\mathbb{R}^{3})$ convergence for the bounded sequence. This allows us to choose the compactly  perturbed functional $J_{cp}=J_{0}=J_{\lambda=0}$ that satisfies the condition $(C1)$ in Lemma \ref{abstract}. By setting another Nehari-Pankov manifold (see (\ref{Ncp}))
 \begin{equation*}\label{Ncp}
\mathcal{N}_{cp} = \{E\in (\mathcal{V}_{\Omega}\oplus\mathcal{W}_{\Omega})\setminus\mathcal{W}_{\Omega}: J'_{cp}(u)|_{\mathbb{R}u\oplus\mathcal{W}_{\Omega}}=0\},
\end{equation*}
and controlling the ground state energy of $J_{\lambda}$ lower than the ground state energy of the perturbed functional $J_{cp}$, i.e.
\begin{equation*}\label{compare c0}
c_{\lambda}=\mathop{\mathrm{inf}}\limits_{\mathcal{N}_{\lambda}}J_{\lambda}\leq \mathop{\mathrm{inf}}\limits_{\mathcal{N}_{cp}}J_{cp}=c_{0},
\end{equation*}
we can obtain the ground state solutions of the curl-curl equation (\ref{nonlocal case}).

It remains to prove the achievement of $c_{0}$. However, as we all know, for the classical elliptic equation (\ref{elliptic}), the sharp constants corresponded to the infimums of the energy level are only attained provided $\Omega=\mathbb{R}^{3}$, even though they are not dependent on the shape of domain. Based on this fact, one can compare the energy levels by the extremal functions. Inspired by the local case in \cite{Mederski2021}, we are motivated to investigate the sharp constant of the Sobolev type inequality involving the curl operator on the entire space $\mathbb{R}^{3}$.
Let $S_{\mathrm{curl,HL}}=S_{\mathrm{curl,HL}}(\mathbb{R}^{3})$ be the largest constant such that the inequality
\begin{equation}\label{inequality1}
\int_{\mathbb{R}^{3}}|\nabla\times u|^{2}dx\geq S_{\mathrm{curl,HL}}\mathop{\mathrm{inf}}\limits_{w\in \mathcal{W}_{\mathbb{R}^{3}}}\left(\int_{\mathbb{R}^{3}}|I_{\alpha/2}\ast|u+w|^{2^{\ast}_{\alpha}}|^{2}dx\right)^{\frac{1}{2^{\ast}_{\alpha}}}.
\end{equation}
holds for any $u\in W^{\alpha,2^{\ast}_{\alpha}}(\mathrm{curl};\mathbb{R}^{3})\setminus\mathcal{W}_{\mathbb{R}^{3}}$. Then the achievement of the sharp constant $S_{\mathrm{curl},HL}$ is related to a certain least energy solution of the limiting problem
\begin{equation}\label{limit equation}
\nabla\times(\nabla\times u)=\left(I_{\alpha}\ast|u|^{2^{\ast}_{\alpha}}\right)|u|^{2^{\ast}_{\alpha}-2}u,~~\mathrm{in}~\mathbb{R}^{3},
\end{equation}
where $u\in W^{\alpha,2^{\ast}_{\alpha}}_{0}(\mathrm{curl};\mathbb{R}^{3})$. By setting the functional
\begin{equation}\label{J(u)}
J(u)=\frac{1}{2}\int_{\mathbb{R}^{3}}|\nabla\times u|^{2}dx-\frac{1}{2\cdot 2^{\ast}_{\alpha}}\int_{\mathbb{R}^{3}}|I_{\alpha/2}\ast|u|^{2^{\ast}_{\alpha}}|^{2}dx,
\end{equation}
and introducing the following Nehari-Pankov manifold (see (\ref{manifold2}))
\begin{equation*}\label{manifold1}
\mathcal{N}:=\left\{u\in W^{\alpha,2^{\ast}_{\alpha}}_{0}(\mathrm{curl};\mathbb{R}^{3})\setminus \mathcal{W}_{\mathbb{R}^{3}}: J'(u)u=0~\mathrm{and}~J'(u)|_{\mathcal{W}_{\mathbb{R}^{3}}}=0\right\},
\end{equation*}
then we have
\begin{thm}\label{attained}
We have the following two conclusions:\\
\noindent $(a)$. $\mathop{\mathrm{inf}}\limits_{\mathcal{N}}  J=\frac{2^{\ast}_{\alpha}-1}{2\cdot2^{\ast}_{\alpha}} S_{\mathrm{curl},HL}^{\frac{2^{\ast}_{\alpha}}{2^{\ast}_{\alpha}-1}}$ and is attained. Moreover, if $u\in \mathcal{N}$ and $J(u)=\mathop{\mathrm{inf}}\limits_{\mathcal{N}}J$, then $u$ is a ground state solution to equation (\ref{limit equation}) and equality holds in (\ref{inequality1}) for this $u$. If $u$ satisfies equality (\ref{inequality1}), then there are unique $t>0$ and $w\in\mathcal{W}_{\mathbb{R}^{3}}$ such that $t(u+w)\in\mathcal{N}$ and $J(t(u+w))=\mathop{\mathrm{inf}}\limits_{\mathcal{N}}J$.

\noindent $(b)$. $S_{\mathrm{curl},HL}> S_{HL}$, where
\begin{equation}\label{SHL}
S_{HL}:=\mathop{\mathrm{inf}}\limits_{u\in\mathcal{D}^{1,2}(\mathbb{R}^{3},\mathbb{R}^{3})\setminus\{0\}}\frac{\int_{\mathbb{R}^{3}}|\nabla u|^{2}dx}{\left(\int_{\mathbb{R}^{3}}|I_{\alpha/2}\ast|u|^{2^{\ast}_{\alpha}}|^{2}dx\right)^{\frac{1}{2^{\ast}_{\alpha}}}}.
\end{equation}
\end{thm}
Note that $S_{HL}$
is the best constant of the combination of the HLS inequality and the Sobolev inequality, see \cite[Lemma 1.2]{Du2019} for example.
It is not clear that whether the sharp constant $S_{\mathrm{curl},HL}$ is independent on shape of the domain $\Omega$ or not. Therefore, we may define another two constants  $S_{\mathrm{curl},HL}(\Omega)$ and $\overline{S}_{\mathrm{curl},HL}(\Omega)$. $S_{\mathrm{curl},HL}(\Omega)$ is the largest possible constant such that the inequality
\begin{equation}\label{inequality2}
\int_{\mathbb{R}^{3}}|\nabla\times u|^{2}dx\geq S_{\mathrm{curl},HL}(\Omega)\mathop{\mathrm{inf}}\limits_{w\in \mathcal{W}_{\mathbb{R}^{3}}}\left(\int_{\mathbb{R}^{3}}|I_{\alpha/2}\ast|u+w|^{2^{\ast}_{\alpha}}|^{2}dx\right)^{\frac{1}{2^{\ast}_{\alpha}}}
\end{equation}
holds for any $u\in W^{\alpha,2^{\ast}_{\alpha}}_{0}(\mathrm{curl};\Omega)\setminus\mathcal{W}_{\mathbb{R}^{3}}$ with a zero extending; $\overline{S}_{\mathrm{curl},HL}(\Omega)$ is another constant such that the inequality
\begin{equation}\label{inequality3}
\int_{\Omega}|\nabla\times u|^{2}dx\geq \overline{S}_{\mathrm{curl},HL}(\Omega)\mathop{\mathrm{inf}}\limits_{w\in \mathcal{W}_{\Omega}}\left(\int_{\Omega}|I_{\alpha/2}\ast|u+w|^{2^{\ast}_{\alpha}}|^{2}dx\right)^{\frac{1}{2^{\ast}_{\alpha}}}.
\end{equation}
holds for any $u\in W^{\alpha,2^{\ast}_{\alpha}}_{0}(\mathrm{curl};\Omega)\setminus\mathcal{W}_{\Omega}$, and $\overline{S}_{\mathrm{curl},HL}(\Omega)$ is
largest with this property. We compare the four constants as follow.
\begin{thm}\label{four constant}
Let $\Omega$ be a bounded domain, either convex or with $\mathcal{C}^{1,1}$ boundary. Then
\begin{equation*}
S_{\mathrm{curl},HL}=S_{\mathrm{curl},HL}(\Omega)\geq \overline{S}_{\mathrm{curl},HL}(\Omega)\geq S_{HL}.
\end{equation*}
\end{thm}
Unfortunately, we don't have any information about the shape of the solutions of (\ref{limit equation}). Hence, the method of taking cut-off functions and comparing the energy levels does not work well any longer. Inspired by the idea in \cite{Mederski2018}, we are going to investigate the energy levels of the ground states. %Before stating the main results, it is useful to introduce some notations.
From \cite{Bartsch2015-1} we know that the spectrum of the curl-curl operator in $W^{2}_{0}(\mathrm{curl};\Omega)$ consists of the eigenvalue $\lambda_{0}=0$ with infinite multiplicity and of a sequence of eigenvalues
\begin{equation*}
0<\lambda_{1}\leq\lambda_{1}\leq\lambda_{2}\leq\cdot\cdot\cdot\leq\lambda_{k}\longrightarrow\infty
\end{equation*}
with  finite multiplicity $m(\lambda_{k})\in\mathbb{N}$.

The main results for the existence are as follow:
\begin{thm}\label{main thm}
Suppose $\Omega$ is a bounded domain, either convex or with $\mathcal{C}^{1,1}$ boundary. Let $\lambda\in(-\lambda_{\nu},-\lambda_{\nu-1}]$ for some $\nu\geq 1$. Then $c_{\lambda}>0$ and the following statements hold:

\noindent $(a) $. If $c_{\lambda}<c_{0}$, then there is ground state solution to (\ref{nonlocal case}) , i.e. $c_{\lambda}$ is attained by a critical point of $J_{\lambda}$. A sufficient condition for this inequality to hold is $\lambda\in\left(-\lambda_{\nu},-\lambda_{\nu}+\bar{S}_{\mathrm{curl},HL}(\Omega)|\mathrm{diam}\Omega|^{-\frac{3\cdot2^{\ast}_{\alpha}-\alpha-3}{2^{\ast}_{\alpha}}}\right)$,
where $|\mathrm{diam}\Omega|=\mathop{max}\limits_{x,y\in\Omega}|x-y|$.

\noindent  $(b)$. There exists $\varepsilon_{\nu}\geq\bar{S}_{\mathrm{curl},HL}(\Omega)|\mathrm{diam}\Omega|^{-\frac{3\cdot2^{\ast}_{\alpha}-\alpha-3}{2^{\ast}_{\alpha}}}$ such that $c_{\lambda}$ is not attained for $\lambda\in(-\lambda_{\nu}+\varepsilon_{\nu},-\lambda_{\nu-1}]$, and $c_{\lambda}=c_{0}$ for $\lambda\in(-\lambda_{\nu}+\varepsilon_{\nu},-\lambda_{\nu-1}]$. We do not exclude that $\varepsilon>\lambda_{\nu}-\lambda_{\nu-1}$, so these intervals may be empty.

\noindent $(c)$. $c_{\lambda}\longrightarrow 0$ as $\lambda\longrightarrow-\lambda^{-}_{\nu}$, and the function
\begin{equation*}
(-\lambda_{\nu},-\lambda_{\nu}+\varepsilon_{\nu}]\cap (-\lambda_{\nu},-\lambda_{\nu-1}]\ni\lambda\mapsto c_{\lambda}\in(0,\infty)
\end{equation*}
is continuous and strictly increasing.

\noindent $(d)$. There exist at least $\sharp\{k:-\lambda_{k}<\lambda<-\lambda_{\nu}+\frac{2^{\ast}_{\alpha}-1}{2\cdot2^{\ast}_{\alpha}}\bar{S}_{\mathrm{curl},HL}(\Omega)|\mathrm{diam}\Omega|^{-\frac{3\cdot2^{\ast}_{\alpha}-\alpha-3}{2^{\ast}_{\alpha}}}\} $ pairs of solutions $\pm u$ to (\ref{nonlocal case}).
\end{thm}

The paper is organized as follow. In Section 2 we introduce some work spaces on bounded domains and entire space $\mathbb{R}^{3}$, and we adapt the concentration compactness lemma for the curl-curl problem with nonlocal nonlinearities. And an abstract critical point theorem is also recalled in this part for the readers' convenience. In Section 3, we show that the sharp constant $S_{\mathrm{curl},HL}$ is attained provided $\Omega=\mathbb{R}^{3}$, and we compare the four constants as we introduced. In the last Section, we are devoted to the proof of Theorem \ref{main thm}.

\section{\bf  Preliminaries and Variational Setting}
\subsection{Preliminaries}\label{setting}

Throughout this paper we assume that $\Omega\subset\mathbb{R}^{3}$ is a bounded domain, either convex or with $\mathcal{C}^{1,1}$ boundary. In some cases $\Omega$ is only required to be a Lipschitz domain, see \cite{Mederski2021} for more details. We shall look for solutions to problem (\ref{nonlocal case}) and (\ref{limit equation}) in $W^{\alpha,2^{\ast}_{\alpha}}_{0}(\mathrm{curl};\Omega)$ and $W^{\alpha,2^{\ast}_{\alpha}}_{0}(\mathrm{curl};\mathbb{R}^{3})$ respectively. Now we are ready to introdcue the definitions of the working spaces.

\subsubsection{Coulomb Space Involving Curl Operator}
\begin{Def}\cite[Definition 1]{Mercuri2016}\label{coulomb space}
Let $N\in\mathbb{N}$, $\alpha\in(0,N)$ and $p\geq1$. We define the Coulomb space $Q^{\alpha,p}(\mathbb{R}^{N},\mathbb{R}^{N})$ as the vector space of measurable functions $u: \mathbb{R}^{N}\longrightarrow\mathbb{R}^{N}$ such that
\begin{equation*}
||u||_{Q^{\alpha,p}(\mathbb{R}^{N},\mathbb{R}^{N})}=\left(\int_{\mathbb{R}^{N}}|I_{\alpha/2}\ast|u|^{p}|^{2}dx\right)^{\frac{1}{2p}}<+\infty.
\end{equation*}
\end{Def}
It is not difficult to see that $||\cdot||_{Q^{\alpha,p}(\mathbb{R}^{N},\mathbb{R}^{N})}$ defines a norm, see Proposition 2.1 in \cite{Mercuri2016}, and the Coulomb space is complete with respect to this norm. By the same way, we also define $Q^{\alpha,p}(\Omega,\mathbb{R}^{N})$ as the Coulomb space on the bounded domain. For the the dual space of $Q^{\alpha,p}(\mathbb{R}^{N},\mathbb{R}^{N})$, it can be characterized by the following proposition, and it is also adopted in $Q^{\alpha,p}(\Omega,\mathbb{R}^{N})$.
\begin{Prop}\cite[Proposition 2.11]{Mercuri2016}\label{dual}
Let $T$ be a distribution, then $T\in \left(Q^{\alpha,p}(\mathbb{R}^{N},\mathbb{R}^{N})\right)'$ if and only if there exists $G(x,y)\in L^{\frac{2p}{2p-1}}(\mathbb{R}^{N},L^{\frac{p}{p-1}}(\mathbb{R}^{N}))$ such that for every $\varphi\in\mathcal{C}_{0}^{\infty}(\mathbb{R}^{N},\mathbb{R}^{N})$,
\begin{equation*}
\langle T,\varphi\rangle=\int_{\mathbb{R}^{N}}\left(\int_{\mathbb{R}^{N}}G(x,y)I_{\alpha/2}(x-y)^{\frac{1}{p}}dy\right)\varphi(x)dx.
\end{equation*}
\end{Prop}
\begin{proof}
By the definition of the Coulomb space $Q^{\alpha,p}(\mathbb{R}^{N},\mathbb{R}^{N})$, one can observe that the map
\begin{equation*}
\mathcal{L}:Q^{\alpha,p}(\mathbb{R}^{N},\mathbb{R}^{N})\longrightarrow L^{2p}(\mathbb{R}^{N},L^{p}(\mathbb{R}^{N}))
\end{equation*}
defined by
\begin{equation*}
\mathcal{L}u(x,y)=(I_{\alpha/2}(x-y))^{\frac{1}{p}}u(y)
\end{equation*}
is a linear isometry from $Q^{\alpha,p}(\mathbb{R}^{N},\mathbb{R}^{N})$ into $L^{2p}(\mathbb{R}^{N},L^{p}(\mathbb{R}^{N}))$. Then any linear functional on $Q^{\alpha,p}(\mathbb{R}^{N})$ can be extended to a linear functional on
$L^{2p}(\mathbb{R}^{N},L^{p}(\mathbb{R}^{N}))$. Namely, there exists $G(x,y)\in L^{\frac{2p}{2p-1}}(\mathbb{R}^{N},L^{\frac{p}{p-1}}(\mathbb{R}^{N}))$ such that
\begin{equation*}
\langle T,\varphi\rangle=\langle G(x,y),\mathcal{L}\varphi\rangle.
\end{equation*}
\end{proof}

For the Coulomb space involving curl operator, we have the following definition.
\begin{Def}\label{curl}
Let $N=3$, $\alpha\in(0,3)$ and $p\geq1$. We define the Coulomb space involving curl operator $W^{\alpha,p}(\mathrm{curl};\mathbb{R}^{3})$ as the vector space of functions $u\in Q^{\alpha,p}(\mathbb{R}^{3},\mathbb{R}^{3})$ such that u is weakly differentiable
in $\mathbb{R}^{3}$, $\nabla\times u\in L^{2}(\mathbb{R}^{3},\mathbb{R}^{3})$ and
\begin{equation*}
||u||_{W^{\alpha,p}(\mathrm{curl};\mathbb{R}^{3})}=\left(\int_{\mathbb{R}^{3}}|\nabla\times u|^{2}dx+\left(\int_{\mathbb{R}^{3}}|I_{\alpha/2}\ast|u|^{p}|^{2}dx\right)^{\frac{1}{p}}\right)^{\frac{1}{2}}
<+\infty.
\end{equation*}
\end{Def}
The function $||\cdot||_{W^{\alpha,p}(\mathrm{curl};\mathbb{R}^{3})}$ defines a norm in view of the Proposition 2.1 in \cite{Mercuri2016}. By the same way, we also define $W^{\alpha,p}(\mathrm{curl};\Omega)$ as the Coulomb space involve curl operator on the bounded spaces, namely,
\begin{equation}\label{work space 2}
W^{\alpha,p}(\mathrm{curl};\Omega):=\{u\in Q^{\alpha,p}(\Omega,\mathbb{R}^{3}):\nabla\times u\in L^{2}(\Omega,\mathbb{R}^{3})\}.
\end{equation}
We are going to prove that $W^{\alpha,p}(\mathrm{curl};\mathbb{R}^{3})$ and $W^{\alpha,p}(\mathrm{curl};\Omega)$ are Banach spaces. The proof of completeness follows by the same arguments as in the proof of Theorem 4.3 in \cite{Kirsch2015} and Proposition 2.2 in \cite{Mercuri2016}. The first ingredient is the following Fatou property for locally converging sequences.
\begin{lem}\label{Fatou}
Let $N=3$, $\alpha\in(0,3)$ and $p\geq1$. If $(u_{n})_{n\in\mathbb{N}}$ is a bounded sequence in $W^{\alpha,p}(\mathrm{curl};\mathbb{R}^{3})$ that converges to a function $u:\mathbb{R}^{3}\longrightarrow \mathbb{R}^{3}$ in $L^{1}_{loc}(\mathbb{R}^{3},\mathbb{R}^{3})$, then $u\in W^{\alpha,p}(\mathrm{curl};\mathbb{R}^{3})$,
\begin{equation}\label{Fatou1}
\int_{\mathbb{R}^{3}}|I_{\alpha/2}\ast| u|^{p}|^{2}dx\leq\mathop{\mathrm{lim~inf}}\limits_{n\longrightarrow\infty}\int_{\mathbb{R}^{3}}|I_{\alpha/2}\ast| u_{n}|^{p}|^{2}dx,
\end{equation}
and
\begin{equation}\label{Fatou2}
\int_{\mathbb{R}^{3}}|\nabla\times u|^{2}dx\leq\mathop{\mathrm{lim~inf}}\limits_{n\longrightarrow\infty}\int_{\mathbb{R}^{3}}|\nabla\times u_{n}|^{2}dx.
\end{equation}
\end{lem}
\begin{proof}
Since $(u_{n})_{n\in\mathbb{N}}\longrightarrow u$ is bounded in $W^{\alpha,p}(\mathrm{curl};\mathbb{R}^{3})$, we have
\begin{equation}\label{Fatou3}
\int_{\mathbb{R}^{3}}|I_{\alpha/2}\ast|u_{n}|^{p}|^{2}dx\leq \infty,
\end{equation}
then by the Fatou lemma, we have
\begin{equation}\label{Fatou4}
\int_{\mathbb{R}^{3}}\mathop{\mathrm{lim~inf}}\limits_{n\longrightarrow\infty}|I_{\alpha/2}\ast| u_{n}|^{p}|^{2}dx\leq\mathop{\mathrm{lim~inf}}\limits_{n\longrightarrow\infty}\int_{\mathbb{R}^{3}}|I_{\alpha/2}\ast| u_{n}|^{p}|^{2}dx.
\end{equation}
By the Fatou lemma again, we have
\begin{equation}\label{Fatou5}
 I_{\alpha}\ast(\mathop{\mathrm{lim~inf}}\limits_{n\longrightarrow\infty} |u_{n}|^{p}) \leq\mathop{\mathrm{lim~inf}}\limits_{n\longrightarrow\infty} I_{\alpha}\ast(|u_{n}|^{p}).
\end{equation}
Since $(u_{n})_{n\in\mathbb{N}}\longrightarrow u$ in $L^{1}_{loc}(\mathbb{R}^{3},\mathbb{R}^{3})$, for almost every $x\in\mathbb{R}^{3}$, we have
\begin{equation}\label{Fatou6}
 I_{\alpha}\ast(\mathop{\mathrm{lim~inf}}\limits_{n\longrightarrow\infty} |u_{n}|^{p})(x)
  \longrightarrow I_{\alpha}\ast(|u|^{p})(x).
\end{equation}
Then (\ref{Fatou1}) follows (\ref{Fatou4}),(\ref{Fatou5}) and (\ref{Fatou6}).

We are going to prove (\ref{Fatou2}). Define $f$ on $\mathcal{D}(\mathbb{R}^{3},\mathbb{R}^{3})$ by
\begin{equation}
\langle f,v\rangle=\int_{\mathbb{R}^{3}}u\cdot (\nabla\times v)~dx,
\end{equation}
since $u_{n}\longrightarrow u$ in $L^{1}_{loc}(\mathbb{R}^{3},\mathbb{R}^{3})$, we have
\begin{equation}
\begin{aligned}
|\langle f,v\rangle|&=|\int_{\mathbb{R}^{3}}u\cdot (\nabla\times v)|=\mathop{\mathrm{lim}}\limits_{n\longrightarrow\infty}|\int_{\mathbb{R}^{3}}u_{n}\cdot (\nabla\times v)~dx|\\
&=\mathop{\mathrm{lim}}\limits_{n\longrightarrow\infty}|\int_{\mathbb{R}^{3}}(\nabla\times u_{n})\cdot v~dx|\leq \mathop{\mathrm{lim~inf}}\limits_{n\longrightarrow\infty}||\nabla\times u_{n}||_{2}\left(\int_{\mathbb{R}^{3}}|v|^{2}dx\right)^{\frac{1}{2}},
\end{aligned}
\end{equation}
where we use the Cauchy-Schwarz inequality. Since $\mathcal{D}(\mathbb{R}^{3},\mathbb{R}^{3})$ is dense in $L^{2}(\mathbb{R}^{3},\mathbb{R}^{3})$, by the the Hahn-Banach theorem, the distribution $f$ can be continuously extend to a linear functional on $L^{2}(\mathbb{R}^{3},\mathbb{R}^{3})$. Therefore, by the Riesz representation theorem, there exists $F\in L^{2}(\mathbb{R}^{3},\mathbb{R}^{3})$ such that for every
$v\in \mathcal{D}(\mathbb{R}^{3},\mathbb{R}^{3})$
\begin{equation}
\int_{\mathbb{R}^{3}}F\cdot v~dx=\langle f,v\rangle=\int_{\mathbb{R}^{3}}u\cdot (\nabla\times v)~dx.
\end{equation}
Setting $\nabla\times u$ as the curl of $u$ in the following distribute sense
\begin{equation}
\int_{\mathbb{R}^{3}}u\cdot(\nabla\times v)~dx=\int_{\mathbb{R}^{3}}(\nabla\times u)\cdot v~dx,
\end{equation}
we can see $F=\nabla\times u\in L^{2}(\mathbb{R}^{3},\mathbb{R}^{3})$ in the weak sense. Choosing $v=\nabla\times u$ we find that
\begin{equation}
\begin{aligned}
\int_{\mathbb{R}^{3}}|\nabla\times u|^{2}~dx&\leq \mathop{\mathrm{lim~inf}}\limits_{n\longrightarrow\infty}||\nabla\times u_{n}||_{2}\left(\int_{\mathbb{R}^{3}}|v|^{2}\right)^{\frac{1}{2}}\\
&\leq \mathop{\mathrm{lim~inf}}\limits_{n\longrightarrow\infty}||\nabla\times u_{n}||_{2}\left(\int_{\mathbb{R}^{3}}|\nabla\times u|^{2}\right)^{\frac{1}{2}}.
\end{aligned}
\end{equation}
Therefore we have
\begin{equation}\label{fatou10}
\int_{\mathbb{R}^{3}}|\nabla\times u|^{2}dx\leq\mathop{\mathrm{lim~inf}}\limits_{n\longrightarrow\infty}\int_{\mathbb{R}^{3}}|\nabla\times u_{n}|^{2}dx.
\end{equation}
\end{proof}

\begin{lem}\label{complete1}
Let N=3, $\alpha\in(0,3)$ and $p\geq1$. The normed spaces $W^{\alpha,p}(\mathrm{curl};\mathbb{R}^{3})$ and $W^{\alpha,p}(\mathrm{curl};\Omega)$ are complete.
\end{lem}
\begin{proof}
Let $(u_{n})_{n\in\mathbb{N}}$ be a Cauchy sequence in $W^{\alpha,p}(\mathrm{curl};\mathbb{R}^{3})$. By the local estimate of the Coulomb energy, $(u_{n})_{n\in\mathbb{N}}$ is also a Cauchy sequence in $L_{loc}^{p}(\mathbb{R}^{3},\mathbb{R}^{3})$. Hence there exists $u\in L_{loc}^{p}(\mathbb{R}^{3},\mathbb{R}^{3})$ such that $(u_{n})_{n\in\mathbb{N}}\longrightarrow u$ in $L_{loc}^{p}(\mathbb{R}^{3},\mathbb{R}^{3})$. In light of Lemma \ref{Fatou}, we conclude that $u\in W^{\alpha,p}(\mathrm{curl};\mathbb{R}^{3})$. Moreover, for every $n\in\mathbb{N}$ the sequence $(u_{n}-u_{m})_{m\in\mathbb{N}}$ converges to $(u_{n}-u)$ in $L_{loc}^{p}(\mathbb{R}^{3},\mathbb{R}^{3})$. Hence, by Lemma \ref{Fatou} again, we have
\begin{equation*}
\begin{aligned}
&\mathop{\mathrm{lim~sup}}\limits_{n\longrightarrow\infty}
\left(\int_{\mathbb{R}^{3}}|\nabla\times u_{n}-\nabla\times u|^{2}dx+\int_{\mathbb{R}^{3}}|I_{\alpha/2}\ast|u_{n}-u|^{p}|^{2}dx\right)\\
&\leq\mathop{\mathrm{lim~sup}}\limits_{n\longrightarrow\infty}\mathop{\mathrm{lim~sup}}\limits_{m\longrightarrow\infty}
\left(\int_{\mathbb{R}^{3}}|\nabla\times u_{n}-\nabla\times u_{m}|^{2}dx+\int_{\mathbb{R}^{3}}|I_{\alpha/2}\ast|u_{n}-u_{m}|^{p}|^{2}dx\right)\\
&\leq\mathop{\mathrm{lim~sup}}\limits_{m,n\longrightarrow\infty}
\left(\int_{\mathbb{R}^{3}}|\nabla\times u_{n}-\nabla\times u_{m}|^{2}dx+\int_{\mathbb{R}^{3}}|I_{\alpha/2}\ast|u_{n}-u_{m}|^{p}|^{2}dx\right)\\
&\leq0.
\end{aligned}
\end{equation*}
This implies $W^{\alpha,p}(\mathrm{curl};\mathbb{R}^{3})$ is complete. The completeness of $W^{\alpha,p}(\mathrm{curl};\Omega)$ can be proved in the same way.
\end{proof}

We also define
\begin{equation*}\label{work space 1}
W^{\alpha,p}_{0}(\mathrm{curl};\Omega)=\mathrm{closure}~\mathrm{of}~\mathcal{C}^{\infty}_{0}(\Omega;\mathbb{R}^{3})~\mathrm{in}~W^{\alpha,p}(\mathrm{curl};\Omega).
\end{equation*}
If $p$ lies in some suitable range, then the two Coulomb spaces are the same for the case $\Omega=\mathbb{R}^{3}$.
\begin{lem}
Let $\alpha\in(0,3)$, then $W^{\alpha,p}(\mathrm{curl};\mathbb{R}^{3})=W^{\alpha,p}_{0}(\mathrm{curl};\mathbb{R}^{3})$ as $\frac{3+\alpha}{3}\leq p\leq 3+\alpha$.
%$\frac{6-\mu}{3}\leq p\leq 3+\mu$ for $0\leq\mu\leq\frac{3}{2}$ and as $\frac{6-\mu}{3}\leq p\leq %6-\mu$ for $\frac{3}{2}\leq \mu\leq 3$.
\end{lem}
\begin{proof}
Let $\eta_{R}\in \mathcal{C}^{\infty}_{0}(\mathbb{R}^{3})$ be such that $|\nabla\eta_{R}|\leq\frac{2}{R}$ for $R\leq|x|\leq 2R$, $\eta_{R}=1$ for $|x|\leq R$ and $\eta_{R}=0$ for $|x|\geq 2R$. Then for $u=(u_{1},u_{2},u_{3})\in W^{\alpha,p}(\rm{curl},\mathbb{R}^{3})$, we have $\eta_{R}u\longrightarrow u$ in $Q^{\alpha,p}(\mathbb{R}^{3},\mathbb{R}^{3})$ as $R\longrightarrow\infty$. Note that
\begin{equation}\label{partial}
\nabla\times(\eta_{R}u_{i})
%=\partial_{i}(\eta_{R}u_{j})-\partial_{j}(\eta_{R}u_{i})
=(\partial_{i}\eta_{R})u_{j}-(\partial_{j}\eta_{R})u_{i}+\eta_{R}(\partial_{i}u_{j}-\partial_{j}u_{i}),~~i\neq j.
\end{equation}
If $p=2$, we have $(\partial_{i}\eta_{R})u_{j}\longrightarrow 0$ in $L^{2}(\mathbb{R}^{3},\mathbb{R}^{3})$ as $\alpha\in(0,3)$. Indeed, we have
\begin{equation*}
\begin{aligned}
\int_{\mathbb{R}^{3}}(\partial\eta_{R})^{2}u_{j}^{2}dx&=\int_{R\leq|x|\leq 2R}(\partial\eta_{R})^{2}u_{j}^{2}dx
\leq\left(\frac{2}{R}\right)^{2}\int_{R\leq|x|\leq 2R}u_{j}^{2}dx\leq\left(\frac{2}{R}\right)^{2}\int_{|x|\leq 2R}u_{j}^{2}dx\\
&\leq CR^{-2}(2R)^{\frac{3-\alpha}{2}}\left(\int_{|x|\leq 2R}|I_{\alpha/2}\ast|u|^{2}|^{2}\right)^{\frac{1}{2}}.
\end{aligned}
\end{equation*}
%where the last inequality is a basic property of Coulomb space, see \cite[Proposition %2.3]{Mercuri2016}.
If $p\neq2$, let  $q$ be such that $\frac{1}{p}+\frac{1}{q}=\frac{1}{2}$, then applying the H\"{o}lder inequality we have
\begin{equation*}
\begin{aligned}
\int_{\mathbb{R}^{3}}(\partial\eta_{R})^{2}u_{j}^{2}dx&\leq\left(\int_{R\leq|x|\leq 2R}|\partial_{i}\eta_{R}|^{q}dx\right)^{\frac{2}{q}}\left(\int_{R\leq|x|\leq 2R}|u_{j}|^{p}dx\right)^{\frac{2}{p}}\\
&\leq C_{1}(R^{3-q})^{\frac{2}{q}}(2R^{\frac{3-\alpha}{2}})^{\frac{2}{p}}\left(\int_{|x|\leq 2R}|I_{\frac{\alpha}{2}}\ast|u|^{p}|^{2}dx\right)^{\frac{1}{p}}\\
&\leq C_{2}R^{\frac{p-(3+\alpha)}{p}}\left(\int_{|x|\leq 2R}|I_{\alpha/2}\ast|u|^{p}|^{2}dx\right)^{\frac{1}{p}}.
\end{aligned}
\end{equation*}
Then, for $p\leq3+\alpha$, we have $(\partial_{i}\eta_{R})u_{j}\longrightarrow 0$ in $L^{2}(\mathbb{R}^{3},\mathbb{R}^{3})$ as $R\longrightarrow \infty$. As $\partial_{i}u_{j}-\partial_{j}u_{i}\in L^{2}(\mathbb{R}^{3})$, it follows that the left-hand side in (\ref{partial}) tends to $\partial_{i}u_{j}-\partial_{j}u_{i}$ in $L^{2}(\mathbb{R}^{3})$ as $R\longrightarrow\infty$.  Hence $\eta_{R}u\longrightarrow u$ in $W^{\alpha,p}(\rm{curl};\mathbb{R}^{3})$ and functions of compact support are dense in $W^{\alpha,p}(\rm{curl};\mathbb{R}^{3})$. The rest of the proof is similar to the \cite[Lemma 2.1]{Mederski2021}.
\end{proof}

\begin{lem}\label{well define}
(i)~$J_{\lambda}(u)$ and $J(u)$ are well defined on $W^{\alpha,2^{\ast}_{\alpha}}_{0}(\mathrm{curl};\Omega)$ and $W^{\alpha,2^{\ast}_{\alpha}}_{0}(\mathrm{curl};\mathbb{R}^{3})$ respectively.

(ii)~Let \begin{equation*}\label{G}
G(u)(x,y)=\left(\frac{1}{|x-y|^{\frac{2^{\ast}_{\alpha}-1}{2^{\ast}_{\alpha}}N-\frac{2\cdot2^{\ast}_{\alpha}-1}{2\cdot2^{\ast}_{\alpha}}\alpha}}\right)|u(y)|^{2^{\ast}_{\alpha}}|u(x)|^{2^{\ast}_{\alpha}-1}.
\end{equation*}
Then, $G(u)(x,y)\in L^{\frac{2\cdot2^{\ast}_{\alpha}}{2\cdot2^{\ast}_{\alpha}-1}}(\mathbb{R}^{3};L^{\frac{2^{\ast}_{\alpha}}{2^{\ast}_{\alpha}-1}}(\mathbb{R}^{3}))$.

(iii)~$J_{\lambda}(u)$ and $J(u)$ are of class $\mathcal{C}^{1}$.
\end{lem}
\begin{proof}
(i)~From the definition of $W^{\alpha,2^{\ast}_{\alpha}}_{0}(\mathrm{curl};\Omega)$ and $W^{\alpha,2^{\ast}_{\alpha}}_{0}(\mathrm{curl};\mathbb{R}^{3})$ , we know that the functionals $J_{\lambda}(u)$ and $J(u)$ are well defined.

(ii) Set
\begin{equation*}
I(u)=\frac{1}{2\cdot2^{\ast}_{\alpha}}\int_{\mathbb{R}^{3}}|I_{\alpha/2}\ast|u|^{2^{\ast}_{\alpha}}|^{2}dx.
\end{equation*}
We claim that $I'(u)\in \left(Q^{\alpha,2^{\ast}_{\alpha}}(\mathbb{R}^{3},\mathbb{R}^{3})\right)'$.
Indeed, for any $\varphi\in Q^{\alpha,2^{\ast}_{\alpha}}(\mathbb{R}^{3},\mathbb{R}^{3})$, we have
\begin{equation}
\begin{aligned}
\langle I'(u),\varphi\rangle&=\int_{\mathbb{R}^{3}}(I_{\alpha} \ast|u|^{2^{\ast}_{\alpha}})|u|^{2^{\ast}_{\alpha}-2}u\cdot\varphi dx
%&=\int_{\mathbb{R}^{3}}(I_{\alpha} %\ast|u|^{2^{\ast}_{\alpha}})^{\frac{2^{\ast}_{\alpha}-1}{2^{\ast}_{\alpha}}}
%(I_{\alpha}\ast|u|^{2^{\ast}_{\alpha}})^{\frac{1}{2^{\ast}_{\alpha}}}|u|^{2^{\ast}_{\alpha}-2}u\cdot\varphi %dx\\
\leq\left(\int_{\mathbb{R}^{3}}(I_{\alpha}\ast|u|^{2^{\ast}_{\alpha}})|u|^{2^{\ast}_{\alpha}}\right)^{\frac{2^{\ast}_{\alpha}-1}{2^{\ast}_{\alpha}}} \cdot\left(\int_{\mathbb{R}^{3}}(I_{\alpha}\ast|u|^{2^{\ast}_{\alpha}})|\varphi|^{2^{\ast}_{\alpha}}\right)^{\frac{1}{2^{\ast}_{\alpha}}}\\ &=||u||^{\frac{2^{\ast}_{\alpha}-1}{2}}_{Q^{\alpha,2^{\ast}_{\alpha}}(\mathbb{R}^{3},\mathbb{R}^{3})}
\cdot\left(\int_{\mathbb{R}^{3}}(I_{\alpha/2}\ast|u|^{2^{\ast}_{\alpha}})(I_{\alpha/2}\ast|\varphi|^{2^{\ast}_{\alpha}})\right)^{\frac{1}{2^{\ast}_{\alpha}}}\\
&\leq||u||^{\frac{2^{\ast}_{\alpha}-1}{2}}_{Q^{\alpha,2^{\ast}_{\alpha}}(\mathbb{R}^{3},\mathbb{R}^{3})}
\cdot\left(\int_{\mathbb{R}^{3}}|I_{\alpha/2}\ast|u|^{2^{\ast}_{\alpha}}|^{2}dx\right)^{\frac{1}{2\cdot2^{\ast}_{\alpha}}}
\cdot\left(\int_{\mathbb{R}^{3}}|I_{\alpha/2}\ast|\varphi|^{2^{\ast}_{\alpha}}|^{2}dx\right)^{\frac{1}{2\cdot2^{\ast}_{\alpha}}}\\
&%=||u||^{\frac{2^{\ast}_{\alpha}-1}{2}}_{Q^{\alpha,2^{\ast}_{\alpha}}(\mathbb{R}^{3},\mathbb{R}^{3})}\cdot
%||u||_{Q^{\alpha,2^{\ast}_{\alpha}}(\mathbb{R}^{3},\mathbb{R}^{3})}\cdot
%||\varphi||_{Q^{\alpha,2^{\ast}_{\alpha}}(\mathbb{R}^{3},\mathbb{R}^{3})}
=||u||^{\frac{2^{\ast}_{\alpha}+1}{2}}_{Q^{\alpha,2^{\ast}_{\alpha}}(\mathbb{R}^{3},\mathbb{R}^{3})}\cdot
||\varphi||_{Q^{\alpha,2^{\ast}_{\alpha}}(\mathbb{R}^{3},\mathbb{R}^{3})}.
\end{aligned}
\end{equation}
%there exist $I'(u): Q^{\alpha,2^{\ast}_{\alpha}}(\mathbb{R}^{3},\mathbb{R}^{3})\longrightarrow(Q^{\alpha,2^{\ast}_{\alpha}}(\mathbb{R}^{3},\mathbb{R}^{3}))'$. Indeed,
Then, by the definition of the functional space on Coulomb space, we have $I'(u)\in \left(Q^{\alpha,2^{\ast}_{\alpha}}(\mathbb{R}^{3},\mathbb{R}^{3})\right)'$.

On the other hand, $G(u)(x,y)$ obviously satisfies that
\begin{equation}
\langle I'(u),\varphi\rangle=\int_{\mathbb{R}^{3}}\left(\int_{\mathbb{R}^{3}}G(u)(x,y)(I_{\alpha/2}(x-y))^{\frac{1}{2^{\ast}_{\alpha}}}dy\right)\varphi(x)dx,
\end{equation}
Therefore, by Proposition \ref{dual}, we have $G(u)(x,y)\in L^{\frac{2\cdot2^{\ast}_{\alpha}}{2\cdot2^{\ast}_{\alpha}-1}}(\mathbb{R}^{3};L^{\frac{2^{\ast}_{\alpha}}{2^{\ast}_{\alpha}-1}}(\mathbb{R}^{3}))$.

(iii)~We are going to show that $I'(u)$ is continuous. For any sequences $u_{n}$, $u$ and $\varphi\in Q^{\alpha,2^{\ast}_{\alpha}}(\mathbb{R}^{3},\mathbb{R}^{3})$, we have
\begin{equation}
\begin{aligned}
&\langle(I'(u_{n})-I'(u)),\varphi\rangle\\
%=\langle (I'(u_{n})-I'(u)),\varphi\rangle
&=\int_{\mathbb{R}^{3}}\left(I_{\alpha}\ast|u_{n}|^{2^{\ast}_{\alpha}}\right)|u_{n}|^{2^{\ast}_{\alpha}-2}u_{n}\cdot\varphi dx-
\int_{\mathbb{R}^{3}}\left(I_{\alpha}\ast|u|^{2^{\ast}_{\alpha}}\right)|u|^{2^{\ast}_{\alpha}-2}u\cdot\varphi dx\\
%&=\int_{\mathbb{R}^{3}}\left(I_{\alpha}\ast|u_{n}|^{2^{\ast}_{\alpha}}\right)|u_{n}|^{2^{\ast}_{\alpha}-2}u_{n}\cdot\varphi %dx-
%\int_{\mathbb{R}^{3}}\left(I_{\alpha}\ast|u|^{2^{\ast}_{\alpha}}\right)|u_{n}|^{2^{\ast}_{\alpha}-2}u_{n}\cdot\varphi %dx\\
%&~~~~+\int_{\mathbb{R}^{3}}\left(I_{\alpha}\ast|u|^{2^{\ast}_{\alpha}}\right)|u_{n}|^{2^{\ast}_{\alpha}-2}u_{n}\cdot\varphi %dx-
%\int_{\mathbb{R}^{3}}\left(I_{\alpha}\ast|u|^{2^{\ast}_{\alpha}}\right)|u|^{2^{\ast}_{\alpha}-2}u\cdot\varphi %dx\\
%&=\int_{\mathbb{R}^{3}}\left(I_{\alpha}\ast(|u_{n}|^{2^{\ast}_{\alpha}}-|u|^{2^{\ast}_{\alpha}})\right)|u_{n}|^{2^{\ast}_{\alpha}-2}u_{n}\cdot\varphi %dx+
%\int_{\mathbb{R}^{3}}\left(I_{\alpha}\ast|u|^{2^{\ast}_{\alpha}}\right)(|u_{n}|^{2^{\ast}_{\alpha}-2}u_{n}-|u|^{2^{\ast}_{\alpha}-2}u)\cdot\varphi %dx\\
&=\int_{\mathbb{R}^{3}}\left(I_{\alpha}\ast(|u_{n}|^{2^{\ast}_{\alpha}}-|u|^{2^{\ast}_{\alpha}})\right)^{\frac{1}{2^{\ast}_{\alpha}}}\left(I_{\alpha}\ast(|u_{n}|^{2^{\ast}_{\alpha}}-|u|^{2^{\ast}_{\alpha}})\right)^{\frac{2^{\ast}_{\alpha}-1}{2^{\ast}_{\alpha}}}|u_{n}|^{2^{\ast}_{\alpha}-2}u_{n}\cdot\varphi dx\\
&~~~~+\int_{\mathbb{R}^{3}}\left(I_{\alpha}\ast|u|^{2^{\ast}_{\alpha}}\right)^{\frac{1}{2^{\ast}_{\alpha}}}\left(I_{\alpha}\ast|u|^{2^{\ast}_{\alpha}}\right)^{\frac{2^{\ast}_{\alpha}-1}{2^{\ast}_{\alpha}}}(|u_{n}|^{2^{\ast}_{\alpha}-2}u_{n}-|u|^{2^{\ast}_{\alpha}-2}u)\cdot\varphi dx\\
&\leq\left(\int_{\mathbb{R}^{3}}\left(I_{\alpha}\ast(|u_{n}|^{2^{\ast}_{\alpha}}-|u|^{2^{\ast}_{\alpha}}) \right)|\varphi|^{2^{\ast}_{\alpha}}dx\right)^{\frac{1}{2^{\ast}_{\alpha}}}\cdot \left(\int_{\mathbb{R}^{3}}\left(I_{\alpha}\ast(|u_{n}|^{2^{\ast}_{\alpha}}-|u|^{2^{\ast}_{\alpha}}) \right)|u_{n}|^{2^{\ast}_{\alpha}}dx\right)^{\frac{2^{\ast}_{\alpha}-1}{2^{\ast}_{\alpha}}}\\
&~~~~+\left(\int_{\mathbb{R}^{3}}\left(I_{\alpha}\ast|u|^{2^{\ast}_{\alpha}} \right)|\varphi|^{2^{\ast}_{\alpha}}dx\right)^{\frac{1}{2^{\ast}_{\alpha}}}\cdot \left(\int_{\mathbb{R}^{3}}\left(I_{\alpha}\ast|u|^{2^{\ast}_{\alpha}} \right)(|u_{n}|^{2^{\ast}_{\alpha}-1}-|u|^{2^{\ast}_{\alpha}-1})^{\frac{2^{\ast}_{\alpha}}{2^{\ast}_{\alpha}-1}}dx\right)^{\frac{2^{\ast}_{\alpha}-1}{2^{\ast}_{\alpha}}}=A_{1}\cdot A_{2}+A_{3}\cdot A_{4}.\\
\end{aligned}
\end{equation}
By the semi-group property and H\"{o}lder inequality, we have
\begin{equation}
\begin{aligned}
A_{1}&=\left(\int_{\mathbb{R}^{3}}\left(I_{\alpha/2}\ast(|u_{n}|^{2^{\ast}_{\alpha}}-|u|^{2^{\ast}_{\alpha}})\right)\left(I_{\alpha/2}\ast|\varphi|^{2^{\ast}_{\alpha}}\right)   dx\right)^{\frac{1}{2^{\ast}_{\alpha}}}\\
&\leq\left(\int_{\mathbb{R}^{3}}|I_{\alpha/2}\ast(|u_{n}|^{2^{\ast}_{\alpha}}-|u|^{2^{\ast}_{\alpha}})|^{2}dx\right)^{\frac{1}{2\cdot2^{\ast}_{\alpha}}}\cdot
\left(\int_{\mathbb{R}^{3}}|I_{\alpha/2}\ast|\varphi|^{2^{\ast}_{\alpha}}|^{2}dx\right)^{\frac{1}{2\cdot2^{\ast}_{\alpha}}}=B_{1}^{\frac{1}{2\cdot2^{\ast}_{\alpha}}}\cdot||\varphi||_{Q^{\alpha,2^{\ast}_{\alpha}}(\mathbb{R}^{3},\mathbb{R}^{3})},
\end{aligned}
\end{equation}
where $B_{1}=\int_{\mathbb{R}^{3}}|I_{\alpha/2}\ast(|u_{n}|^{2^{\ast}_{\alpha}}-|u|^{2^{\ast}_{\alpha}})|^{2}dx$.
Recalling the mean value theorem, we have
\begin{equation}
||u_{n}|^{2^{\ast}_{\alpha}}-|u|^{2^{\ast}_{\alpha}}|=C(2^{\ast}_{\alpha})(|u|+\theta|u_{n}-u|)^{2^{\ast}_{\alpha}-1}|u_{n}-u|=|\xi|^{2^{\ast}_{\alpha}-1}|u_{n}-u|~~\mathrm{for}~0\leq\theta\leq1.
\end{equation}
Therefore, by linearity of the convolution and by positivity of the Riesz-kernel, we deduce that
\begin{equation}
\begin{aligned}
B_{1}&
%=\int_{\mathbb{R}^{3}}|I_{\alpha/2}\ast(|u_{n}|^{2^{\ast}_{\alpha}}-|u|^{2^{\ast}_{\alpha}})|^{2}dx
=\int_{\mathbb{R}^{3}}|I_{\alpha/2}\ast(|\xi|^{2^{\ast}_{\alpha}-1}|u_{n}-u|)|^{2}dx=\int_{\mathbb{R}^{3}}I_{\alpha}\ast(|\xi|^{2^{\ast}_{\alpha}-1}|u_{n}-u|)\cdot (|\xi|^{2^{\ast}_{\alpha}-1}|u_{n}-u|)dx\\
%&=\int_{\mathbb{R}^{3}}\left(I_{\alpha}\ast(|\xi|^{2^{\ast}_{\alpha}-1}|u_{n}-u|)\right)^{\frac{2^{\ast}_{\alpha}-1}{2^{\ast}_{\alpha}}}\left(I_{\alpha}\ast(|\xi|^{2^{\ast}_{\alpha}-1}|u_{n}-u|)\right)^{\frac{1}{2^{\ast}_{\alpha}}}\cdot %(|\xi|^{2^{\ast}_{\alpha}-1}|u_{n}-u|)dx\\
&\leq \left( \int_{\mathbb{R}^{3}}\left(I_{\alpha}\ast(|\xi|^{2^{\ast}_{\alpha}-1}|u_{n}-u|)\right) |\xi|^{2^{\ast}_{\alpha}}dx\right)^{\frac{2^{\ast}_{\alpha}-1}{2^{\ast}_{\alpha}}} \cdot
\left( \int_{\mathbb{R}^{3}}\left(I_{\alpha}\ast(|\xi|^{2^{\ast}_{\alpha}-1}|u_{n}-u|)\right) |u_{n}-u|^{2^{\ast}_{\alpha}}dx\right)^{\frac{1}{2^{\ast}_{\alpha}}}\\
&=\left( \int_{\mathbb{R}^{3}}\left(I_{\alpha/2}\ast(|\xi|^{2^{\ast}_{\alpha}-1}|u_{n}-u|)\right) \cdot  \left(I_{\alpha/2}\ast|\xi|^{2^{\ast}_{\alpha}}\right)dx\right)^{\frac{2^{\ast}_{\alpha}-1}{2^{\ast}_{\alpha}}}\\
&~~~~~~~~~~~~~~~~~~~~\cdot\left( \int_{\mathbb{R}^{3}}\left(I_{\alpha/2}\ast(|\xi|^{2^{\ast}_{\alpha}-1}|u_{n}-u|)\right) \cdot \left(I_{\alpha/2}\ast|u_{n}-u|^{2^{\ast}_{\alpha}}\right)dx\right)^{\frac{1}{2^{\ast}_{\alpha}}}\\
&\leq\left( \int_{\mathbb{R}^{3}}|I_{\alpha/2}\ast(|\xi|^{2^{\ast}_{\alpha}-1}|u_{n}-u|)|^{2} dx\right)^{\frac{1}{2}\frac{2^{\ast}_{\alpha}-1}{2^{\ast}_{\alpha}}}\cdot \left( \int_{\mathbb{R}^{3}}|I_{\alpha/2}\ast|\xi|^{2^{\ast}_{\alpha}}|^{2} dx\right)^{\frac{1}{2}\frac{2^{\ast}_{\alpha}-1}{2^{\ast}_{\alpha}}}\\
&~~~~~~~~~~~~~~~~~~~~\cdot\left( \int_{\mathbb{R}^{3}}|I_{\alpha/2}\ast(|\xi|^{2^{\ast}_{\alpha}-1}|u_{n}-u|)|^{2} dx\right)^{\frac{1}{2}\frac{1}{2^{\ast}_{\alpha}}}\cdot \left( \int_{\mathbb{R}^{3}}|I_{\alpha/2}\ast(|u_{n}-u|^{2^{\ast}_{\alpha}})|^{2} dx\right)^{\frac{1}{2}\frac{1}{2^{\ast}_{\alpha}}}\\
%&=B_{1}^{\frac{1}{2}\frac{2^{\ast}_{\alpha}-1}{2^{\ast}_{\alpha}}}\cdot %||\xi||^{2^{\ast}_{\alpha}-1}_{Q^{\alpha,2^{\ast}_{\alpha}}(\mathbb{R}^{3},\mathbb{R}^{3})}\cdot
%B_{1}^{\frac{1}{2}\frac{1}{2^{\ast}_{\alpha}}}\cdot||u_{n}-u||_{Q^{\alpha,2^{\ast}_{\alpha}}(\mathbb{R}^{3},\mathbb{R}^{3})}\\
&\leq B_{1}^{\frac{1}{2}}\left(||u_{n}||^{2^{\ast}_{\alpha}-1}_{Q^{\alpha,2^{\ast}_{\alpha}}(\mathbb{R}^{3},\mathbb{R}^{3})}+
||u||^{2^{\ast}_{\alpha}-1}_{Q^{\alpha,2^{\ast}_{\alpha}}(\mathbb{R}^{3},\mathbb{R}^{3})}\right)
\cdot||u_{n}-u||_{Q^{\alpha,2^{\ast}_{\alpha}}(\mathbb{R}^{3},\mathbb{R}^{3})}.\\
\end{aligned}
\end{equation}
This implies
\begin{equation}
B_{1}\leq\left(||u_{n}||^{2(2^{\ast}_{\alpha}-1)}_{Q^{\alpha,2^{\ast}_{\alpha}}(\mathbb{R}^{3},\mathbb{R}^{3})}+
||u||^{2(2^{\ast}_{\alpha}-1)}_{Q^{\alpha,2^{\ast}_{\alpha}}(\mathbb{R}^{3},\mathbb{R}^{3})}\right)
\cdot||u_{n}-u||^{2}_{Q^{\alpha,2^{\ast}_{\alpha}}(\mathbb{R}^{3},\mathbb{R}^{3})}.
\end{equation}
Thus,
\begin{equation}
A_{1}\leq \left(||u_{n}||^{\frac{2^{\ast}_{\alpha}-1}{2^{\ast}_{\alpha}}}_{Q^{\alpha,2^{\ast}_{\alpha}}(\mathbb{R}^{3},\mathbb{R}^{3})}+||u||^{\frac{2^{\ast}_{\alpha}-1}{2^{\ast}_{\alpha}}}_{Q^{\alpha,2^{\ast}_{\alpha}}(\mathbb{R}^{3},\mathbb{R}^{3})}\right)\cdot||u_{n}-u||^{\frac{1}{2^{\ast}_{\alpha}}}_{Q^{\alpha,2^{\ast}_{\alpha}}(\mathbb{R}^{3},\mathbb{R}^{3})}
\cdot||\varphi||_{Q^{\alpha,2^{\ast}_{\alpha}}(\mathbb{R}^{3},\mathbb{R}^{3})}.
\end{equation}
Similarly, we have
\begin{equation}
\begin{aligned}
A_{2}&\leq \left(||u_{n}||^{\frac{(2^{\ast}_{\alpha}-1)^{2}}{2^{\ast}_{\alpha}}}_{Q^{\alpha,2^{\ast}_{\alpha}}(\mathbb{R}^{3},\mathbb{R}^{3})}+||u||^{\frac{(2^{\ast}_{\alpha}-1)^{2}}{2^{\ast}_{\alpha}}}_{Q^{\alpha,2^{\ast}_{\alpha}}(\mathbb{R}^{3},\mathbb{R}^{3})}\right)\cdot||u_{n}-u||^{\frac{2^{\ast}_{\alpha}-1}{2^{\ast}_{\alpha}}}_{Q^{\alpha,2^{\ast}_{\alpha}}(\mathbb{R}^{3},\mathbb{R}^{3})}
\cdot||u_{n}||^{2^{\ast}_{\alpha}-1}_{Q^{\alpha,2^{\ast}_{\alpha}}(\mathbb{R}^{3},\mathbb{R}^{3})}\\
A_{3}&\leq||u||_{Q^{\alpha,2^{\ast}_{\alpha}}(\mathbb{R}^{3},\mathbb{R}^{3})}||\varphi||_{Q^{\alpha,2^{\ast}_{\alpha}}(\mathbb{R}^{3},\mathbb{R}^{3})}\\
%A_{4}&\leq %\left(||u_{n}||^{2^{\ast}_{\alpha}-2}_{Q^{\alpha,2^{\ast}_{\alpha}}(\mathbb{R}^{3},\mathbb{R}^{3})}+||u||^{2^{\ast}_{\alpha}-2}_{Q^{\alpha,2^{\ast}_{\alpha}}(\mathbb{R}^{3},\mathbb{R}^{3})}\right)\cdot||u_{n}-u||_{Q^{\alpha,2^{\ast}_{\alpha}}(\mathbb{R}^{3},\mathbb{R}^{3})}
%\cdot||u||^{2^{\ast}_{\alpha}-1}_{Q^{\alpha,2^{\ast}_{\alpha}}(\mathbb{R}^{3},\mathbb{R}^{3})}.\\
A_{4}&\leq \left(||u_{n}||^{2^{\ast}_{\alpha}-2}_{Q^{\alpha,2^{\ast}_{\alpha}}(\mathbb{R}^{3},\mathbb{R}^{3})}+||u||^{2^{\ast}_{\alpha}-2}_{Q^{\alpha,2^{\ast}_{\alpha}}(\mathbb{R}^{3},\mathbb{R}^{3})}\right)\cdot||u_{n}-u||_{Q^{\alpha,2^{\ast}_{\alpha}}(\mathbb{R}^{3},\mathbb{R}^{3})}
\cdot||u||^{2^{\ast}_{\alpha}-1}_{Q^{\alpha,2^{\ast}_{\alpha}}(\mathbb{R}^{3},\mathbb{R}^{3})}.
\end{aligned}
\end{equation}
Therefore, for any $u_{n}\longrightarrow u$ in $Q^{\alpha,2^{\ast}_{\alpha}}(\mathbb{R}^{3},\mathbb{R}^{3})$, we have $\langle(I'(u_{n})-I'(u)),\varphi\rangle\longrightarrow 0$. This implies that $I(u)$ is $\mathcal{C}^{1}$. Therefore, $J_{\lambda}(u)$ and $J(u)$ are of class $\mathcal{C}^{1}$.
\end{proof}

To apply the concentration compactness arguments, we need to introduce the following Coulomb-Sobolev space.
\begin{Def}
Let $\Omega\subset\mathbb{R}^{N}$, $\alpha\in(0,N)$ and $p\geq1$. We define $W^{1,\alpha,p}(\Omega)$ as the scalar space of measurable functions $u:\Omega\longrightarrow \mathbb{R}$ such that $u \in Q^{\alpha,p}(\Omega)$ and $u$ is weakly differentiable in $\Omega$, $Du\in Q^{\alpha,p}(\Omega,\mathbb{R}^{N})$ and
\begin{equation*}
||u||_{W^{1,\alpha,p}(\Omega)}=\left( \left(\int_{\Omega}|I_{\alpha/2}\ast|u|^{p}|^{2}dx\right)^{\frac{1}{2}}+\left(\int_{\Omega}|I_{\alpha/2}\ast|Du|^{p}|^{2}dx\right)^{\frac{1}{2}}\right)^{\frac{1}{p}}<+\infty.
\end{equation*}
\end{Def}
%It differs to the classical Coulomb-Sobolev space $E^{\alpha,p}(\mathbb{R}^{N})$ in \cite{Mercuri2016}.

We are going to prove that the Coulomb-Sobolev space $W^{1,\alpha,p}(\Omega)$ is a Banach space. Firstly, We have the following Fatou property for locally converging sequence.
\begin{lem}\label{fatou}
%Let $N\in\mathbb{N}\setminus\{1\}$, $\alpha\in(0,N)$ and $p\geq1$.
Let $N\in\mathbb{N}$, $\alpha\in(0,N)$ and $p\geq1$. If $(u_{n})_{n\in\mathbb{N}}$ is a bounded sequence in $W^{1,\alpha,p}(\Omega)$ such that $u_{n}\longrightarrow u$ in $L^{1}_{loc}(\Omega)$ and $Du_{n}\longrightarrow g$ in $L^{1}_{loc}(\Omega,\mathbb{R}^{N})$, then $g=Du$ and $u\in W^{1,\alpha,p}(\Omega)$,
\begin{equation}\label{fatou1}
\int_{\Omega}|I_{\alpha/2}\ast| u|^{p}|^{2}dx\leq\mathop{\mathrm{lim~inf}}\limits_{n\longrightarrow\infty}\int_{\Omega}|I_{\alpha/2}\ast| u_{n}|^{p}|^{2}dx,
\end{equation}
and
\begin{equation}\label{fatou2}
\int_{\Omega}|I_{\alpha/2}\ast|Du|^{p}|^{2}dx\leq\mathop{\mathrm{lim~inf}}\limits_{n\longrightarrow\infty}\int_{\Omega}|I_{\alpha/2}\ast|Du_{n}|^{p}|^{2}dx.
\end{equation}
\end{lem}
\begin{proof}
The proof of (\ref{fatou1}) follows the same argument in the proof of (\ref{Fatou1}). We are going to prove (\ref{fatou2}). For $v\in \mathcal{C}_{0}^{\infty}(\Omega)$ we conclude that
\begin{equation}
\int_{\Omega}u_{n}\cdot \mathrm{div}v dx=-\int_{\Omega}\nabla u_{n}\cdot~v dx.
\end{equation}
Since $u_{n}\longrightarrow u$ in $L^{1}_{loc}(\Omega)$, we have
\begin{equation}
\int_{\Omega}u_{n}~\mathrm{div}v~ dx\longrightarrow\int_{\Omega}u~\mathrm{div}v~dx.
\end{equation}
Since $Du_{n}\longrightarrow g$ in $L^{1}_{loc}(\Omega,\mathbb{R}^{N})$, we have
\begin{equation}
-\int_{\Omega}\nabla u_{n}\cdot v~dx\longrightarrow -\int_{\Omega}g\cdot v~ dx.
\end{equation}
Setting $Du$ as the weak derivation of $u$ in the following distribute sense
\begin{equation}
\int_{\Omega}u~\mathrm{div}v~dx=-\int_{\Omega}Du\cdot v~dx,
\end{equation}
we can see $g=Du\in L^{p}(\Omega,\mathbb{R}^{N})$ in the weak sense, and $Du_{n}\longrightarrow Du$ in $L^{1}_{loc}(\Omega,\mathbb{R}^{N})$. Based on this fact, we can obtain (\ref{fatou2}) by the same analysis in the proof of (\ref{Fatou1}).
\end{proof}

\begin{lem}
Let $N\in\mathbb{N}$, $\alpha\in(0,N)$ and $p\geq1$. The normed space $W^{1,\alpha,p}(\Omega)$ is complete.
\end{lem}
\begin{proof}
Let $(u_{n})_{n\in\mathbb{N}}$ be a Cauchy sequence in $W^{1,\alpha,p}(\Omega)$. By the local estimate of the Coulomb energy, $(u_{n})_{n\in\mathbb{N}}$ and $(Du_{n})_{n\in\mathbb{N}}$ are also the Cauchy sequences in $L_{loc}^{p}(\Omega)$. Hence there exists $u\in L_{loc}^{p}(\Omega)$ such that $(u_{n})_{n\in\mathbb{N}}\longrightarrow u$ in $L_{loc}^{p}(\Omega)$ and $g\in L_{loc}^{p}(\Omega,\mathbb{R}^{N})$ such that $(Du_{n})_{n\in\mathbb{N}}\longrightarrow g$ in $L_{loc}^{p}(\Omega,\mathbb{R}^{N})$. In light of Lemma \ref{fatou}, we conclude that $u\in W^{1,\alpha,p}(\Omega)$. Moreover, for every $n\in\mathbb{N}$ the sequence $(u_{n}-u_{m})_{m\in\mathbb{N}}$ converges to $(u_{n}-u)$ in $L_{loc}^{p}(\Omega)$. Hence, by Lemma \ref{fatou} again, we have
\begin{equation*}
\begin{aligned}
&\mathop{\mathrm{lim~sup}}\limits_{n\longrightarrow\infty}\left(\int_{\Omega}|I_{\alpha/2}\ast|Du_{n}-Du|^{p}|^{2}dx+\int_{\Omega}|I_{\alpha/2}\ast|u_{n}-u|^{p}|^{2}dx\right)\\
&\leq\mathop{\mathrm{lim~sup}}\limits_{n\longrightarrow\infty}\mathop{\mathrm{lim~sup}}\limits_{m\longrightarrow\infty}\left(\int_{\Omega}|I_{\alpha/2}\ast|Du_{n}-Du_{m}|^{p}|^{2}dx+\int_{\Omega}|I_{\alpha/2}\ast|u_{n}-u_{m}|^{p}|^{2}dx\right)\\
&\leq\mathop{\mathrm{lim~sup}}\limits_{m,n\longrightarrow\infty}\left(\int_{\Omega}|I_{\alpha/2}\ast|Du_{n}-Du_{m}|^{p}|^{2}dx+\int_{\Omega}|I_{\alpha/2}\ast|u_{n}-u_{m}|^{p}|^{2}dx\right)\\
&\leq0.
\end{aligned}
\end{equation*}
This implies $W^{1,\alpha,p}(\Omega)$ is complete.
\end{proof}
We show that the Coulomb-Sobolev space $W^{1,\alpha,p}(\Omega)$ can be naturally identified with the completion of the set of the test functions $\mathcal{C}_{0}^{\infty}(\Omega)$ under the norm $||\cdot||_{W^{1,\alpha,p}}$.
\begin{lem}\label{dense}
Let $N\in\mathbb{N}$, $\alpha\in(0,N)$ and $p\geq1$. The space of test function $\mathcal{C}_{0}^{\infty}(\Omega)$ is dense in $W^{1,\alpha,p}(\Omega)$.
\end{lem}
\begin{proof}
Since the test function $\mathcal{C}_{0}^{\infty}(\Omega)$ is dense in $Q^{\alpha,p}(\Omega)$, see Proposition 2.6 in \cite{Mercuri2016}, then, by Lemma \ref{fatou} the conclusion also holds in $W^{1,\alpha,p}(\Omega)$.
\end{proof}

Similar to the Poincar\'{e} inequality for the local case, we have the following Poincar\'{e} inequality for the nonlocal case.
 %which in particular implies $W^{1,\alpha,p}(\Omega)\subset Q^{\alpha,p}_{\mathrm{loc}}(\Omega)$.
\begin{lem}\label{poincare}
For all $N\in\mathbb{N}$ and $\alpha\in(0,N)$, there exist $p\in(\frac{N-\alpha}{2},\infty)$ if $\alpha\in(0,N-2)$, while $p\in[N,\infty)$ if $\alpha\in[N-2,N)$, such that for every $a\in\Omega$ and $\rho>0$
\begin{equation*}
\int_{B_{\rho}(a)}|I_{\alpha/2}\ast|u|^{p}|^{2}dx\leq
C\rho^{\frac{N-\alpha}{2}}\bigg(\int_{B_{\rho}(a)}|I_{\alpha/2}\ast|Du|^{p}|^{2}dx\bigg)^{\frac{1}{p}}.
\end{equation*}
%Particularly, we have %$p=2^{\ast}_{\alpha}=\frac{2N-\mu}{N-2}=6-\mu\in(3,6)\subset(\frac{\mu}{2},\infty)$ as $\mu\in(0,3)$ and %$N=3$.
\end{lem}
\begin{proof}
By the HLS inequality (\ref{HLS1}), we have
\begin{equation}\label{261}
\int_{B_{\rho}(a)}|I_{\alpha/2}\ast|u|^{p}|^{2}dx\leq
C_{1}(\alpha,\rho,p,N)\left(\int_{B_{\rho}(a)}|u(x)|^{\frac{2Np}{N+\alpha}}dx\right)^{\frac{N+\alpha}{N}}.
\end{equation}
If $\alpha\in(0,N-2)$ and $p\in(\frac{N-\alpha}{2},N)\subset(1,N)$, then we have
\begin{equation}\label{1h}
\left(\int_{B_{\rho}(a)}|u(x)|^{\frac{2Np}{N+\alpha}}dx\right)^{\frac{N+\alpha}{N}}\leq C_{2}(\alpha,\rho,p,N)\left(\int_{B_{\rho}(a)}|u(x)|^{\frac{Np}{N-p}}dx\right)^{\frac{N-p}{Np}}\leq C_{3}(\alpha,\rho,p,N)\left(\int_{B_{\rho}(a)}|Du|^{p}dx\right)^{\frac{2}{p}}.
\end{equation}
On the other hand, if $\alpha\in(0,N-2)$ and $p\in[N,\infty)$, we know there exists $h\in[\frac{2Np}{N+\alpha+2p},N)$ such that
\begin{equation}\label{h}
\left(\int_{B_{\rho}(a)}|u(x)|^{\frac{2Np}{N+\alpha}}dx\right)^{\frac{N+\alpha}{N}}\leq C_{4}(\alpha,\rho,p,N)\left(\int_{B_{\rho}(a)}|Du|^{h}dx\right)^{\frac{2}{h}}\leq C_{5}(\alpha,\rho,p,N)\left(\int_{B_{\rho}(a)}|Du|^{p}dx\right)^{\frac{2}{p}},
\end{equation}
where the H\"{o}lder inequality was applied. Consequently, for $\alpha\in[N-2,N)$, there also exist $p\in[N,\infty)$ and $h\in[\frac{2Np}{N+\alpha+2p},N)$ such that (\ref{h}) holds.

Then the conclusion follows from \eqref{261}, \eqref{1h}, \eqref{h} and the local estimate of Coulomb energy \cite[Proposition 2.3]{Mercuri2016}, which says that
\begin{equation*}
\left(\int_{B_{\rho}(a)}|Du|^{p}dx\right)^{\frac{2}{p}}\leq C\rho^{\frac{N-\alpha}{2}}\bigg(\int_{B_{\rho}(a)}|I_{\alpha/2}\ast|Du|^{p}|^{2}dx\bigg)^{\frac{1}{p}}.
\end{equation*}

\end{proof}

To establish the Helmholtz decomposition, we also define the following Coulomb-Sobolev space.
\begin{Def}
Let $\Omega\subset\mathbb{R}^{3}$, $\alpha\in(0,3)$ and $p\in(1,\infty)$. We define $W^{1,\alpha,p}_{0}(\mathbb{R}^{3})$ and $W^{1,\alpha,p}_{0}(\Omega)$ as the completion of $C_{0}^{\infty}(\mathbb{R}^{3})$ and $C_{0}^{\infty}(\Omega)$ with respect to the norm
\begin{equation*}
||w||_{W^{1,\alpha,p}_{0}(\mathbb{R}^{3})}=|\nabla w|_{Q^{\alpha,p}(\mathbb{R}^{3})},~~
||w||_{W^{1,\alpha,p}_{0}(\Omega)}=|\nabla w|_{Q^{\alpha,p}(\Omega)}.
\end{equation*}
%Namely, %$W^{1,\alpha,p}_{0}(\mathbb{R}^{3})=\overline{C_{0}^{\infty}(\mathbb{R}^{3})}^{|\nabla\cdot|_{Q^{\alpha,p}(\mathbb{R}^{3})}}$.
\end{Def}
\begin{Prop}
$W^{1,\alpha,p}_{0}(\mathbb{R}^{3})$ is linearly isometric to
\begin{equation*}
\nabla W^{1,\alpha,p}_{0}(\mathbb{R}^{3}):=\{\nabla w\in Q^{\alpha,p}(\mathbb{R}^{3},\mathbb{R}^{3}): w\in W^{1,\alpha,p}_{0}(\mathbb{R}^{3})\},
\end{equation*}
and  $W^{1,\alpha,p}_{0}(\Omega)$ is linearly isometric to
\begin{equation*}
\nabla W^{1,\alpha,p}_{0}(\Omega):=\{\nabla w\in Q^{\alpha,p}(\Omega,\mathbb{R}^{3}): w\in W^{1,\alpha,p}_{0}(\Omega)\}.
\end{equation*}
\end{Prop}
\begin{proof}
Set the map $\nabla:W^{1,\alpha,p}_{0}(\mathbb{R}^{3})\longrightarrow \nabla W^{1,\alpha,p}_{0}(\mathbb{R}^{3})$. Since the Coulomb space is complete, the map is obviously injective and surjective. We also easily check that the map is isometric by the definition of $W^{1,\alpha,p}_{0}(\mathbb{R}^{3})$, this implies our conclusion.
\end{proof}

\subsubsection{Helmholtz decomposition}
\ \ \\

Let $\mathcal{D}^{1,2}(\mathbb{R}^{3},\mathbb{R}^{3})$ denote the completion of $\mathcal{C}^{\infty}_{0}(\mathbb{R}^{3},\mathbb{R}^{3})$ with respect to the norm $|\nabla\cdot|_{2}$. Recall the subspace $\mathcal{V}_{\mathbb{R}^{3}}$ and $\mathcal{W}_{\mathbb{R}^{3}}$ of $W^{\alpha,2^{\ast}_{\alpha}}_{0}(\mathrm{curl};\mathbb{R}^{3})$ in the introduction, we have the following Helmholtz decomposition on $W^{\alpha,2^{\ast}_{\alpha}}_{0}(\mathrm{curl};\mathbb{R}^{3})$.
\begin{lem}\label{Helmholtz}
$\mathcal{V}_{\mathbb{R}^{3}}$ and $\mathcal{W}_{\mathbb{R}^{3}}$ are closed subspaces of $W^{\alpha,2^{\ast}_{\alpha}}_{0}(\mathrm{curl};\mathbb{R}^{3})$ and
\begin{equation}\label{Helmholtz1}
W^{\alpha,2^{\ast}_{\alpha}}_{0}(\mathrm{curl};\mathbb{R}^{3})=\mathcal{V}_{\mathbb{R}^{3}}\oplus\nabla W^{1,\alpha,2^{
\ast}_{\alpha}}_{0}(\mathbb{R}^{3})=\mathcal{V}_{\mathbb{R}^{3}}\oplus\mathcal{W}_{\mathbb{R}^{3}}. ~~(direct~sum)
\end{equation}
Moreover, $\mathcal{V}_{\mathbb{R}^{3}}\subset \mathcal{D}^{1,2}(\mathbb{R}^{3},\mathbb{R}^{3})$ and the norms $|\nabla\cdot|_{2}$ and $||\cdot||_{W^{\alpha,2^{\ast}_{\alpha}}(\mathrm{curl};\mathbb{R}^{3})}$ are equivalent in $\mathcal{V}_{\mathbb{R}^{3}}$.
\end{lem}
\begin{proof}
By the HLS inequality in Proposition \ref{HLS}, there is a continuous embedding
\begin{equation*}
L^{2^{\ast}}(\mathbb{R}^{3},\mathbb{R}^{3})\hookrightarrow Q^{\alpha,2^{\ast}_{\alpha}}(\mathbb{R}^{3},\mathbb{R}^{3}).
\end{equation*}
Then the conclusion follows from the argument in \cite[Lemma 3.2]{Mederski2015}. Indeed, Since $W^{1,\alpha,2^{\ast}_{\alpha}}_{0}(\mathbb{R}^{3})$ is a complete space, then $\nabla W^{1,\alpha,2^{\ast}_{\alpha}}_{0}(\mathbb{R}^{3})$ is a closed subspace of $Q^{\alpha,2^{\ast}_{\alpha}}(\mathbb{R}^{3},\mathbb{R}^{3})$. Moreover $\mathrm{cl}\mathcal{V}_{\mathbb{R}^{3}}\cap\nabla W^{1,\alpha,2^{\ast}_{\alpha}}_{0}(\mathbb{R}^{3})=\{0\}$ in $Q^{\alpha,2^{\ast}_{\alpha}}(\mathbb{R}^{3},\mathbb{R}^{3})$, hence $\mathcal{V}_{\mathbb{R}^{3}}\cap\nabla W^{1,\alpha,2^{\ast}_{\alpha}}_{0}(\mathbb{R}^{3})=\{0\}$ in $W^{\alpha,2^{\ast}_{\alpha}}(\mathrm{curl};\mathbb{R}^{3})$. In view of the Helmholtz decomposition, and smooth function $\varphi\in\mathcal{C}_{0}^{\infty}(\mathbb{R}^{3},\mathbb{R}^{3})$ can be written as
\begin{equation*}
\varphi=\varphi_{1}+\nabla\varphi_{2}
\end{equation*}
such that $\varphi_{1}\in\mathcal{D}^{1,2}(\mathbb{R}^{3},\mathbb{R}^{3})\cap\mathcal{C}^{\infty}(\mathbb{R}^{3},\mathbb{R}^{3})$, $\mathrm{div}(\varphi_{1})=0$ and $\varphi_{2}\in\mathcal{C}^{\infty}(\mathbb{R}^{3})$ is the Newton potential of $\mathrm{div}(\varphi)$. Since $\varphi$ has compact support, then $\nabla\varphi_{2}\in L^{6}(\mathbb{R}^{3},\mathbb{R}^{3})\subset Q^{\alpha,2^{\ast}_{\alpha}}(\mathbb{R}^{3},\mathbb{R}^{3})$ and $\varphi_{1}=\varphi-\nabla\varphi_{2}\in\mathcal{V}_{\mathbb{R}^{3}}$. Observe that $\nabla\times\nabla\varphi_{1}=-\triangle\varphi_{1}$, hence
\begin{equation*}
|\nabla\times u|_{2}=|\nabla u|_{2}=||u||_{\mathcal{D}^{1,2}(\mathbb{R}^{3},\mathbb{R}^{3})}
\end{equation*}
for any $u\in\mathcal{V}_{\mathbb{R}^{3}}$. By the Sobolev embedding we have $\mathcal{V}_{\mathbb{R}^{3}}$ is continuously embedded in $L^{6}(\mathbb{R}^{3},\mathbb{R}^{3})$ and by the HLS inequality also in $Q^{\alpha,2^{\ast}_{\alpha}}(\mathbb{R}^{3},\mathbb{R}^{3})$. Therefore the norms $||\cdot||_{\mathcal{D}^{1,2}(\mathbb{R}^{3},\mathbb{R}^{3})}$ and $||\cdot||_{W^{\alpha,2^{\ast}_{\alpha}}}(\mathrm{curl};\mathbb{R}^{3})$ are equivalent on $\mathcal{V}_{\mathbb{R}^{3}}$ and by the density argument we get the decomposition (\ref{Helmholtz1}).
\end{proof}

For the bounded domains case, we recall the definition of $\mathcal{V}'_{\Omega}$ in \cite{Bartsch2015-1}, that is
\begin{equation*}
\begin{aligned}
\mathcal{V}'_{\Omega}&=\{v\in W^{2}_{0}(\mathrm{curl};\Omega):\int_{\Omega}\langle v,\varphi\rangle dx=0~\mathrm{for~every}~\varphi\in C^{\infty}_{0}(\Omega,\mathbb{R}^{3})~\mathrm{with}~\nabla\times\varphi=0\}.\\
%&=\{v\in W^{2}_{0}(\mathrm{curl};\Omega):\mathrm{div}~v=0~\mathrm{in~the~sense~of~distributions}\}
\end{aligned}
\end{equation*}
Indeed, if $\varphi$ is supported in a ball, we have $\varphi=\nabla\psi$
for some $\psi\in\mathcal{C}_{0}^{\infty}(\Omega)$, hence we have $\mathrm{div}(v)=0$ for $v\in\mathcal{V}'_{\Omega}$. This implies that
%we set $H(\mathrm{curl};\Omega)=\{u\in L^{2}(\Omega,\mathbb{R}^{3},\nabla\times u\in L^{2}(\Omega,\mathbb{R}^{3}))\}$ and %$H_{0}(\mathrm{curl};\Omega)=\overline{C_{0}^{\infty}(\Omega,\mathbb{R}^{3})}^{||\cdot||_{H(\mathrm{curl};\Omega)}}$, then
\begin{equation*}\label{XN}
\begin{aligned}
\mathcal{V}'_{\Omega}&=\{v\in W_{0}^{2}(\mathrm{curl};\Omega):\mathrm{div}~v=0~\mathrm{in~the~sense~of~distributions}\}\\
&\subset \{v\in W^{2}_{0}(\mathrm{curl};\Omega):\mathrm{div}~v\in L^{2}(\Omega,\mathbb{R}^{3})\}=:X_{N}(\Omega).
\end{aligned}
\end{equation*}
Furthermore, since $\Omega$ is a bounded domain, either convex or with $\mathcal{C}^{1,1}$ boundary, $X_{N}(\Omega)$ is continuously embedded in $H^{1}(\Omega,\mathbb{R}^{3})$, see \cite{Amrouche1998}. Therefore in view of the Rellich's theorem $\mathcal{V}'_{\Omega}$ is compactly embedded in $L^{2}(\Omega,\mathbb{R}^{3})$ and continuously in $L^{6}(\Omega,\mathbb{R}^{3})$, so is $Q^{\alpha,2^{\ast}_{\alpha}}(\Omega,\mathbb{R}^{3})$.
%Therefore, $\mathcal{V}_{\Omega}$ is compactly embedded in $Q^{\alpha,p}(\Omega,\mathbb{R}^{3})$ for %$\frac{3+\alpha}{3}\leq p<2^{\ast}_{\alpha}=3+\alpha$ and continuously in %$Q^{\alpha,2^{\ast}_{\alpha}}(\Omega,\mathbb{R}^{3})$.
This implies in particular that $\mathcal{V}'_{\Omega}\subset \mathcal{V}_{\Omega}$. On the other hand, since $W^{\alpha,2^{\ast}_{\alpha}}_{0}(\mathrm{curl};\Omega)\subset W^{2}_{0}(\mathrm{curl};\Omega)$, we have $\mathcal{V}_{\Omega}\subset \mathcal{V}'_{\Omega}$. Therefore, we can see that $\mathcal{V}_{\Omega}=\mathcal{V}'_{\Omega}$
%\begin{equation*}
%\begin{aligned}
%\mathcal{V}_{\Omega}&=\{v\in W^{2}_{0}(\mathrm{curl};\Omega):\int_{\Omega}\langle v,\varphi\rangle %dx=0~\mathrm{for~every}~\varphi\in %C^{\infty}_{0}(\Omega,\mathbb{R}^{3})~\mathrm{with}~\nabla\times\varphi=0\}\\
%&=\{v\in W^{2}_{0}(\mathrm{curl};\Omega):\mathrm{div}~v=0~\mathrm{in~the~sense~of~distributions}\}
%\end{aligned}
%\end{equation*}
is a Hilbert space with  inner product
\begin{equation*}
(v,z)=\int_{\Omega}\langle \nabla\times v,\nabla\times z\rangle dx=\int_{\Omega}\langle\nabla v,\nabla z\rangle dx.
\end{equation*}
Also, one can easily observe that $\mathcal{V}_{\Omega}$ is a closed linear subspace of $W^{\alpha,2^{\ast}_{\alpha}}_{0}(\mathrm{curl};\Omega)$. Therefore, by Theorem 4.21 (c) in \cite{Kirsch2015}, we have the following Helmholtz decomposition
\begin{equation}\label{bounded helmholtz}
W^{\alpha,2^{\ast}_{\alpha}}_{0}(\mathrm{curl};\Omega)=\mathcal{V}_{\Omega}\oplus\mathcal{W}_{\Omega}. ~~(\mathrm{direct~sum})
\end{equation}
and that
\begin{equation}\label{orth}
\int_{\Omega}\langle v,w\rangle dx=0~~\mathrm{if}~v\in\mathcal{V}_{\Omega},~w\in \mathcal{W}_{\Omega},
\end{equation}
which means that $\mathcal{V}_{\Omega}$ and $\mathcal{W}_{\Omega}$ are orthogonal in $L^{2}(\Omega,\mathbb{R}^{3})$. Then the norm
\begin{equation*}
||v+w||:=((v,v)+|w|_{Q^{\alpha,2^{\ast}_{\alpha}}}^{2})^{\frac{1}{2}},~v\in\mathcal{V}_{\Omega},~w\in \mathcal{W}_{\Omega}
\end{equation*}
is equivalent to $||\cdot||_{W^{\alpha,2^{\ast}_{\alpha}}_{0}(\mathrm{curl};\Omega)}$.

For the setting of boundary condition, according to \cite[Theorem 3.33]{Monk2003}, there is a continuous tangential trace operator $\gamma_{t}:W^{2}(\mathrm{curl};\Omega)\longrightarrow H^{-1/2}(\partial\Omega)$ such that
\begin{equation*}
\gamma_{t}(u)=\nu\times u|_{\partial\Omega}~~~~\mathrm{for}~\mathrm{any}~u\in C^{\infty}(\overline{\Omega},\mathbb{R}^{3})
\end{equation*}
and
\begin{equation*}
W^{2}_{0}(\mathrm{curl};\Omega)=\{u\in W^{2}(\mathrm{curl};\Omega):\gamma_{t}(u)=0\}.
\end{equation*}
 Hence the vector field $u\in W^{\alpha,2^{\ast}_{\alpha}}_{0}(\mathrm{curl};\Omega)=\mathcal{V}_{\Omega}\oplus\mathcal{W}_{\Omega}\subset W^{2}_{0}(\mathrm{curl};\Omega)$ satisfies the boundary condition in (\ref{nonlocal case}).

On the other hand, $\mathcal{W}_{\Omega}$ contains all gradient vectors fields, i.e. $\nabla W^{1,\alpha,2^{\ast}_{\alpha}}_{0}(\Omega)\subset \mathcal{W}_{\Omega}$. However, for some general domains, $\{w\in \mathcal{W}_{\Omega}:\mathrm{div}(w)=0\}$ may be nontrivial (harmonic field) and hence $\nabla W^{1,\alpha,2^{\ast}_{\alpha}}_{0}(\Omega)\subsetneq\mathcal{W}_{\Omega}$, see \cite[pp.4314-4315]{Bartsch2017}. While in the topology domains as we supposed, we have the following conclusion, which is a trivial extended from Lemma 2.3 in \cite{Mederski2021}.
\begin{lem} \label{gradient}
There holds $\mathcal{W}_{\Omega}=W^{\alpha,2^{\ast}_{\alpha}}_{0}(\mathrm{curl};\Omega)\cap \mathcal{W}_{\mathbb{R}^{3}}=W^{\alpha,2^{\ast}_{\alpha}}_{0}(\mathrm{curl};\Omega)\cap \nabla W^{1,\alpha,2^{\ast}_{\alpha}}(\Omega)$. If $\partial\Omega$ is connected, then $\mathcal{W}_{\Omega}=\nabla W^{1,\alpha,2^{\ast}_{\alpha}}_{0}(\Omega)$. If $\Omega$ is unbounded, $\mathcal{W}_{\Omega}=W^{\alpha,2^{\ast}_{\alpha}}_{0}(\mathrm{curl};\Omega)\cap \mathcal{W}_{\mathbb{R}^{3}}$ still holds.
\end{lem}

\subsection{Concentration Compactness Lemma}
In view of the Helmholtz decomposition, the work space is decomposed into a Hilbert space $\mathcal{V}_{\Omega}$ and a Banach space $\mathcal{W}_{\Omega}$. For a bounded sequence in the work space, one can obtain the a.e convergency in $\mathcal{V}_{\Omega}$ by the Rellich compactness theorem, which is important to the weak-weak$^{\ast}$ continuity of $J'(u)$. While in the subspace $\mathcal{W}_{\Omega}$, $w_{n}=\nabla p_{n}\rightharpoonup \nabla p=w$ can not deduce the a.e convergency. By setting the convex nonlinearity satisfied the coercive condition, Merderski \cite{Mederski2021} connected the subspaces $\mathcal{V}_{\Omega}$ and $\mathcal{W}_{\Omega}$ by the global minimum argument, then the a.e. convergency on $\mathcal{W}_{\Omega}$ can be recovered by the second concentration compactness lemma, see Lions \cite{Lions1985}. Since the nonlinearity becomes a nonlocal term, we make some minor modifications to the concentration compactness lemma.

In this subsection, We work in some subspaces of $Q^{\alpha,2^{\ast}_{\alpha}}(\Omega,\mathbb{R}^{3})$ and $Q^{\alpha,2^{\ast}_{\alpha}}(\mathbb{R}^{3},\mathbb{R}^{3})$.
Let $Z\subset \mathcal{V}_{\Omega}$ be a finite-dimension subspace of $Q^{\alpha,2^{\ast}_{\alpha}}(\Omega,\mathbb{R}^{N})$ such that $Z\cap\mathcal{W}_{\Omega}=\{0\}$ and put
\begin{equation*}
\widetilde{\mathcal{W}}:=\mathcal{W}_{\Omega}\oplus Z.
\end{equation*}
Correspondingly, in $\mathbb{R}^{3}$, we put $Z=\{0\}$ and $\widetilde{\mathcal{W}}=\mathcal{W}_{\mathbb{R}^{3}}$. For simplicity, we only show the discussion on bounded domains $\Omega$, and the case in the entire space $\mathbb{R}^{3}$ is similar. Note that we always assume that $v\in\mathcal{V}_{\mathbb{R}^{3}}\subset\mathcal{D}^{1,2}(\mathbb{R}^{3},\mathbb{R}^{3})$ but not $\mathcal{V}_{\Omega}$, we then have
\begin{lem}\label{global mini}
Assume $F(u)=(I_{\alpha}\ast|u|^{2^{\ast}_{\alpha}})|u|^{2^{\ast}_{\alpha}}$ and $f(u)=\partial_{u}F(u)$, then $F(u)$ is uniformly strictly convex with respect to $u\in\mathbb{R}^{N}$, i.e. for any compact $A\subset(\mathbb{R}^{3}\times\mathbb{R}^{3})\setminus\{(u,u):u\in\mathbb{R}^{3}\}$
\begin{equation}\label{strictly}
\mathop{\mathrm{inf}}\limits_{x\in\Omega,(u_{1},u_{2})\in A}\left(\frac{1}{2}\left(F(u_{1})+F(u_{2})\right)-F\left(\frac{u_{1}+u_{2}}{2}\right)\right)>0;
\end{equation}
Moreover, for any $v\in\mathcal{V}_{\mathbb{R}^{3}}$ we find a unique $\widetilde{w}_{\Omega}(v)\in\widetilde{\mathcal{W}}$ such that
\begin{equation}\label{mini1}
\int_{\Omega}F(v+\widetilde{w}_{\Omega}(v))dx\leq\int_{\Omega}F(v+\widetilde{w})dx ~~\mathrm{for}~\mathrm{all}~\widetilde{w}\in\widetilde{\mathcal{W}}.
\end{equation}
In other word,
\begin{equation}\label{mini2}
\int_{\Omega}\langle f(v+\widetilde{w}),\zeta\rangle dx=0~~\mathrm{for}~\mathrm{all}~\zeta\in\widetilde{\mathcal{W}}~\mathrm{if}~\mathrm{and}~\mathrm{only}~\mathrm{if}~\widetilde{w}=\widetilde{w}_{\Omega}(v).
\end{equation}
\end{lem}
\begin{proof}
The uniformly convexity of $F(u)$ follows from the Proposition 2.8 in \cite{Mercuri2016}. Now, we prove that $F(u)$ is strictly convex. Set $I(u)=\int_{\mathbb{R}^{3}}F(u)dx$ and  $u(x)=(u_{1},u_{2},u_{3})$, then for any $(s_{1},s_{2},s_{3})\in \mathbb{R}^{3}$ we have
\begin{equation*}
I(u)=\int_{\Omega}|I_{\alpha/2}\ast|u|^{2^{\ast}_{\alpha}}|^{2}dx=\int_{\Omega}\left[I_{\alpha/2}\ast\left(\mathop{\sum}\limits_{i=1}^{3}s_{i}|u_{i}|^{2^{\ast}_{\alpha}}\right)\right]^{2}dx
=\int_{\Omega}\left[\mathop{\sum}\limits_{i=1}^{3}s_{j}\left(I_{\alpha/2}\ast |u_{i}|^{2^{\ast}_{\alpha}}\right)\right]^{2}dx
\end{equation*}
Set
\begin{equation*}
g(s_{1},s_{2},s_{3})=\left[\mathop{\sum}\limits_{i=1}^{3}s_{j}\left(I_{\alpha/2}\ast |u_{i}|^{2^{\ast}_{\alpha}}\right)\right]^{2}=h(L(s_{1},s_{2},s_{3})),
\end{equation*}
where $h(t)=t^{2}$ ia a strict convex function and
\begin{equation*}
L(s_{1},s_{2},s_{3})=\mathop{\sum}\limits_{i=1}^{3}s_{j}\left(I_{\alpha/2}\ast |u_{i}|^{2^{\ast}_{\alpha}}\right).
\end{equation*}
is a linear functional. Then, for each $x\in\mathbb{R}^{3}$, $g(s_{1},s_{2},s_{3})$ is convex.

Indeed, fix $\lambda\in(0,1)$ and $(s_{1},s_{2},s_{3}),(r_{1},r_{2},r_{3})\in\mathbb{R}^{3}$, we have
\begin{equation*}
\begin{aligned}
&g((1-\lambda)(s_{1},s_{2},s_{3})+\lambda(r_{1},r_{2},r_{3}))
=h(L((1-\lambda)(s_{1},s_{2},s_{3})+\lambda(r_{1},r_{2},r_{3})))\\
&=h((1-\lambda)L(s_{1},s_{2},s_{3})+\lambda L(r_{1},r_{2},r_{3}))\leq(1-\lambda)h(L(s_{1},s_{2},s_{3}))+\lambda h(L(r_{1},r_{2},r_{3}))\\
&=(1-\lambda)g(s_{1},s_{2},s_{3})+\lambda g(r_{1},r_{2},r_{3}).
\end{aligned}
\end{equation*}
Moreover, since $L$ is an injective function, we deduce that $g$ is strictly convex. Hence, $I(u)$ is strictly convex, so is $F(u)$. On the other hand, $I(u)$ is coercive in $Q^{\alpha,2^{\ast}_{\alpha}}(\Omega,\mathbb{R}^{3})$. Then, by the global minimum theorem, we have (\ref{mini1}) and (\ref{mini2}).
\end{proof}

Denote the space of finite measures in $\mathbb{R}^{3}$ by $\mathcal{M}(\mathbb{R}^{3})$. Then we have the following concentration-compactness lemma, see \cite[Lemma 3.1]{Mederski2021} for the local case.
\begin{lem}\label{concentration}
Assume $F(u)=(I_{\alpha}\ast|u|^{2^{\ast}_{\alpha}})|u|^{2^{\ast}_{\alpha}}$. Suppose  $(v_{n})\subset \mathcal{V}_{\mathbb{R}^{3}}$, $v_{n}\rightharpoonup v_{0}$ in $\mathcal{V}_{\mathbb{R}^{3}}$, $v_{n}\longrightarrow v_{0}$ a.e. in $\mathbb{R}^{3}$, $|\nabla v_{n}|^{2}\rightharpoonup \mu$ and $\left(I_{\alpha}\ast|v_{0}|^{2^{\ast}_{\alpha}}\right)|v_{0}|^{2^{\ast}_{\alpha}}\rightharpoonup\rho$ in $\mathcal{M}(\mathbb{R}^{3})$. Then there exists an at most countable set $I\subset \mathbb{R}^{3}$ and nonnegative weights $\{\mu_{x}\}_{x\in I}$, $\{\rho_{x}\}_{x\in I}$ such that
\begin{equation*}
\mu\geq |\nabla v_{0}|^{2}+\mathop{\Sigma}\limits_{x\in I}\mu_{x}\delta_{x},~~\rho=\left(I_{\alpha}\ast|v_{0}|^{2^{\ast}_{\alpha}}\right)|v_{0}|^{2^{\ast}_{\alpha}}+\mathop{\Sigma}\limits_{x\in I}\rho_{x}\delta_{x},
\end{equation*}
and passing to a subsequence, $\widetilde{w}_{\Omega}(v_{n})\rightharpoonup \widetilde{w}_{\Omega}(v_{0})$ in $\widetilde{\mathcal{W}}$, $\widetilde{w}_{\Omega}(v_{n})\longrightarrow \widetilde{w}_{\Omega}(v_{0})$ a.e. in $\Omega$ and in $L^{p}_{loc}(\Omega)$ for any $1\leq p\leq2^{\ast}_{\alpha}$.
\end{lem}
\begin{Rem}
If $\Omega=\mathbb{R}^{3}$, $\widetilde{\mathcal{W}}=\mathcal{W}_{\mathbb{R}^{3}}$, we have the same conclusion, that is $\widetilde{w}_{\mathbb{R}^{3}}(v_{n})\rightharpoonup \widetilde{w}_{\mathbb{R}^{3}}(v_{0})$ in $\widetilde{\mathcal{W}}$, $\widetilde{w}_{\mathbb{R}^{3}}(v_{n})\longrightarrow \widetilde{w}_{\mathbb{R}^{3}}(v_{0})$ a.e. in $\mathbb{R}^{3}$.
\end{Rem}
\begin{proof}
Step 1. Let $\varphi\in \mathcal{C}^{\infty}_{0}(\mathbb{R}^{3})$, then by the definition of $S_{HL}$ in (\ref{SHL}), we have
\begin{equation*}
\begin{aligned}
\left(\int_{\mathbb{R}^{3}}|I_{\alpha/2}\ast|\varphi(v_{n}-v_{0})|^{2^{\ast}_{\alpha}}|^{2}dx\right)^{\frac{1}{2^{\ast}_{\alpha}}} S_{HL}&\leq \int_{\mathbb{R}^{3}}|\nabla(\varphi (v_{n}-v_{0}))|^{2}dx.
\end{aligned}
\end{equation*}
This means that
\begin{equation*}
\begin{aligned}
\left(\int_{\mathbb{R}^{3}}|\varphi|^{2\cdot2^{\ast}_{\alpha}}\left(I_{\alpha}\ast|v_{n}-v_{0}|^{2^{\ast}_{\alpha}}\right)|(v_{n}-v_{0})|^{2^{\ast}_{\alpha}}dx\right)^{\frac{1}{2^{\ast}_{\alpha}}} S_{HL}&\leq \int_{\mathbb{R}^{3}}|\varphi|^{2}|\nabla(v_{n}-v_{0})|^{2}dx+o(1).
\end{aligned}
\end{equation*}
Using the Brezis-Lieb lemma for the nonlocal case on the left-hand side, see \cite[pp.1226]{Gao2018}, we then obtain
\begin{equation}\label{sobolev}
\left(\int_{\mathbb{R}^{3}}|\varphi|^{2\cdot 2^{\ast}_{\alpha}}d\overline{\rho}\right)^{\frac{1}{2^{\ast}_{\alpha}}}S_{HL}\leq \left(\int_{\mathbb{R}^{3}}|\varphi|^{2}d\overline{\mu}\right)^{1/2},
\end{equation}
where $\overline{\mu}:=\mu-|\nabla v_{0}|^{2}$ and $\overline{\rho}=\rho-\left(I_{\alpha}\ast|v_{0}|^{2^{\ast}_{\alpha}}\right)|v_{0}|^{2^{\ast}_{\alpha}}$. Set $I=\{x\in \mathbb{R}^{3}:\mu(\{x\})>0\}$. Since $\mu$ is finite and $\mu,\overline{\mu}$ have the same singular set, $I$ is at most countable and $\mu\geq|\nabla v_{0}|^{2}+\Sigma_{x\in I}\mu_{x}\delta_{x}$. As in the proof of Lemma 2.5 in \cite{Gao2017} it follows from (\ref{sobolev}) that $\overline{\rho}=\Sigma_{x\in I}\rho_{x}\delta_{x}$. So $\mu$ and $\rho$ are as claimed.

Step 2. To recover the a.e. convergency of the sequence on $\mathcal{W}_{\Omega}$, we consider the global minimum argument which connects $\mathcal{V}_{\mathbb{R}^{3}}$ and $\mathcal{W}_{\Omega}$. Using (\ref{mini1}) we infer that
\begin{equation}\label{mini3}
\begin{aligned}
|v_{n}+\widetilde{w}_{\Omega}(v_{n})|^{2\cdot2^{\ast}_{\alpha}}_{Q^{\alpha,2^{\ast}_{\alpha}}}&\leq\int_{\Omega}F(v_{n}+\widetilde{w}_{\Omega}(v_{n}))\leq\int_{\Omega}F(v_{n})dx\leq |v_{n}|_{Q^{\alpha,2^{\ast}_{\alpha}}}^{2\cdot2^{\ast}_{\alpha}}.
\end{aligned}
\end{equation}
Since the right-hand side above is bounded, so is $(|\widetilde{w}_{\Omega}(v_{n})|_{Q^{\alpha,2^{\ast}_{\alpha}}})$. Hence, by the uniform convexity and reflexivity of Coulomb space, see \cite[Section 2.4.1]{Mercuri2016}, up to a subsequence, $\widetilde{w}_{\Omega}(v_{n})\rightharpoonup\widetilde{w}_{0}$ for some $\widetilde{w}_{0}$.

In the following we are going to prove that $\widetilde{w}_{\Omega}(v_{n})\longrightarrow \widetilde{w}_{0}$ a.e. in $\Omega$ after taking subsequence. The convexity of $F$ in $u$ implies that
\begin{equation*}
	F\left(\frac{u_{1}+u_{2}}{2}\right)\geq F(u_{1})+\left\langle f(u_{1}),\frac{u_{2}-u_{1}}{2}\right\rangle,
\end{equation*}
applying (\ref{strictly}), we obtain for any $k\geq 1$ and $|u_{1}-u_{2}|\geq \frac{1}{k}$, $|u_{1}|,|u_{2}|\leq k$ that
\begin{equation}\label{mk}
	m_{k}\leq\frac{1}{2}(F(u_{1})+F(u_{2}))-F\left(\frac{u_{1}+u_{2}}{2}\right)\leq\frac{1}{4}\langle f(u_{1})-f(u_{2}),u_{1}-u_{2}\rangle,
\end{equation}
where
\begin{equation*}
	m_{k}:=\mathop{\mathrm{inf}}\limits_{x\in\Omega,u_{1},u_{2}\in\mathbb{R}^{3}}\frac{1}{2}(F(u_{1})+F(u_{2}))-F\left(\frac{u_{1}+u_{2}}{2}\right)>0~~\mathrm{for}~~\frac{1}{k}\leq|u_{1}-u_{2}|,~|u_{1}|,|u_{2}|\leq k.
\end{equation*}
Now we decompose by  $\widetilde{w}_{\Omega}(v_{n})=w_{n}+z_{n}$, $\widetilde{w}_{0}=w_{0}+z_{0}$ where $w_{n}$, $w_{0}\in\mathcal{W}_{\Omega}$ and $z_{n}$, $z_{0}\in Z$. Obviously, since $Z$ is a finite dimension space, we may assume $z_{n}\longrightarrow z_{0}$ in $Z$ and a.e. in $\Omega$.  Notice that $v_{n}+\widetilde{w}_{\Omega}(v_{n})$ is bounded in $Q^{\alpha,2^{\ast}_{\alpha}}(\Omega,\mathbb{R}^{3})$, we may introduce
\begin{equation}\label{omega}
\Omega_{n,k}:=\{x\in\Omega:|v_{n}+\widetilde{w}_{\Omega}(v_{n})-v_{0}-w_{0}-z_{0}|\geq\frac{1}{k}~\mathrm{and}~|v_{n}+\widetilde{w}_{\Omega}(v_{n})|,|v_{0}+w_{0}+z_{0}|\leq k\}.
\end{equation}
Then, by (\ref{omega}) and (\ref{mk}), we have
\begin{equation}\label{measure1}
\begin{aligned}
&4m_{k}\int_{\Omega_{n,k}}|\varphi|^{2\cdot2^{\ast}_{\alpha}}dx\\
&\leq\int_{\Omega}|\varphi|^{2\cdot2^{\ast}_{\alpha}}\langle f(v_{n}+\widetilde{w}_{\Omega}(v_{n}))-f(v_{0}+w_{0}+z_{0}),v_{n}+\widetilde{w}_{\Omega}(v_{n})-v_{0}-w_{0}-z_{0}  \rangle dx\\
&=\int_{\Omega}|\varphi|^{2\cdot2^{\ast}_{\alpha}}\langle f(v_{n}+\widetilde{w}_{\Omega}(v_{n}))-f(v_{0}+w_{0}+z_{0}),v_{n}-v_{0} \rangle dx\\
&~~~+\int_{\Omega}|\varphi|^{2\cdot2^{\ast}_{\alpha}}\langle f(v_{n}+\widetilde{w}_{\Omega}(v_{n}))-f(v_{0}+w_{0}+z_{0}),\widetilde{w}_{\Omega}(v_{n})-w_{0}-z_{0}\rangle dx=I_{1}+I_{2}.
\end{aligned}
\end{equation}
 Since $|v_{n}+\widetilde{w}_{\Omega}(v_{n})|\leq k$ and $|v_{0}+w_{0}+z_{0}|\leq k$ on $\Omega_{n,k}$, we have $|v_{n}+\widetilde{w}_{\Omega}(v_{n})|\leq C_{1}|v_{n}|$ and $|v_{0}+w_{0}+z_{0}|\leq C_{2}|v_{0}|$.
Then, by the similar estimation in (iii) of Lemma \ref{well define} we have
\begin{equation}\label{measure4}
\begin{aligned}
I_{1}&=\int_{\Omega}|\varphi|^{2\cdot2^{\ast}_{\alpha}}\Big\langle (I_{\alpha}\ast|v_{n}+\widetilde{w}_{\Omega}(v_{n})|^{2^{\ast}_{\alpha}})|v_{n}+\widetilde{w}_{\Omega}(v_{n})|^{2^{\ast}_{\alpha}-2}(v_{n}+\widetilde{w}_{\Omega}(v_{n}))\\
&~~~~~~~~~~~~~~~~~~-(I_{\alpha}\ast|v_{0}+w_{0}+z_{0}|^{2^{\ast}_{\alpha}})|v_{0}+w_{0}+z_{0}|^{2^{\ast}_{\alpha}-2}(v_{0}+w_{0}+z_{0}),v_{n}-v_{0}           \Big\rangle dx\\
&=\int_{\Omega}|\varphi|^{2\cdot2^{\ast}_{\alpha}}\Big\langle (I_{\alpha}\ast|v_{n}+\widetilde{w}_{\Omega}(v_{n})|^{2^{\ast}_{\alpha}})|v_{n}+\widetilde{w}_{\Omega}(v_{n})|^{2^{\ast}_{\alpha}-2}(v_{n}+z_{n})\\
&~~~~~~~~~~~~~~~~~~-(I_{\alpha}\ast|v_{0}+w_{0}+z_{0}|^{2^{\ast}_{\alpha}})|v_{0}+w_{0}+z_{0}|^{2^{\ast}_{\alpha}-2}(v_{0}+z_{0}),v_{n}-v_{0}           \Big\rangle dx\\
&\leq C\int_{\Omega}|\varphi|^{2\cdot2^{\ast}_{\alpha}}\Big\langle (I_{\alpha}\ast|v_{n}|^{2^{\ast}_{\alpha}})|v_{n}|^{2^{\ast}_{\alpha}-2}v_{n}-(I_{\alpha}\ast|v_{0}|^{2^{\ast}_{\alpha}})|v_{0}|^{2^{\ast}_{\alpha}-2}v_{0},v_{n}-v_{0}           \Big\rangle dx\\
&\leq C\left(\int_{\Omega}|\varphi|^{2\cdot2^{\ast}_{\alpha}}(I_{\alpha}\ast|v_{n}-v_{0}|^{2^{\ast}_{\alpha}})|v_{n}-v_{0}|^{2^{\ast}_{\alpha}}dx\right)^{\frac{1}{2^{\ast}_{\alpha}}} =C\left(\int_{\Omega}|\varphi|^{2\cdot2^{\ast}_{\alpha}}d\overline{\rho}\right)^{\frac{1}{2^{\ast}_{\alpha}}}.
\end{aligned}
\end{equation}
where we use the fact that $Z$ is a finite dimension space and $\int_{\Omega}\langle v,w\rangle dx=0$ see (\ref{orth}).

Next, we are going to show that $I_{2}=o(1)$. Fix $l\geq 1$. In view of Lemma \ref{dense}, Lemma \ref{poincare} and Lemma 1.1 in \cite{Leinfelder1983}, there exists $\xi_{n}\in W^{1,\alpha,2^{\ast}_{\alpha}}(B_{l})$ such that $w_{n}=\nabla\xi_{n}$ and we may assume without loss of generality that $\int_{B_{l}}\xi_{n}dx=0$. Then by the Poincar\'{e} inequality in Lemma \ref{poincare}
\begin{equation*}
||\xi_{n}||_{W^{1,\alpha,2^{\ast}_{\alpha}}(B_{l})}\leq  C|w_{n}|_{Q^{\alpha,2^{\ast}_{\alpha}}(B_{l},\mathbb{R}^{3})}\leq C|w_{n}|_{Q^{\alpha,2^{\ast}_{\alpha}}(\mathbb{R}^{3},\mathbb{R}^{3})},
\end{equation*}
and passing to a subsequence, $\xi_{n}\rightharpoonup\xi$ for some $\xi\in W^{1,\alpha,2^{\ast}_{\alpha}}(B_{l})$. So by the natural compactly embedding, $\xi_{n}\longrightarrow\xi$ in $Q^{\alpha,2^{\ast}_{\alpha}}(B_{l})$. Now take any $\varphi\in C^{\infty}_{0}(B_{l})$. Since $\nabla(|\varphi|^{2\cdot2^{\ast}_{\alpha}}(\xi_{n}-\xi))\in \widetilde{\mathcal{W}}$, in view of (\ref{mini2}) we get
\begin{equation*}
\int_{\Omega}\langle f(v_{n}+\widetilde{w}_{\Omega}(v_{n})),\nabla(|\varphi|^{2\cdot2^{\ast}_{\alpha}}(\xi_{n}-\xi))\rangle dx=0.
\end{equation*}
That is
\begin{equation*}
\begin{aligned}
\int_{\Omega}|\varphi|^{2\cdot2^{\ast}_{\alpha}}\langle &f(v_{n}+\widetilde{w}_{\Omega}(v_{n})),w_{n}-\nabla\xi\rangle dx=\int_{\Omega}\langle f(v_{n}+\widetilde{w}_{\Omega}(v_{n})),\nabla(|\varphi|^{2\cdot2^{\ast}_{\alpha}}) (\xi-\xi_{n})\rangle dx
\end{aligned}
\end{equation*}
where the right-hand side tends to 0 as $n\longrightarrow \infty$. Since $w_{n}\rightharpoonup\nabla\xi$ in $Q^{\alpha,2^{\ast}_{\alpha}}(B_{l})$,
\begin{equation*}
\int_{\Omega}|\varphi|^{2\cdot2^{\ast}_{\alpha}}\langle f(v_{0}+\nabla\xi+z_{0}),w_{n}-\nabla\xi\rangle dx=o(1).
\end{equation*}
Hence, recalling that $\widetilde{w}_{\Omega}(v_{n})=w_{n}+z_{n}$ and $z_{n}\longrightarrow z_{0}$, we obtain
\begin{equation}\label{o(1)}
I_{2}=\int_{\Omega}|\varphi|^{2\cdot2^{\ast}_{\alpha}}\langle f(v_{n}+\widetilde{w}_{\Omega}(v_{n}))-f(v_{0}+\nabla\xi+z_{0}),\widetilde{w}_{\Omega}(v_{n})-\nabla\xi-z_{0}\rangle dx=o(1).
\end{equation}

%Here we have used the fact that $\int_{\Omega}a(x)|v_{n}-v_{0}|^{2}dx\longrightarrow 0$ if %$v_{n}\rightharpoonup v_{0}$ in $Q^{\alpha,2^{\ast}_{\alpha}}(\Omega,\mathbb{R}^{3})$.
Since $\varphi\in C^{\infty}_{0}(B_{l})$ is arbitrary, it follows from (\ref{measure1}) and (\ref{o(1)}) that
\begin{equation}\label{measure}
4m_{k}|\Omega_{n,k}\cap E|\leq(\overline{\rho}(E))^{1/2^{\ast}_{\alpha}}+o(1)
\end{equation}
for any Borel set $E\subset B_{l}$. On the other hand, we can find an open set $E_{k}\supset I$ such that $|E_{k}|<\frac{1}{2^{k+1}}$. Then, taking $E=B_{l}\setminus E_{k}$ in (\ref{measure}), we have $4m_{k}|\Omega_{n,k}\cap(B_{l}\setminus E_{k})|=o(1)$ as $n\longrightarrow \infty$ because $\mathrm{supp}(\overline{\rho})\subset I$; hence we can find a sufficiently large $n_{k}$ such that $|\Omega_{n_{k},k}\cap B_{l}|<\frac{1}{2^{k}}$ and we obtain
\begin{equation*}
|\mathop{\cap}\limits_{j=1}^{\infty}\mathop{\cup}\limits_{k=j}^{\infty}\Omega_{n_{k},k}\cap B_{l}|\leq \mathop{\mathrm{lim}}\limits_{j\longrightarrow\infty}\mathop{\Sigma}\limits_{k=j}^{\infty}|\Omega_{n_{k},k}\cap B_{l}|\leq\mathop{\mathrm{lim}}\limits_{j\longrightarrow\infty}\frac{1}{2^{j-1}}=0.
\end{equation*}
According to the fact that $\widetilde{w}_{\Omega}(v_{n})\rightharpoonup\widetilde{w}_{0}$, one can employ the diagonal procedure and hence find a subsequence of $\widetilde{w}_{\Omega}(v_{n})$ which converges to $\widetilde{w}_{0}$ a.e. in $\Omega=\cup_{l=1}^{\infty}B_{l}$.

Let $p\in[1,2^{\ast}_{\alpha}]$. For $\Omega'\subset\Omega$ such that $|\Omega'|<+\infty$ we have
\begin{equation*}
\begin{aligned}
&\int_{\Omega'}|v_{n}-v_{0}+\widetilde{w}_{\Omega}(v_{n})-\widetilde{w}_{0}|^{p}dx
\leq |\Omega'|^{1-\frac{p}{2^{\ast}_{\alpha}}} \left(\int_{\Omega}|v_{n}-v_{0}+\widetilde{w}_{\Omega}(v_{n})-\widetilde{w}_{0}|^{2^{\ast}_{\alpha}}dx\right)^{\frac{p}{2^{\ast}_{\alpha}}}\\
&\leq |\Omega'|^{1-\frac{p}{2^{\ast}_{\alpha}}}|\mathrm{diam}\Omega|^{\frac{3-\alpha}{2}\cdot\frac{p}{2^{\ast}_{\alpha}}}\left(\int_{\Omega}|I_{\alpha/2}\ast|v_{n}-v_{0}+\widetilde{w}_{\Omega}(v_{n})-\widetilde{w}_{0}|^{2^{\ast}_{\alpha}}|^{2}\right)^{\frac{1}{2}\cdot\frac{p}{2^{\ast}_{\alpha}}},
\end{aligned}
\end{equation*}
where $\mathrm{diam}\Omega=\mathop{max}\limits_{x,y\in\Omega}|x-y|$. Hence by the Vitali convergence theorem, $v_{n}-v_{0}+\widetilde{w}_{\Omega}(v_{n})-\widetilde{w}_{0}\longrightarrow0 $ in $L^{p}_{loc}(\Omega)$ after passing to a subsequence.

Step 3. We show that $\widetilde{w}_{\Omega}(v_{0})=\widetilde{w}_{0}$. Take any $\widetilde{w}\in\widetilde{\mathcal{W}}$ and observe that by the Vitali convergence theorem,
\begin{equation*}
0=\int_{\Omega}\langle f(v_{n}+\widetilde{w}_{\Omega}(v_{n})),\widetilde{w}\rangle dx\longrightarrow \int_{\Omega}\langle f(v_{0}+\widetilde{w}_{0}),\widetilde{w}\rangle dx,
\end{equation*}
up to a subsequence. Now (\ref{mini2}) implies that $\widetilde{w}_{0} =\widetilde{w}_{\Omega}(v_{0})$ which completes the proof.
\end{proof}

\subsection{Abstract Critical Point Theory}

For readers convenient, we end this section with recalling the abstract critical point lemma, see \cite[Section 4]{Bartsch2015-1} and \cite[Section 3]{Mederski2018} for more details. Let $X$ be a reflexive Banach space with norm $||\cdot||$ and with a topological direct sum decomposition $X=X^{+}\oplus\widetilde{X}$, where $X^{+}$ is a Hilbert space with a scalar product. For $u\in X$ we denote by $u^{+}\in X^{+}$ and $\widetilde{u}\in\widetilde{X}$ the corresponding summands so that $u=u^{+}+\widetilde{u}$. We may assume that $\langle u,v\rangle=||u||^{2}$ for any $u\in X^{+}$ and that $||u||^{2}=||u^{+}||^{2}+||\widetilde{u}||^{2}$. The topology $\mathcal{T}$ on $X$ is defined as the product of the norm topology in $X^{+}$ and the weak topology in $\widetilde{X}$. Thus $u_{n}\mathop{\longrightarrow}\limits^{\mathcal{T}}u$ is equivalent to $u^{+}_{n}\longrightarrow u^{+}$ and $\widetilde{u}\rightharpoonup\widetilde{u}$.

Let $J\in \mathcal{C}^{1}(X,\mathbb{R})$ be a functional on $X$ of the form
\begin{equation*}
J(u)=\frac{1}{2}||u^{+}||-I(u)~~\mathrm{for}~~u=u^{+}+\widetilde{u}\in X^{+}\oplus\widetilde{X}
\end{equation*}
such that the following assumptions hold

\noindent (A1)\quad $I\in \mathcal{C}^{1}(X,\mathbb{R})$ and $I(u)\geq I(0)=0$ for any $u\in X$.

\noindent (A2)\quad $I$ is $\mathcal{T}$-sequentially lower semi-continuous: $u_{n}\mathop{\longrightarrow}\limits^{\mathcal{T}}u\Longrightarrow \mathrm{lim}~\mathrm{inf}~I(u_{n})\geq I(u)$.

\noindent (A3)\quad If $u_{n}\mathop{\longrightarrow}\limits^{\mathcal{T}} u$ and $I(u_{n})\longrightarrow I(u)$ then $u_{n}\longrightarrow u$.

\noindent (A4)\quad There exists $r>0$ such that $a:=\mathop{\mathrm{inf}}\limits_{u\in X^{+}:||u||=r}J(u)>0$.

\noindent (B1)\quad $||u^{+}||+I(u)\longrightarrow\infty$ as $||u||\longrightarrow \infty$.

\noindent (B2)\quad $I(t_{n}u_{n})/t_{n}^{2}\longrightarrow\infty$ if $t_{n}\longrightarrow \infty$ and $u^{+}_{n}\longrightarrow u^{+}$ for some $u^{+}\neq 0$ as $n\longrightarrow\infty$.

\noindent (B3)\quad $\frac{t^{2}-1}{2}I'(u)(u)+tI'(u)(v)+I(u)-I(tu+v)<0$ for every $u\in\mathcal{N},t>0,v\in X$ such that $u\neq tu+v$.

We defined the following Nehari-Pankov
\begin{equation*}\label{define}
\mathcal{N}:=\{u\in X\setminus\widetilde{X}:J'(u)|_{\mathbb{R}u\oplus \widetilde{X}}=0\}.
\end{equation*}
Correspondingly, we defined the (PS)$^{\mathcal{T}}_{c}$ condition for $J$.
\begin{Def}\label{PSC}
We say that $J$ satisfies the $(PS)_{c}^{\mathcal{T}}$ condition in $\mathcal{N}$ if every $(PS)_{c}$ sequence $(u_{n})\in \mathcal{N}$ has a subsequence which convergence in the $\mathcal{T}$ topology:
\begin{equation*}
u_{n}\in\mathcal{N},~~ J(u_{n})\longrightarrow 0,~~ J'(u_{n})\longrightarrow c~~~~\Longrightarrow~~~~u_{n}\mathop{\longrightarrow}\limits^{\mathcal{T}}u\in X~~along~a~subsequence.
\end{equation*}
\end{Def}

We also recall the compactly perturbed problem with respect to another decomposition of $X$. Namely, suppose that
\begin{equation}\label{decomposition 2}
\widetilde{X}=X^{0}\oplus X^{1},
\end{equation}
where $X^{0},X^{1}$ are closed in $\widetilde{X}$, and $X^{0}$ is a Hilbert space. For $u\in \widetilde{X}$ we denote $u^{0}\in X^{0}$ and $u^{1}\in X^{1}$ the corresponding summands so that $u=u^{0}+u^{1}$. We use the same notation for the scalar product in $X^{+}\oplus X^{0}$ and $\langle u,u\rangle=||u||^{2}=||u^{+}||^{2}+||u^{0}||^{2}$ for any $u=u^{+}+u^{0}\in X^{+}\oplus X^{0}$, hence $X^{+}$ and $X^{0}$ are orthogonal. We consider another functional $J_{cp}\in\mathcal{C}^{1}(X,\mathbb{R})$ of the form
\begin{equation*}
J_{cp}=\frac{1}{2}||u^{+}+u^{0}||^{2}-I_{cp}(u)~~\mathrm{for}~u=u^{+}+u^{0}+u^{1}\in X^{+}\oplus X^{0}\oplus X^{1}.
\end{equation*}
We define the corresponding Nehari-Pankov manifold for $J_{cp}$
\begin{equation*}
\mathcal{N}_{cp}:=\{u\in X\setminus X^{1}:J'_{cp}(u)|_{\mathbb{R}u\oplus X^{1}}=0\},
\end{equation*}
and assume that $J_{cp}$ satisfies all corresponding assumption (A1)-(A4), (B1)-(B3), where we replace $X^{+}\oplus X^{0}$, $X^{1}$ and $I_{cp}$ instead of $X^{+}, X$ and $I$ respectively. Moreover, we enlist new additional conditions:

(C1)\quad $J_{cp}(u_{n})-J_{u_{n}}\longrightarrow 0$ if $(u_{n})\subset \mathcal{N}_{cp}$ is bounded and $(u^{+}_{n}+u^{0}_{n})\rightharpoonup 0$. Moreover there is $M>0$ such that $J_{cp}(u)-J(u)\leq M||u^{+}+u^{0}||^{2}$ for $u\in \mathcal{N}_{cp}$.

(C2)\quad $I(t_{n}u_{n})\setminus t^{2}_{n}\longrightarrow \infty$ and $t_{n}\longrightarrow\infty$ and $(I(tu^{+}_{n}))_{n}$ is bounded away from 0 for any $t>1$.

(C3)\quad $J'$ is weak-to-weak$^{\ast}$ continuous on $\mathcal{N}$, i.e. if $(u_{n})_{n}\subset \mathcal{N}$, $u_{n}\rightharpoonup u$, then $J'(u_{n})\mathop{\rightharpoonup}\limits^{\ast} J'(u)$ in $X^{\ast}$. Moreover $J$ is weakly sequentially lower semi-continuous on $\mathcal{N}$, i.e. if $(u_{n})_{n}\subset \mathcal{N}$, $u_{n}\rightharpoonup u$ and $u\in\mathcal{N}$, then $\mathop{\mathrm{lim}~\mathrm{inf}}\limits_{n\longrightarrow \infty}J(u_{n})\geq J(u)$.

There we present the abstract critical point theorem:
\begin{lem}\cite[Theorem 3.2]{Mederski 2018}\label{abstract}:
 Let $J\in \mathcal{C}^{1}(X,\mathbb{R})$ be coercive on $\mathcal{N}$ and let $J_{cp}\in\mathcal{C}^{1}(X,\mathbb{R})$ be coercive on $\mathcal{N}_{cp}$. Suppose that $J$ and $J_{cp}$ satisfy (A1)-(A4),(B1)-(B3) and set $c=\mathop{\mathrm{inf}}\limits_{\mathcal{N}}J$ and $d=\mathop{\mathrm{inf}}\limits_{\mathcal{N}_{cp}}J_{cp}$. Then the following statements hold:

 (a)\quad If (C1)-(C2) hold and $\beta<d$, then any $(PS)_{\beta}$-sequence in $\mathcal{N}$ contains a weakly convergent subsequence with a nontrivial limit point.

 (b)\quad If (C1)-(C3) hold and $c<d$, then $c$ is achieved by a critical point (ground state) of $J$.

 (c)\quad Suppose that $J$ is even and satisfies the $(PS)_{\beta}^{\mathcal{T}}$-condition in $\mathcal{N}$ for any $\beta<\beta_{0}$ for some fixed $\beta_{0}\in(c,\infty]$. Let
 \begin{equation*}
m(\mathcal{N},\beta_{0})=sup\{\gamma(J^{-1}((0,\beta)\cap\mathcal{N}):\beta<\beta_{0}\}\in \mathbb{N}_{0},
 \end{equation*}
 where $\gamma$ stands for the Krasnoselskii genus for closed and symmetric subsets of $X$. Then $J$ has at least $m(\mathcal{N},\beta_{0})$ pairs of critical points $u$ and $-u$ such that $u\neq0$ and $c\leq J(u)<\beta_{0}$.
\end{lem}

\section{\bf  Sharp constant $S_{\mathrm{curl},HL}(\mathbb{R}^{3})$}
\subsection{Proof of Theorem \ref{attained}}

In this subsection, we consider functional (\ref{J(u)}), which is associated to equation (\ref{limit equation}), and we work on the following Nehari-Pankov manifold
\begin{equation*}\label{manifold2}
\mathcal{N}:=\left\{u\in W^{\alpha,2^{\ast}_{\alpha}}_{0}(\mathrm{curl};\mathbb{R}^{3})\setminus \mathcal{W}_{\mathbb{R}^{3}}: J'(u)u=0~\mathrm{and}~J'(u)|_{\mathcal{W}_{\mathbb{R}^{3}}}=0\right\}.
\end{equation*}
\begin{lem}
There exists a continuous mapping $m:\mathcal{V}_{\mathbb{R}^{3}}\setminus\{0\}\longrightarrow\mathcal{N}$.
\end{lem}
\begin{proof}
By Lemma \ref{Helmholtz}, $W_{0}^{\alpha,2^{\ast}_{\alpha}}(\mathrm{curl};\mathbb{R}^{3})=\mathcal{V}_{\mathbb{R}^{3}}\oplus\mathcal{W}_{\mathbb{R}^{3}}$. It follows from (\ref{mini1}) and (\ref{mini2}) that if $v\in\mathcal{V}_{\mathbb{R}^{3}}$ and $\widetilde{w}_{\mathbb{R}^{3}}(v)\in\widetilde{\mathcal{W}}=\mathcal{W}_{\mathbb{R}^{3}}$, then we have $J'(v+\widetilde{w}_{\mathbb{R}^{3}}(v))|_{\mathcal{W}_{\mathbb{R}^{3}}}=0$. And as
\begin{equation}\label{t}
J(t(v+\widetilde{w}(v)))=\frac{t^{2}}{2}\int_{\mathbb{R}^{3}}|\nabla v|^{2}dx-\frac{t^{2\cdot2^{\ast}_{\alpha} }}{2\cdot 2^{\ast}_{\alpha}}\int_{\mathbb{R}^{3}}|I_{\alpha/2}\ast|v+\widetilde{w}_{\mathbb{R}^{3}}(v)|^{2^{\ast}_{\alpha}}|^{2}dx,
\end{equation}
there is a unique $t(v)>0$ such that
\begin{equation}\label{m}
t(v)(v+\widetilde{w}_{\mathbb{R}^{3}}(v))\in\mathcal{N}~~\mathrm{for}~v\in\mathcal{V}_{\mathbb{R}^{3}}\setminus\{0\}.
\end{equation}
Setting $m(v):=t(v)(v+\widetilde{w}_{\mathbb{R}^{3}}(v))$, we then note that
\begin{equation}\label{note}
J(m(v))\geq J(t(v+\widetilde{w}_{\mathbb{R}^{3}}))~~\mathrm{for}~\mathrm{all}~t>0~\mathrm{and}~\widetilde{w}_{\mathbb{R}^{3}}\in\mathcal{W}_{\mathbb{R}^{3}}.
\end{equation}
Since $J(m(v))\geq J(v)$ and there exist $a,r>0$ such that $J(v)\geq a$ if $||v||=r$, this implies that $\mathcal{N}$ is bounded away from $\mathcal{W}_{\mathbb{R}^{3}}$ and hence closed. Therefore,  by the similar analysis in \cite[Lemma 4.4 ]{Mederski2021}, the mapping $m$ is continuous.
\end{proof}
\begin{lem}
Set $\mathcal{S}:=\{v\in\mathcal{V}_{\mathbb{R}^{3}}:||v||=1\}$, there exist a $(PS)_{c}$ sequence $(v_{n})$ for $J\circ m$, and a $(PS)_{c}$ sequence $(m(v_{n}))$ for $J$ on $\mathcal{N}$.
\end{lem}
\begin{proof}
By the continuity of mapping $m$, we easily observe that $m|_{\mathcal{S}}:\mathcal{S}\longrightarrow\mathcal{N}$ is a homeomorphism with the inverse $u=v+m(v)\mapsto\frac{v}{||v||}$. Recall the argument in \cite[Proposition 4.4(b)]{Mederski2015}, we know that $J\circ m|_{\mathcal{S}}:\mathcal{S}\longrightarrow
\mathcal{R}$ is of class $\mathcal{C}^{1}$ and is bounded from below by the constant $a>0$. By the Ekeland variational principle, there is a $(PS)_{c}$ sequence $(v_{n})\subset \mathcal{S}$ such that
\begin{equation}\label{minimizer}
(J\circ m)(v_{n})\longrightarrow \mathop{\mathrm{inf}}\limits_{\mathcal{S}}J\circ m=\mathop{\mathrm{inf}}\limits_{\mathcal{N}}J\geq a>0.
\end{equation}
Again, by the argument in \cite[Proposition 4.4(b)]{Mederski2015}, we have $(m(v_{n}))$ is a $(PS)_{c}$ sequence for $J$ on $\mathcal{N}$.
\end{proof}

\begin{proof}[\bf Complete of the Proof of Theorem \ref{attained}]
Firstly, we prove part (a). Taking a minimizing sequence $(u_{n})=(m(v_{n}))\subset \mathcal{N}$ and set $u_{n}=t(v_{n})(v_{n}+\widetilde{w}_{\mathbb{R}^{3}}(v_{n}))=v'_{n}+\widetilde{w}_{\mathbb{R}^{3}}(v'_{n})\in\mathcal{V}_{\mathbb{R}^{3}}\oplus\mathcal{W}_{\mathbb{R}^{3}}$. Then we have
\begin{equation*}
J(u_{n})=J(u_{n})-\frac{1}{2\cdot 2^{\ast}_{\alpha}}J'(u_{n})u_{n}=\frac{2^{\ast}_{\alpha}-1}{2\cdot2^{\ast}_{\alpha}}|\nabla\times u_{n}|^{2}_{2}=\frac{2^{\ast}_{\alpha}-1}{2\cdot2^{\ast}_{\alpha}}|\nabla v'_{n}|^{2}_{2}.
\end{equation*}
Since the norm $|\nabla\cdot|_{2}$ is an equivalent norm in $\mathcal{V}_{\mathbb{R}^{3}}$, it follows that $J(u_{n})$ is coercive on $\mathcal{N}$, hence $(v'_{n})$ is bounded. On the other hand, we also have
\begin{equation*}
J(u_{n})=J(u_{n})-\frac{1}{2}J'(u_{n})u_{n}=\frac{2^{\ast}_{\alpha}-1}{2\cdot2^{\ast}_{\alpha}}\int_{\mathbb{R}^{3}}|I_{\alpha/2}\ast|u_{n}|^{2^{\ast}_{\alpha}}|^{2}dx.
\end{equation*}
By (\ref{minimizer}), $J(u_{n})$ is bounded away from $0$, so is $|u_{n}|_{Q^{\alpha,2^{\ast}_{\alpha}}}\nrightarrow 0$, and hence by (\ref{mini3}), we also have $|v'_{n}|_{Q^{\alpha,2^{\ast}_{\alpha}}}\nrightarrow 0$.

Denote $T_{s,y}(v'):=s^{1/2}v'(s\cdot+y)$, where $s>0$, $y\in\mathbb{R}^{3}$.  Then, passing to a subsequence and using the argument in \cite[Theorem 1]{Solimini1995}, we have $\overline{v}_{n}=T_{s_{n},y_{n}}(v'_{n})\rightharpoonup v_{0}$ for some $v_{0}\neq0$, where $(s_{n})\subset \mathbb{R}^{+}$ and $(y_{n})\subset\mathbb{R}^{3}$. Taking subsequence again, we also have $\overline{v}_{n}\longrightarrow v_{0}$ a.e. in $\mathbb{R}^{3}$ and in view of the concentration-compactness Lemma \ref{concentration}, we deduce $\widetilde{w}_{\mathbb{R}^{3}}(\overline{v}_{n})\rightharpoonup \widetilde{w}_{\mathbb{R}^{3}}(v_{0})$ and $\widetilde{w}_{\mathbb{R}^{3}}(\overline{v}_{n})\longrightarrow \widetilde{w}_{\mathbb{R}^{3}}(v_{0})$ a.e. in $\mathbb{R}^{3}$. Setting $u:=v_{0}+\widetilde{w}_{\mathbb{R}^{3}}(v_{0})$ and assume without loss of generality that $s_{n}=1$ and $y_{n}=0$, then by Lemma 4.6 in \cite{Mederski2021}, we have $u_{n}\rightharpoonup u$ and $u_{n}\longrightarrow u$ a.e. in $\mathbb{R}^{3}$. Moreover, by Lemma \ref{well define}
%the HLS inequality, the Riesz potential defines a linear continuous map from %$L^{\frac{2N}{2N-\mu}}(\mathbb{R}^{3},\mathbb{R}^{3})$ to %$L^{\frac{2N}{\mu}}(\mathbb{R}^{3},\mathbb{R}^{3})$, see \cite[pp.1226]{Gao2018}.
we have
\begin{equation*}\label{weak3}
(I_{\alpha}\ast|u_{n}|^{2^{\ast}_{\alpha}})|u_{n}|^{2^{\ast}_{\alpha}-2}u_{n}
\rightharpoonup(I_{\alpha}\ast|u|^{2^{\ast}_{\alpha}})|u|^{2^{\ast}_{\alpha}-2}u~~\mathrm{in}~~(Q^{\alpha,2^{\ast}_{\alpha}}(\mathbb{R}^{3},\mathbb{R}^{3}))',
\end{equation*}
Therefore, for any $z\in W^{\alpha,2^{\ast}_{\alpha}}_{0}(\mathrm{curl};\mathbb{R}^{3})$, using weak and a.e. convergence, we have
\begin{equation*}\label{weak}
\langle J'(u_{n}),z\rangle=\int_{\mathbb{R}^{3}}\langle\nabla\times u_{n},z\rangle dx-\int_{\mathbb{R}^{3}}\left\langle\left(I_{\alpha}\ast|u_{n}|^{2^{\ast}_{\alpha}}\right)|u_{n}(x)|^{2^{\ast}_{\alpha}-2}u_{n}(x),z\right\rangle dx
\longrightarrow \langle J'(u),z\rangle.
\end{equation*}
This implies that $u$ is a solution to (\ref{limit equation}). Using Fatou's lemma, we deduce that
\begin{equation*}
\begin{aligned}
\mathop{\mathrm{inf}}\limits_{\mathcal{N}}J&=J(u_{n})+o(1)=J(u_{n})-\frac{1}{2}J'(u_{n})u_{n}+o(1)\\
&=\frac{2^{\ast}_{\alpha}-1}{2\cdot2^{\ast}_{\alpha}}
\int_{\mathbb{R}^{3}}|I_{\alpha/2}\ast|u_{n}|^{2^{\ast}_{\alpha}}|^{2}dx+o(1)\geq \frac{2^{\ast}_{\alpha}-1}{2\cdot2^{\ast}_{\alpha}}\int_{\mathbb{R}^{3}}|I_{\alpha/2}\ast|u|^{2^{\ast}_{\alpha}}|^{2}dx
+o(1)\\
&=J(u)-\frac{1}{2}J'(u)u+o(1)=J(u)+o(1).
\end{aligned}
\end{equation*}
Hence $J(u)\leq \mathop{\mathrm{inf}}\limits_{\mathcal{N}} J\leq J(u)$ and as a solution, $u\in\mathcal{N}$.

Next, we show $\mathop{\mathrm{inf}}\limits_{\mathcal{N}} J=\frac{2^{\ast}_{\alpha}-1}{2\cdot2^{\ast}_{\alpha}} S_{\mathrm{curl},HL}^{\frac{2^{\ast}_{\alpha}}{2^{\ast}_{\alpha}-1}}$, where $S_{\mathrm{curl},HL}$ is the sharp constant in (\ref{inequality1}), which can be rewritten as follow
\begin{equation}\label{ScurlHL}
S_{\mathrm{curl},HL}=\mathop{\mathrm{inf}}\limits_{\mathop{u\in W^{\alpha,2^{\ast}_{\alpha}}_{0}(\mathrm{curl};\mathbb{R}^{3})}\limits_{\nabla\times u\neq0}}\frac{\int_{\mathbb{R}^{3}}|\nabla\times u|^{2}dx}{\left(\int_{\mathbb{R}^{3}}|I_{\alpha/2}\ast|u+\widetilde{w}_{\mathbb{R}^{3}}(u)|^{2^{\ast}_{\alpha}}|^{2}dx\right)^{\frac{1}{2^{\ast}_{\alpha}}}}.
\end{equation}
In fact, by (\ref{mini1}), it is clear that a minimize $\widetilde{w}_{\mathbb{R}^{3}}(u)$ exists uniquely for any $u\in W^{\alpha,2^{\ast}_{\alpha}}_{0}(\mathrm{curl};\Omega)$, not only $u\in\mathcal{V}_{\mathbb{R}^{3}}$. So by Lemma \ref{Helmholtz}, $u+\widetilde{w}_{\mathbb{R}^{3}}(u)=v+\widetilde{w}_{\mathbb{R}^{3}}(v)\in \mathcal{V}_{\mathbb{R}^{3}}\oplus\mathcal{W}_{\mathbb{R}^{3}}$ for some $v\in\mathcal{V}_{\mathbb{R}^{3}}$ and therefore
\begin{equation}\label{equlity}
\begin{aligned}
\mathop{\mathrm{inf}}\limits_{w\in\mathcal{W}_{\mathbb{R}^{3}}}\int_{\mathbb{R}^{3}}|I_{\alpha/2}\ast|u+w|^{2^{\ast}_{\alpha}}|^{2}dx
&=\int_{\mathbb{R}^{3}}|I_{\alpha/2}\ast|u+\widetilde{w}_{\mathbb{R}^{3}}(u)|^{2^{\ast}_{\alpha}}|^{2}dx=\int_{\mathbb{R}^{3}}|I_{\alpha/2}\ast|v+\widetilde{w}_{\mathbb{R}^{3}}(v)|^{2^{\ast}_{\alpha}}|^{2}dx.
\end{aligned}
\end{equation}
On the other hand, since $u+\widetilde{w}_{\mathbb{R}^{3}}(u)\in\mathcal{N}$, $J'(u)u=0$, i.e.
\begin{equation*}
\int_{\mathbb{R}^{3}}|\nabla\times u|^{2}dx=\int_{\mathbb{R}^{3}}|I_{\alpha/2}\ast|u+\widetilde{w}_{\mathbb{R}^{3}}(u)|^{2^{\ast}_{\alpha}}|^{2}dx.
\end{equation*}
Then we can easily calculate that
\begin{equation*}
\mathop{\mathrm{inf}}\limits_{\mathcal{N}} J=\frac{2^{\ast}_{\alpha}-1}{2\cdot2^{\ast}_{\alpha}}
\int_{\mathbb{R}^{3}}|\nabla\times u|^{2}dx
=\frac{2^{\ast}_{\alpha}-1}{2\cdot2^{\ast}_{\alpha}} S_{\mathrm{curl},HL}^{\frac{2^{\ast}_{\alpha}}{2^{\ast}_{\alpha}-1}}.
\end{equation*}

As we can see, if $u$ satisfies equality (\ref{inequality1}), then $t(u)(u+\widetilde{w}_{\mathbb{R}^{3}}(u))\in\mathcal{N}$ and is a minimizer for $J|_{\mathcal{N}}$ and the corresponding point $v$ in $\mathcal{S}$ is a minimizer for $J\circ m|_{\mathcal{S}}$, see (\ref{minimizer}). Hence $v$ is a critical point of $J\circ m|_{\mathcal{S}}$ and $m(v)=u$ is a critical point of $J$. This completes the proof of (a).

(b) To compare the constants $S_{\mathrm{curl},HL}$ and $S_{HL}$, see  (\ref{ScurlHL}) and (\ref{SHL}), we firstly claim that $S_{\mathrm{curl},HL}\geq S_{HL}$. In fact, by (\ref{equlity}) and $\mathrm{div}v=0$, we have
\begin{equation}\label{ScurlHL1}
S_{\mathrm{curl},HL}=\mathop{\mathrm{inf}}\limits_{v\in\mathcal{V}_{\mathbb{R}^{3}}\setminus\{0\}}\frac{\int_{\mathbb{R}^{3}}|\nabla v|^{2}dx}{\left(\int_{\mathbb{R}^{3}}|I_{\alpha/2}\ast|v+\widetilde{w}_{\mathbb{R}^{3}}(v)|^{2^{\ast}_{\alpha}}|^{2}dx\right)^{\frac{1}{2^{\ast}_{\alpha}}}}.
\end{equation}
Then given $\varepsilon>0$, we can find $v\neq0$ such that
\begin{equation*}
\int_{\mathbb{R}^{3}}|\nabla v|^{2}dx\leq (S_{\mathrm{curl},HL}+\varepsilon)\left(\int_{\mathbb{R}^{3}}|I_{\alpha/2}\ast|v+\widetilde{w}_{\mathbb{R}^{3}}(v)|^{2^{\ast}_{\alpha}}|^{2}dx\right)^{\frac{1}{2^{\ast}_{\alpha}}}.
\end{equation*}
Since $\widetilde{w}_{\mathbb{R}^{3}}(v)$ is a minimizer, we deduce that
\begin{equation*}
\int_{\mathbb{R}^{3}}|\nabla v|^{2}dx\leq (S_{\mathrm{curl},HL}+\varepsilon)\left(\int_{\mathbb{R}^{3}}|I_{\alpha/2}\ast|v|^{2^{\ast}_{\alpha}}|^{2}dx\right)^{\frac{1}{2^{\ast}_{\alpha}}}.
\end{equation*}
Then by the definition of $S_{HL}$, see (\ref{SHL}), we have
\begin{equation*}
\int_{\mathbb{R}^{3}}|\nabla v|^{2}dx\leq\frac{(S_{\mathrm{curl},HL}+\varepsilon)}{S_{HL}}\int_{\mathbb{R}^{3}}|\nabla v|^{2}dx.
\end{equation*}
Hence we get our claim by $S_{\mathrm{curl},HL}+\varepsilon\geq S_{HL}$.

Secondly, we exclude the case $S_{\mathrm{curl},HL}=S_{HL}$. Otherwise, all inequalities above become equalities with $\varepsilon=0$. Then by the form of the optimal function to $S_{HL}$, see (\ref{form}), up to multiplicative constants, we get the contradiction with $\mathrm{div}v=0$. These show that $S_{\mathrm{curl},HL}>S_{HL}$.
\end{proof}

\subsection{Proof of Theorem \ref{four constant}}

To compare the sharp constants $S_{\mathrm{curl},HL}(\mathbb{R}^{3})$ and $\bar{S}_{\mathrm{curl},HL}(\Omega)$, we have introduced another constant $S_{\mathrm{curl},HL}(\Omega)$. Recall from Section 2.1 that we have the following Helmholtz decomposition in entire space $\mathbb{R}^{3}$ and in the bounded domain $\Omega$:
\begin{equation*}
W^{\alpha,2^{\ast}_{\alpha}}_{0}(\mathrm{curl};\mathbb{R}^{3})=\mathcal{V}_{\mathbb{R}^{3}}\oplus\mathcal{W}_{\mathbb{R}^{3}}
~~\mathrm{and}~~W^{\alpha,2^{\ast}_{\alpha}}_{0}(\mathrm{curl};\Omega)=\mathcal{V}_{\Omega}\oplus\mathcal{W}_{\Omega}.
\end{equation*}
%Then for $u\in W^{\alpha,2^{\ast}_{\alpha}}_{0}(\mathrm{curl};\Omega)$, we denote %$\widetilde{w}_{\Omega}(u)$ by the minimizer of
%\begin{equation*}
%\int_{\Omega}|I_{\alpha/2}\ast|u+w|^{2^{\ast}_{\alpha}}|^{2}dx,~~w\in\mathcal{W}_{\Omega},
%\end{equation*}
%see (\ref{equlity}).
Then, as (\ref{ScurlHL}), we note that $S_{\mathrm{curl},HL}(\Omega)$ (see (\ref{inequality2})) can be characterized as
\begin{equation}\label{SCURLHL}
\begin{aligned}
S_{\mathrm{curl},HL}(\Omega)&=\mathop{\mathrm{inf}}\limits_{\mathop{u\in W^{\alpha,2^{\ast}_{\alpha}}_{0}(\mathrm{curl};\Omega)}\limits_{\nabla\times u\neq0}}\mathop{\mathrm{sup}}\limits_{w\in\mathcal{W}_{\mathbb{R}^{3}}}\frac{\int_{\mathbb{R}^{3}}|\nabla\times u|^{2}dx}{\left(\int_{\mathbb{R}^{3}}|I_{\alpha/2}\ast|u+w|^{2^{\ast}_{\alpha}}|^{2}dx\right)^{\frac{1}{2^{\ast}_{\alpha}}}}\\
&=\mathop{\mathrm{inf}}\limits_{\mathop{u\in W^{\alpha,2^{\ast}_{\alpha}}_{0}(\mathrm{curl};\Omega)}\limits_{\nabla\times u\neq0}}\frac{\int_{\mathbb{R}^{3}}|\nabla\times u|^{2}dx}{\left(\int_{\mathbb{R}^{3}}|I_{\alpha/2}\ast|u+\widetilde{w}_{\mathbb{R}^{3}}(u)|^{2^{\ast}_{\alpha}}|^{2}dx\right)^{\frac{1}{2^{\ast}_{\alpha}}}},
\end{aligned}
\end{equation}
where $u\in W^{\alpha,2^{\ast}_{\alpha}}_{0}(\mathrm{curl};\Omega)$ is extended by $0
$ outside $\Omega$.
%but $v$ and $w$ need not be 0 there.
For constant $\bar{S}_{\mathrm{curl},HL}(\Omega)$ in domains $\Omega\neq\mathbb{R}^{3}$, it also can be characterized as
\begin{equation}\label{BARSCURLHL}
\begin{aligned}
\bar{S}_{\mathrm{curl},HL}(\Omega)&=\mathop{\mathrm{inf}}\limits_{\mathop{u\in W^{\alpha,2^{\ast}_{\alpha}}_{0}(\mathrm{curl};\Omega)}\limits_{\nabla\times u\neq0}}\mathop{\mathrm{sup}}\limits_{w\in\mathcal{W}_{\Omega}}\frac{\int_{\Omega}|\nabla\times u|^{2}dx}{\left(\int_{\Omega}|I_{\alpha/2}\ast|u+w|^{2^{\ast}_{\alpha}}|^{2}dx\right)^{\frac{1}{2^{\ast}_{\alpha}}}}\\
&=\mathop{\mathrm{inf}}\limits_{\mathop{u\in W^{\alpha,2^{\ast}_{\alpha}}_{0}(\mathrm{curl};\Omega)}\limits_{\nabla\times u\neq0}}\frac{\int_{\Omega}|\nabla\times u|^{2}dx}{\left(\int_{\Omega}|I_{\alpha/2}\ast|u+\widetilde{w}_{\Omega}(u)|^{2^{\ast}_{\alpha}}|^{2}dx\right)^{\frac{1}{2^{\ast}_{\alpha}}}}.
\end{aligned}
\end{equation}

To compare these sharp constants, we introduce the following set
\begin{equation}\label{N Omega}
\mathcal{N}_{\Omega}:=\{u\in W^{\alpha,2^{\ast}_{\alpha}}_{0}(\mathrm{curl};\Omega)\setminus\mathcal{W}_{\Omega}
;J'(u)u=0~\mathrm{and}~J'(u)|_{\mathcal{W}_{\Omega}}=0\}.
\end{equation}
According to the argument in \cite[Lemma 4.2]{Mederski2021}, we have $tu+\widetilde{w}_{\mathbb{R}^{3}}(tu)=t(u+\widetilde{w}_{\mathbb{R}^{3}}(u))$ , then we may assume without loss of generality that $u+\widetilde{w}_{\mathbb{R}^{3}}(u)\in\mathcal{N}$ in (\ref{SCURLHL}). By the maximality and uniqueness of $\widetilde{w}_{\Omega}(u)$, we easily deduce that the mapping $u\mapsto \widetilde{w}_{\Omega}(u)$ is also continuous. Therefore, we may assume that $u+\widetilde{w}_{\Omega}(u)\in\mathcal{N}_{\Omega}$ in (\ref{BARSCURLHL}). Then easily calculate that
\begin{equation}\label{three}
\mathop{\mathrm{inf}}\limits_{\mathcal{N}} J=\frac{2^{\ast}_{\alpha}-1}{2\cdot2^{\ast}_{\alpha}} S_{\mathrm{curl},HL}^{\frac{2^{\ast}_{\alpha}}{2^{\ast}_{\alpha}-1}},~~~
\mathop{\mathrm{inf}}\limits_{\mathcal{N}} J|_{W^{\alpha,2^{\ast}_{\alpha}}_{0}(\mathrm{curl};\Omega)}=\frac{2^{\ast}_{\alpha}-1}{2\cdot2^{\ast}_{\alpha}} S^{\frac{2^{\ast}_{\alpha}}{2^{\ast}_{\alpha}-1}}_{\mathrm{curl},HL}(\Omega),~~~
\mathop{\mathrm{inf}}\limits_{\mathcal{N}_{\Omega}} J=\frac{2^{\ast}_{\alpha}-1}{2\cdot2^{\ast}_{\alpha}} \bar{S}^{\frac{2^{\ast}_{\alpha}}{2^{\ast}_{\alpha}-1}}_{\mathrm{curl},HL}(\Omega).
\end{equation}
\begin{lem}\label{compare1}
$S_{\mathrm{curl},HL}(\Omega)\geq S_{\mathrm{curl},HL}$,~~$S_{\mathrm{curl},HL}(\Omega)\geq \bar{S}_{\mathrm{curl},HL}(\Omega)$.
\end{lem}
\begin{proof}
In view of Lemma \ref{gradient}, $W^{\alpha,2^{\ast}_{\alpha}}_{0}(\mathrm{curl};\Omega)\subset W^{\alpha,2^{\ast}_{\alpha}}_{0}(\mathrm{curl};\mathbb{R}^{3})$, we can easily observe from (\ref{SCURLHL}) and (\ref{ScurlHL}) that $S_{\mathrm{curl},HL}(\Omega)\geq S_{\mathrm{curl},HL}$. Similarly, since $\mathcal{W}_{\Omega}\subset\mathcal{W}_{\mathbb{R}^{3}}$, we can deduce that $S_{\mathrm{curl},HL}(\Omega)\geq \bar{S}_{\mathrm{curl},HL}(\Omega)$ from (\ref{SCURLHL}) and (\ref{BARSCURLHL}).
\end{proof}
To complete Theorem \ref{four constant}, we shall need the following inequality, which corresponds to the condition (B3), and the proof follows a similar argument in \cite[Lemma 4.1]{Mederski2018}.
\begin{lem}\label{tn}
If $u\in W^{\alpha,2^{\ast}_{\alpha}}_{0}(\mathrm{curl};\Omega)\setminus\{0\},u\in\mathcal{W}_{\Omega}$ and $t\geq0$, then
\begin{equation*}
J(u)\geq J(tu+w)-J'(u)\left[\frac{t^{2}-1}{2}u+tw\right].
\end{equation*}
Moreover, strict inequality holds provided $t=1$ and $w=0$. $(\Omega=\mathbb{R}^{3}~admitted.)$
\begin{proof}
By an explicit computation and using $\nabla\times w=0$, we show that
\begin{equation*}
J(u)-J(tu+w)+J'(u)\left[\frac{t^{2}-1}{2}u+tw\right]=\int_{\Omega}\varphi(t,x)dx,
\end{equation*}
where
\begin{equation*}
\begin{aligned}
\varphi(t,x)=&-\left\langle\left(I_{\alpha}\ast|u|^{2^{\ast}_{\alpha}}\right)|u(x)|^{2^{\ast}_{\alpha}-2}u(x) ,\frac{t^{2}-1}{2}u(x)+tw(x)\right\rangle-\frac{1}{2\cdot2^{\ast}_{\alpha}}\int_{\Omega}|I_{\alpha/2}\ast|u|^{2^{\ast}_{\alpha}}|^{2}dx
\\
&+\frac{1}{2\cdot2^{\ast}_{\alpha}}\int_{\Omega}|I_{\alpha/2}\ast|tu+w|^{2^{\ast}_{\alpha}}|^{2}dx.
\end{aligned}
\end{equation*}
It is easy to check that $\varphi(0,x)>0$ as $t=0$ and $\varphi(t,x)\longrightarrow\infty$ as $t\longrightarrow\infty$. Therefore, if there exist $t$ such that $\varphi(t,x)\leq0$, then there exists $t_{0}>0$ such that $\partial_{t}\varphi(t_{0},x)=0$, namely
\begin{equation*}
\begin{aligned}
\partial_{t}\varphi(t_{0},x)=-&\left\langle\left(I_{\alpha}\ast|u|^{2^{\ast}_{\alpha}}\right)|u(x)|^{2^{\ast}_{\alpha}-2}u(x) ,t_{0}u(x)+w(x)\right\rangle\\
&+\left\langle\left(I_{\alpha}\ast|t_{0}u+w|^{2^{\ast}_{\alpha}}\right)|t_{0}u(x)+w(x)|^{2^{\ast}_{\alpha}-2}\left(t_{0}u(x)+w(x)\right),u(x)\right\rangle=0,
\end{aligned}
\end{equation*}
then either $\langle u,t_{0}u+w\rangle$=0, i.e. $-\langle u,w\rangle=t_{0}\langle u,u\rangle=t_{0}|u|^{2}$, or $|u|=|t_{0}+w|$, i.e., $-t_{0}\langle u,w\rangle=\frac{t_{0}^{2}-1}{2}|u|^{2}+\frac{1}{2}|w|^{2}$.
In the first case, we obtain that
\begin{equation*}
\begin{aligned}
\varphi(t_{0},x)&=\left(\frac{t_{0}^{2}+1}{2}-\frac{1}{2\cdot2^{\ast}_{\alpha}}\right)\int_{\Omega}|I_{\alpha/2}\ast|u|^{2^{\ast}_{\alpha}}|^{2}dx+\frac{1}{2\cdot2^{\ast}_{\alpha}}\int_{\mathbb{R}{3}}|I_{\alpha/2}\ast|t_{0}u+w|^{2^{\ast}_{\alpha}}|^{2}dx>0.
\end{aligned}
\end{equation*}
And in the second case, we deduce that
\begin{equation*}
\begin{aligned}
\varphi(t_{0},x)&=\frac{1}{2}\int_{\Omega}(I_{\alpha}\ast|u|^{2^{\ast}_{\alpha}})|u(x)|^{2^{\ast}_{\alpha}-2}|w(x)|^{2}dx\geq0.
\end{aligned}
\end{equation*}
Hence $\varphi(t,x)\geq0$ for all $t\geq0$ and the inequality is strict if $w\neq0$. If $w=0$, we can see
\begin{equation*}
\varphi(t,x)=\left(\frac{t^{2\cdot2^{\ast}_{\alpha}}}{2\cdot2^{\ast}_{\alpha}}-\frac{t^{2}}{2}+\frac{1}{2}-\frac{1}{2\cdot2^{\ast}_{\alpha}}\right)\int_{\Omega}(I_{\alpha}\ast|u|^{2^{\ast}_{\alpha}})|u(x)|^{2^{\ast}_{\alpha}-2}|w(x)|^{2}dx>0
\end{equation*}
provided $t\neq1$.
\end{proof}

\end{lem}
\begin{lem}\label{compare2}
$S_{\mathrm{curl},HL}(\Omega)\leq S_{\mathrm{curl},HL}$.
\end{lem}
\begin{proof}
By Theorem \ref{attained}(a), $u$ is a minimizer for $J$ on $\mathcal{N}$, then we can find a sequence $(u_{n})\subset \mathcal{C}_{0}^{\infty}(\mathbb{R}^{3},\mathbb{R}^{3})$ such that $u_{n}\longrightarrow u$. By the Helmholtz decomposition, we have $u_{n}=v_{n}+w_{m}, v_{n}\in\mathcal{V}_{\mathbb{R}^{3}},w_{n}\in\mathcal{W}_{\mathbb{R}^{3}}$. Since  $u_{n}=v_{n}+w_{n}\longrightarrow u=v_{0}+\widetilde{w}_{\mathbb{R}^{3}}(v_{0})$ and therefore $v_{n}\longrightarrow v_{0}$, $w_{n}\longrightarrow \widetilde{w}_{\mathbb{R}^{3}}(v_{0})$. So $v_{0}\neq0$ and $v_{n}$ are bounded away from $0$ in $Q^{\alpha,2^{\ast}_{\alpha}}(\mathbb{R}^{3},\mathbb{R}^{3})$ due to $u\in \mathcal{N}$.

Assume without loss of generality that $0\in\Omega$. There exist $\lambda_{n}$ such that $\overline{u}_{n}$ given by $\overline{u}_{n}(x):=\lambda_{n}^{1\setminus2}u_{n}(\lambda_{n}x)$ are supported in $\Omega$, that is $\overline{u}_{n}(x)\in W^{\alpha,2^{\ast}_{\alpha}}_{0}(\mathrm{curl};\Omega)$. Set $\widetilde{w}_{\mathbb{R}^{3}}(\overline{u}_{n})\in\mathcal{W}_{\mathbb{R}^{3}}$ and choose $t_{n}$ so that $t_{n}(\overline{u}_{n}+\widetilde{w}_{\mathbb{R}^{3}}(\overline{u}_{n}))\in\mathcal{N}$, then
\begin{equation}\label{tn1}
t_{n}^{2}=\frac{\left(\int_{\mathbb{R}^{3}}|\nabla\times \overline{u}_{n}|^{2}dx\right)^{\frac{1}{2^{\ast}_{\alpha}-1}}}
{\left(\int_{\mathbb{R}^{3}}|I_{\alpha/2}\ast|\overline{u}_{n}+\widetilde{w}_{\mathbb{R}^{3}}(\overline{u}_{n})|^{2^{\ast}_{\alpha}}|^{2}dx\right)^{\frac{1}{2^{\ast}_{\alpha}-1}}}.
\end{equation}
Since the Riesz potential is invariant with respect to translation, we have $||\overline{u}_{n}||=||u_{n}||$ and
\begin{equation*}
\begin{aligned}
\int_{\mathbb{R}^{3}}|I_{\alpha/2}\ast|\overline{u}_{n}+\widetilde{w}_{\mathbb{R}^{3}}(\overline{u}_{n})|^{2^{\ast}_{\alpha}}|^{2}dx&=\int_{\mathbb{R}^{3}}|I_{\alpha/2}\ast|u_{n}+\widetilde{w}_{\mathbb{R}^{3}}(u_{n})|^{2^{\ast}_{\alpha}}|^{2}dx\\
&=\int_{\mathbb{R}^{3}}|I_{\alpha/2}\ast|v_{n}+\widetilde{w}_{\mathbb{R}^{3}}(v_{n})|^{2^{\ast}_{\alpha}}|^{2}dx.
\end{aligned}
\end{equation*}
Therefore, as $(u_{n})$ is bounded, we have $(\overline{u}_{n})$ and $(\widetilde{w}_{\mathbb{R}^{3}}(\overline{u}_{n}))$ are bounded away from $0$, so is $||\overline{u}_{n}(x)+\widetilde{w}_{\mathbb{R}^{3}}(\overline{u}_{n})(x)||_{Q^{\alpha,2^{\ast}_{\alpha}}}$. Then we deduce that $(t_{n})$ is bounded, hence so is $(t_{n}^{2})$. Moreover, since $J(\overline{u}_{n})=J(u_{n})\longrightarrow \frac{2^{\ast}_{\alpha}-1}{2\cdot2^{\ast}_{\alpha}} S_{\mathrm{curl},HL}^{\frac{2^{\ast}_{\alpha}}{2^{\ast}_{\alpha}-1}}$ and $||J'(\overline{u}_{n})||=||J'(u_{n})||\longrightarrow 0$, it follows from Lemma \ref{tn} that
\begin{equation*}
\begin{aligned}
\frac{2^{\ast}_{\alpha}-1}{2\cdot2^{\ast}_{\alpha}} S_{\mathrm{curl},HL}^{\frac{2^{\ast}_{\alpha}}{2^{\ast}_{\alpha}-1}}
&=\mathop{\mathrm{lim}}\limits_{n\longrightarrow\infty}J(\overline{u}_{n})\geq\mathop{\mathrm{lim}}\limits_{n\longrightarrow\infty}\left(J(t_{n}(\overline{u}_{n}+\widetilde{w}_{\mathbb{R}^{3}}(\overline{u}_{n})))-J'(\overline{u}_{n})\left[\frac{t^{2}_{n}-1}{2}\overline{u}_{n}
+t^{2}_{n}\widetilde{w}_{\mathbb{R}^{3}}(\overline{u}_{n}))\right]\right)\\
&=\mathop{\mathrm{lim}}\limits_{n\longrightarrow\infty}J(t_{n}(\overline{u}_{n}+\widetilde{w}_{\mathbb{R}^{3}}(\overline{u}_{n})))\geq\frac{2^{\ast}_{\alpha}-1}{2\cdot2^{\ast}_{\alpha}} S^{\frac{2^{\ast}_{\alpha}}{2^{\ast}_{\alpha}-1}}_{\mathrm{curl},HL}(\Omega).
\end{aligned}
\end{equation*}
The last inequality follows from the fact that $\overline{u}_{n}\in W^{\alpha,2^{\ast}_{\alpha}}_{0}(\mathrm{curl};\Omega)$.
\end{proof}

\begin{proof}[\bf Complete of the Proof of Theorem \ref{four constant}]
Repeating the proof of Theorem \ref{attained} (b) with obvious changes, namely, change the domain $\mathbb{R}^{3}$ into $\Omega$, change $S_{\mathrm{curl},HL}$ into $\bar{S}_{\mathrm{curl},HL}(\Omega)$, we have $\bar{S}_{\mathrm{curl},HL}(\Omega)\geq S_{HL}$. Since the optimal function for $\bar{S}_{\mathrm{curl},HL}(\Omega)$ is not found in our process, we can not exclude the case $\bar{S}_{\mathrm{curl},HL}(\Omega)=S_{HL}$. As a consequently, we complete the proof of Theorem \ref{four constant} by Lemma \ref{compare1} and Lemma \ref{compare2}.
\end{proof}

\section{\bf Proof of Theorem \ref{main thm}}
%For the Brezis-Nirenberg problem on a bounded domain $\Omega$, we may recover the compactness for %$(PS)_{c}$ sequence under the critical energy level. While the energy level is varying with respect %to the parameter $\lambda$. Therefore, in this section we drop the subscript $\Omega$ from notation %and replace it by $\lambda$ ($J_{\lambda}$, $\mathcal{N}_{\lambda}$ etc), also, we note %$\mathcal{V}$, $\mathcal{W}$ for $\mathcal{V}$, $\mathcal{W}$.

According to the spectrum analysis of the curl-curl operator in the introduction, for $\lambda\leq0$, we find two closed and orthogonal subspaces $\mathcal{V}_{\Omega}^{+}$ and $\widetilde{\mathcal{V}}_{\Omega}$ of $\mathcal{V}_{\Omega}$ such that the quadratic form $Q:\mathcal{V}_{\Omega}\longrightarrow\mathbb{R}$ given by
\begin{equation}\label{quad}
Q(v):=\int_{\Omega}(|\nabla\times v|^{2}+\lambda|v|^{2})dx=\int_{\Omega}(|\nabla v|^{2}+\lambda|v|^{2})dx
\end{equation}
is positive defined on $\mathcal{V}_{\Omega}^{+}$ and negative semi-definite on $\widetilde{\mathcal{V}}_{\Omega}$ where $\mathrm{dim}\widetilde{\mathcal{V}}_{\Omega}<\infty$. Writing $u=v+w=v^{+}+\widetilde{v}+w\in\mathcal{V}_{\Omega}^{+}\oplus\widetilde{\mathcal{V}}_{\Omega}\oplus\mathcal{W}_{\Omega}$, the functional $J_{\lambda}$ (see (\ref{functional})) can be expressed as
%we have \begin{equation*}
%Q(v)=Q(v^{+})+Q(\widetilde{v}),
%\end{equation*}
%and the functional $J_{\lambda}$ (see (\ref{functional})) can be expressed as
%\begin{equation*}
%\begin{aligned}
%J_{\lambda}(u)&=\frac{1}{2}\int_{\Omega}|\nabla\times %v|^{2}dx+\frac{\lambda}{2}\int_{\Omega}|v+w|^{2}dx-\frac{1}{2\cdot2^{\ast}_{\alpha}}\int_{\Omega}|I_{\alpha/2}\ast|u|^{2^{\ast}_{\alpha}}|^{2}dx\\
%&=\frac{1}{2}Q(v^{+})+\frac{1}{2}Q(\widetilde{v})+\frac{\lambda}{2}\int_{\Omega}|w|^{2}dx
%-\frac{1}{2\cdot2^{\ast}_{\alpha}}\int_{\Omega}|I_{\alpha/2}\ast|u|^{2^{\ast}_{\alpha}}|^{2}dx.\\
%\end{aligned}
%\end{equation*}
%According to the abstract critical point theorem, we rewrite it as the following form
\begin{equation*}
\begin{aligned}
J_{\lambda}(u)
%&=\frac{1}{2}\int_{\Omega}|\nabla\times v|^{2}dx+\frac{\lambda}{2}\int_{\Omega}|v+w|^{2}dx-\frac{1}{2\cdot2^{\ast}_{\alpha}}\int_{\Omega}|I_{\alpha/2}\ast|u|^{2^{\ast}_{\alpha}}|^{2}dx\\
&=\frac{1}{2}||v^{+}||^{2}+\frac{1}{2}||\widetilde{v}||^{2}+\frac{\lambda}{2}\int_{\Omega}(|v|^{2}+|w|^{2})dx-\frac{1}{2\cdot2^{\ast}_{\alpha}}\int_{\Omega}|I_{\alpha/2}\ast|u|^{2^{\ast}_{\alpha}}|^{2}dx\\
&=\frac{1}{2}||v^{+}||^{2}-I_{\lambda}(v+w),
\end{aligned}
\end{equation*}
where
\begin{equation*}
I_{\lambda}(v+w)=-\frac{1}{2}||\widetilde{v}||^{2}-\frac{\lambda}{2}\int_{\Omega}(|v|^{2}+|w|^{2})dx+\frac{1}{2\cdot2^{\ast}_{\alpha}}\int_{\Omega}|I_{\alpha/2}\ast|u|^{2^{\ast}_{\alpha}}|^{2}dx.
\end{equation*}

Similarly as in \cite{Bartsch2015-1}, we shall show that $J_{\lambda}$ satisfies the assumptions (A1)-(A4), (B1)-(B3) and (C1)-(C3) from Section 2.2.

\begin{lem}\label{check}
Conditions (A1)-(A4), (B1)-(B3) and (C2) in Lemma \ref{abstract} hold for $J_{\lambda}$.
\end{lem}
\begin{proof}
i) By Lemma \ref{well define}, we have $I_{\lambda}$ is of class $\mathcal{C}^{1}$. Since $Q(v)$ is negative on $\widetilde{\mathcal{V}}_{\Omega}$, $2^{\ast}_{\alpha}$ is a upper critical index, we have $I_{\lambda}(u)\geq I_{\lambda}(0)=0$ for any $u\in W^{\alpha,2^{\ast}_{\alpha}}_{0}(\mathrm{curl};\Omega)$.

ii) Since $I_{\lambda}$ is convex, $I_{\lambda}$ is $\mathcal{T}-$sequentialy lower semicontinuous. Hence, (A2) holds.

iii) We easily check (A3), since $u_{n}\rightharpoonup u_{0}$ in $Q^{\alpha,2^{\ast}_{\alpha}}(\Omega,\mathbb{R}^{3})$, and $I_{\lambda}(u_{n})\longrightarrow I_{\lambda}(u_{0})$ imply $|u_{n}|_{Q^{\alpha,2^{\ast}_{\alpha}}}\longrightarrow |u_{0}|_{Q^{\alpha,2^{\ast}_{\alpha}}}$, thus $u_{n}\longrightarrow u_{0}$ in $Q^{\alpha,2^{\ast}_{\alpha}}(\Omega,\mathbb{R}^{3})$.

iv) Since $\mathcal{V}_{\Omega}$ is a Hilbert space, the HLS inequality is still valid on there, then for any $u\in\mathcal{V}_{\Omega}^{+}$, we have
\begin{equation*}
\begin{aligned}
J(u)=J(v,0)&=\frac{1}{2}||v||_{\mathcal{V}_{\Omega}}^{2}+\frac{\lambda}{2}|v|^{2}_{2}-\frac{1}{2\cdot2^{\ast}_{\alpha}}\int_{\Omega}|I_{\alpha/2}\ast|v|^{2^{\ast}_{\alpha}}|^{2}dx\\
&\geq\frac{1}{2}||v||_{\mathcal{V}_{\Omega}}^{2}+\frac{\lambda}{2}|v|^{2}_{2}-\frac{1}{2\cdot2^{\ast}_{\alpha}}\left(\int_{\Omega}|v|^{6}dx\right)^{\frac{1}{6}\ast2}\\
&\geq\frac{\delta}{2}||v||_{\mathcal{V}_{\Omega}}^{2}-\varepsilon|v|^{2}-c_{\varepsilon}|v|^{3}_{6}\geq\frac{\delta}{4}||v||^{2}_{\mathcal{V}_{\Omega}}-C_{1}||v||_{\mathcal{V}_{\Omega}}^{3}
\end{aligned}
\end{equation*}
for some $\delta, C_{1}>0$.

v) Condition (B1) follows from Lemma 5.1 (c) in \cite{Bartsch2015-1}. Suppose that $(||v^{+}_{n}||_{\mathcal{V}_{\Omega}})_{n}$ is bounded and $||(v_{n},w_{n})||\longrightarrow\infty$ as $n\longrightarrow\infty$. Since $\mathrm{dim}(\widetilde{\mathcal{V}}_{\Omega})<\infty$ there holds $|v_{n}+w_{n}|_{Q^{\alpha,2^{\ast}_{\alpha}}}\longrightarrow\infty$. Moreover by the orthogonality $\mathcal{V}_{\Omega}^{+}\perp\widetilde{\mathcal{V}}_{\Omega}$ in $L^{2}(\Omega,\mathbb{R}^{3})$ and $\mathcal{V}_{\Omega}\perp\mathcal{W}_{\Omega}$ in $Q^{\alpha,2^{\ast}_{\alpha}}(\Omega,\mathbb{R}^{3})$, we have
\begin{equation}\label{orth}
||\widetilde{v}_{n}||_{\mathcal{V}_{\Omega}}^{2}\leq C_{1}|\widetilde{v}_{n}|_{2}^{2}\leq C_{1}|v_{n}|^{2}_{2}\leq C_{1}|v_{n}+ w_{n}|_{2}^{2}\leq C_{2}|v_{n}+ w_{n}|_{Q^{\alpha,2^{\ast}_{\alpha}}}^{2}
\end{equation}
for some $0<C_{1}<C_{2}$. This implies
\begin{equation*}
\begin{aligned}
I(v_{n},w_{n})&=-\frac{1}{2}||\widetilde{v}_{n}||-\frac{\lambda}{2}|v_{n}+ w_{n}|^{2}_{2}+\frac{1}{2\cdot2^{\ast}_{\alpha}}|v_{n}+ w_{n}|_{Q^{\alpha,2^{\ast}_{\alpha}}}^{2^{\ast}_{\alpha}}\\
&\geq -\frac{C_{2}}{2}|v_{n}+ w_{n}|_{Q^{\alpha,2^{\ast}_{\alpha}}}^{2}+\frac{1}{2\cdot2^{\ast}_{\alpha}}|v_{n}+ w_{n}|_{Q^{\alpha,2^{\ast}_{\alpha}}}^{2^{\ast}_{\alpha}}\longrightarrow\infty,
\end{aligned}
\end{equation*}
because $|v_{n}+\nabla w_{n}|_{Q^{\alpha,2^{\ast}_{\alpha}}}\longrightarrow\infty$.

vi) This part we check condition (B2) and (C2). By (\ref{orth}), we have
\begin{equation*}
\begin{aligned}
&I(t_{n}(v_{n}+w_{n}))\\
&=\frac{1}{2}||t_{n}\widetilde{v}_{n}||_{\mathcal{V}_{\Omega}}^{2}-\frac{\lambda}{2}|t_{n}(v_{n}+w_{n})|^{2}_{2}+\frac{1}{2\cdot2^{\ast}_{\alpha}}\int_{\Omega}|I_{\alpha/2}\ast|t_{n}(v_{n}+w_{n})|^{2^{\ast}_{\alpha}}|^{2}dx\\
&\geq-\frac{1}{2}t^{2}_{n}||\widetilde{v}_{n}||^{2}_{\mathcal{V}_{\Omega}}-\frac{\lambda}{2}t^{2}_{n}|v_{n}+w_{n}|^{2}_{2}+t_{n}^{2\cdot2^{\ast}_{\alpha}}\frac{1}{2\cdot2^{\ast}_{\alpha}}\int_{\Omega}|I_{\alpha/2}\ast|v_{n}+w_{n}|^{2^{\ast}_{\alpha}}|^{2}dx\\
&\geq-\frac{C_{2}}{2}t^{2}_{n}||{v}_{n}+w_{n}||^{2}_{Q_{\alpha,2^{\ast}_{\alpha}}}+t_{n}^{2\cdot2^{\ast}_{\alpha}}\frac{1}{2\cdot2^{\ast}_{\alpha}}||{v}_{n}+w_{n}||^{2\cdot2^{\ast}_{\alpha}}_{Q_{\alpha,2^{\ast}_{\alpha}}}.
\end{aligned}
\end{equation*}
Then
\begin{equation*}
I(t_{n}(v_{n}+w_{n}))/ t^{2}_{n}\geq-\frac{C_{2}}{2}||\widetilde{v}_{n}||^{2}_{Q_{\alpha,2^{\ast}_{\alpha}}}+t_{n}^{2\cdot2^{\ast}_{\alpha}-2}\frac{1}{2\cdot2^{\ast}_{\alpha}}||\widetilde{v}_{n}||^{2\cdot2^{\ast}_{\alpha}}_{Q_{\alpha,2^{\ast}_{\alpha}}}.
\end{equation*}
If $||(v_{n},w_{n})||\longrightarrow\infty$ then $I(t_{n}(v_{n}+w_{n}))/ t^{2}_{n}\longrightarrow\infty$. If $(||(v_{n},w_{n})||)_{n}$ is bounded. Then $(|v_{n}+w_{n}|_{Q^{\alpha,2^{\ast}_{\alpha}}})_{n}$ is bounded. If $|v_{n}+w_{n}|_{Q^{\alpha,2^{\ast}_{\alpha}}}\longrightarrow0$, then $|v_{n}+w_{n}|_{2}\longrightarrow0$ and by the orthogonality in $L^{2}(\Omega,\mathbb{R}^{3})$ which contradicts $u_{0}\neq0$. Therefore $\frac{t_{n}^{2\cdot2^{\ast}_{\alpha}-2}}{2\cdot2^{\ast}_{\alpha}}||\widetilde{v}_{n}||^{2\cdot2^{\ast}_{\alpha}}_{Q_{\alpha,2^{\ast}_{\alpha}}}\longrightarrow\infty$ as $n\longrightarrow\infty$ and again $I(t_{n}(v_{n}+w_{n}))/ t^{2}_{n}\longrightarrow\infty$.

vii) Condition (B3) follows from Lemma \ref{tn} by changing $J(u)$ into $J_{\lambda}(u)$.
\end{proof}

To apply the Concentration compactness lemma, we set $\widetilde{\mathcal{W}}:=\widetilde{\mathcal{V}}_{\Omega}\oplus\mathcal{W}_{\Omega}$ with $\widetilde{w}_{\Omega}=\widetilde{v}+w$, where $\widetilde{\mathcal{V}}_{\Omega}=Z$, see Section 2.2. On the other hand, we shall extend $\mathcal{V}_{\Omega}^{+}$ into $\mathcal{V}_{\mathbb{R}^{3}}$, which is a closed subspaces of $\mathcal{D}^{1,2}(\mathbb{R}^{3},\mathbb{R}^{3})$. Indeed, let $U$ be a bounded domain in $\mathbb{R}^{3}$, $\bar{\Omega}\subset U$. Since $\mathcal{V}_{\Omega}\subset H^{1}(\Omega,\mathbb{R}^{3})$, then each $v\in\mathcal{V}_{\Omega}$ may be extended to $v'\in H^{1}_{0}(U,\mathbb{R}^{3})$ such that $v'|_{\Omega} =v$. This extension is bounded as a mapping from $\mathcal{V}_{\Omega}$ to $H^{1}_{0}(U,\mathbb{R}^{3})$.  Since
\begin{equation*}
\mathcal{V}':=\{v'\in H^{1}_{0}(U,\mathbb{R}^{3}):v'|_{\Omega}\in\mathcal{V}_{\Omega}\}
\end{equation*}
is a closed subspace of $H^{1}_{0}(U,\mathbb{R}^{3})$, and hence of $\mathcal{D}^{1,2}(\mathbb{R}^{3},\mathbb{R}^{3})$, we then can apply Lemma \ref{concentration} with $\mathcal{V}_{\Omega}^{+}$ replacing $\mathcal{V}_{\mathbb{R}^{3}}$. Set the generalized Nehari-Pankov manifold as follow
\begin{equation}\label{N lambda}
\mathcal{N}_{\lambda}:=\{u\in W^{\alpha,2^{\ast}_{\alpha}}_{0}(\mathrm{curl};\Omega)\setminus(\widetilde{\mathcal{V}}_{\Omega}\oplus\mathcal{W}_{\Omega})
:J'_{\lambda}(u)|_{\mathcal{R}u\oplus\widetilde{\mathcal{V}}_{\Omega}\oplus\widetilde{\mathcal{W}}_{\Omega}}=0\}.
\end{equation}
\begin{lem}\label{weak4}
$J'$ is  weak-to-weak$^{\ast}$ continuous on $\mathcal{N_{\lambda}}$ and condition (C3) in Lemma \ref{abstract} holds.
\end{lem}
\begin{proof}
%Let $(u_{n})$ be a $(PS)_{\beta}-$sequence such that $(u_{n})\subset\mathcal{N}_{\lambda}$. By Lemma \ref{coercive}, $(u_{n})$ is bounded and we can assume
Suppose that $u_{n}\rightharpoonup u_{0}$ in $W^{\alpha,2^{\ast}_{\alpha}}_{0}(\mathrm{curl};\Omega)$. Set $u_{n}=m_{\lambda}(v^{+}_{n})=v^{+}+\widetilde{w}_{\Omega}(v^{+}_{n})$. Since $\mathcal{V}_{\Omega}^{+}$ and $\widetilde{\mathcal{W}}_{\Omega}$ are complementary subspaces, $v^{+}_{n}$ is bounded in $\mathcal{V}_{\Omega}^{+}$. Then passing to a subsequence, we have $v^{+}_{n}\rightharpoonup v^{+}_{0}$ in $\mathcal{V}_{\Omega}^{+}$ , $v^{+}_{n}\longrightarrow v^{+}_{0}$ in $L^{2}(\Omega,\mathbb{R}^{3})$ and a.e. in $\Omega$. Therefore, by the Concentration Compactness Lemma \ref{concentration}, we have $\widetilde{w}_{\Omega}(v^{+}_{n})\longrightarrow \widetilde{w}_{\Omega}(v^{+}_{0})$ in $L^{2}(\Omega,\mathbb{R}^{3})$ and also a.e. in $\Omega$. Hence, we also have $u_{n}\longrightarrow u_{0}$ a.e. in $\Omega$. Then by the Viltali convergence principle, $J'_{\lambda}$ is  weak-to-weak$^{\ast}$ continuous. Moreover, by the lower semi-continuity of $I_{\lambda}$, (C3) holds.
%as in the proof of (\ref{weak}), we have $J^{'}_{\lambda}(u_{0})=0$, this implies that $u_{0}$ is a solution for (\ref{nonlocal case}).
\end{proof}
Now, we set a compactly perturbed problem. Take $X^{0}:=\widetilde{V}_{\Omega}$, $X^{1}:=\mathcal{W}_{\Omega}$ and let us consider the functional $J_{cp}: X =\mathcal{V}_{\Omega}\oplus\mathcal{W}_{\Omega}\longrightarrow \mathbb{R}$ given by
\begin{equation*}\label{Jcp}
J_{cp}(u)=J_{0}(u)=\frac{1}{2}\int_{\Omega}|\nabla\times u|^{2}dx+\frac{1}{2\cdot2^{\ast}_{\alpha}}\int_{\Omega}|I_{\alpha/2}\ast|u|^{2^{\ast}_{\alpha}}|^{2}dx.
\end{equation*}
Moreover we define the corresponding Nehari-Pankov manifold
\begin{equation}\label{Ncp}
\mathcal{N}_{cp} = \{E\in (\mathcal{V}_{\Omega}\oplus\mathcal{W}_{\Omega})\setminus\mathcal{W}_{\Omega}: J'_{cp}(u)|_{\mathbb{R}u\oplus\mathcal{W}_{\Omega}}=0\}.
\end{equation}
Observe that as in Lemma \ref{check} we show that $J_{cp}$ satisfies the corresponding condition (A1)-(A4) and (B1)-(B3).
\begin{lem}\label{C1}
Condition (C1) in Lemma \ref{abstract}  holds.
\end{lem}
\begin{proof}
For any bounded sequence $u_{n}\subset\mathcal{N}_{cp}$, we have $u_{n}\rightharpoonup u$ in $\mathcal{N}_{cp}$. By the concentrated compactness Lemma \ref{concentration}, we have $u_{n}\longrightarrow u$ in $L^{2}(\Omega,\mathbb{R}^{3})$. Since $J_{\lambda}(u)-J_{cp}(u)=\frac{\lambda}{2}\int_{\Omega}|u|^{2}dx$, we have condition (C1) holds.
\end{proof}

\begin{lem}\label{coercive}
$J_{\lambda}$ is coercive on $\mathcal{N}_{\lambda}$ and $J_{cp}$ is coercive on $\mathcal{N}_{cp}$.
\end{lem}
\begin{proof}
The proof is similar to Lemma 4.6 in \cite{Mederski2018}. Let $u_{n}=v_{n}+w_{n}\in\mathcal{N}_{\lambda}$ and suppose that $||u_{n}||\longrightarrow\infty$. Observe that
\begin{equation*}
J_{\lambda}(u_{n})=J_{\lambda}(u_{n})-\frac{1}{2}J'_{\lambda}(u_{n})(u_{n})=(\frac{1}{2}-\frac{1}{2\cdot2^{\ast}_{\alpha}})|u_{n}|_{Q^{\alpha,2^{\ast}_{\alpha}}}^{2\cdot2^{\ast}_{\alpha}}\geq C_{1}|w_{n}|_{Q^{\alpha,2^{\ast}_{\alpha}}}^{2\cdot2^{\ast}_{\alpha}}
\end{equation*}
for some constant $C_{1}>0$, since $\mathcal{W}_{\Omega}$ is closed,  $\mathrm{cl}\mathcal{V}_{\Omega}\cap\mathcal{W}_{\Omega}=\{0\}$ in $Q^{\alpha,2^{\ast}_{\alpha}}(\Omega,\mathbb{R}^{3})$ and the projection $\mathrm{cl}\mathcal{V}_{\Omega}\oplus\mathcal{W}_{\Omega}=\{0\}$ onto $\mathcal{W}_{\Omega}$ is continuous. Hence, if $|u_{n}|_{Q^{\alpha,2^{\ast}_{\alpha}}}\longrightarrow\infty$, then $J_{\lambda}(u_{n})\longrightarrow\infty$ as $n\longrightarrow\infty$. Suppose that $|u_{n}|_{Q^{\alpha,2^{\ast}_{\alpha}}}$ is bounded, Then $||v_{n}||\longrightarrow\infty$ and
\begin{equation*}
\begin{aligned}
J_{\lambda}(u_{n})&=J_{\lambda}(u_{n})-\frac{1}{2\cdot2^{\ast}_{\alpha}}J'_{\lambda}(u_{n})(u_{n})=(\frac{1}{2}-\frac{1}{2\cdot2^{\ast}_{\alpha}})\left(\int_{\Omega}|\nabla\times v_{n}|^{2}dx+\lambda\int_{\Omega}|v_{n}+w_{n}|^{2}dx\right)\\
&\geq(\frac{1}{2}-\frac{1}{2\cdot2^{\ast}_{\alpha}})\left(\int_{\Omega}|\nabla\times v_{n}|^{2}dx+\lambda C_{2}|u_{n}|_{Q^{\alpha,2^{\ast}_{\alpha}}}^{2}\right),
\end{aligned}
\end{equation*}
for some constant $C_{2}>0$. Thus $J_{\lambda}(u_{n})\longrightarrow\infty$. Similarly, we show that $J_{cp}$ is coercive on $\mathcal{N}_{cp}$.
\end{proof}
\begin{lem}\label{ground}
Let $c_{\lambda}<c_{0}$ and $(u_{n})_{n\in\mathbb{N}}\subset\mathcal{N}_{\lambda}$ be a Palais-Smale sequence at $c_{\lambda}$, i.e. $J_{\lambda}(u_{n})\longrightarrow c_{\lambda}$ and $J'_{\lambda}(u_{n})\longrightarrow 0$ as $n\longrightarrow\infty$. Then $u_{n}\rightharpoonup u_{0}\neq0$ for some $u_{0}$ in $W^{\alpha,2^{\ast}_{\alpha}}(\mathrm{curl};\Omega)$. Moreover, $c_{\lambda}$ is achieved by a critical point of $J_{\lambda}$.
\end{lem}
\begin{proof}
The conclusion follows from Lemma \ref{check}, Lemma \ref{C1}, Lemma \ref{weak4}, Lemma \ref{coercive} and Lemma \ref{abstract}(a)(b).
\end{proof}
As we introduced before, we shall verify the $(PS)_{c}$ condition. Similar to \cite[Lemma 6.4]{Mederski2021} we need the following version of the Brezis-Lieb lemma. Setting
\begin{equation*}
N(u)=\left(I_{\alpha}\ast|u|^{2^{\ast}_{\alpha}}\right)|u(x)|^{2^{\ast}_{\alpha}-2}u(x),
\end{equation*}
then we have the following lemma.
\begin{lem}\label{Brezis}
Suppose $(u_{n})$ is bounded in $Q^{\alpha,2^{\ast}_{\alpha}}(\Omega,\mathbb{R}^{3})$ and $u_{n}\longrightarrow u$ a.e. in $\Omega$. Then
\begin{equation*}
N(u_{n})-N(u_{n}-u)\longrightarrow N(u)~~in~(Q^{\alpha,2^{\ast}_{\alpha}}(\Omega,\mathbb{R}^{3}))'~as~n\longrightarrow\infty.
\end{equation*}
\end{lem}
\begin{proof}
By the proof of Lemma \ref{well define}, we have
 $N(u):Q^{\alpha,2^{\ast}_{\alpha}}(\Omega,\mathbb{R}^{3})\longrightarrow(Q^{\alpha,2^{\ast}_{\alpha}}(\Omega,\mathbb{R}^{3}))'$. Therefore, it turns to prove that $G(u_{n})-G(u_{n}-u)\longrightarrow G(u)$ in $L^{\frac{2\cdot2^{\ast}_{\alpha}}{2\cdot2^{\ast}_{\alpha}-1}}(\Omega,L^{\frac{2^{\ast}_{\alpha}}{2^{\ast}_{\alpha}-1}}(\Omega))$ . Since $u_{n}\longrightarrow u$ a.e. in $\Omega$, we have $G(u_{n})-G(u_{n}-u)\longrightarrow G(u)$ a.e. in $\Omega$. Since $u_{n}$ is bounded in $Q^{\alpha,2^{\ast}_{\alpha}}(\Omega,\mathbb{R}^{3})$, we have $G(u_{n})$ is bounded in $L^{\frac{2\cdot2^{\ast}_{\alpha}}{2\cdot2^{\ast}_{\alpha}-1}}(\Omega,L^{\frac{2^{\ast}_{\alpha}}{2^{\ast}_{\alpha}-1}}(\Omega))$ , so is $G(u_{n})-G(u_{n}-u)$. Then we have $G(u_{n})-G(u_{n}-u)\rightharpoonup G(u)$. Therefore, we only need to prove that
 \begin{equation*}
 |G(u_{n})-G(u_{n}-u)|_{L^{\frac{2\cdot2^{\ast}_{\alpha}}{2\cdot2^{\ast}_{\alpha}-1}}(\Omega,L^{\frac{2^{\ast}_{\alpha}}{2^{\ast}_{\alpha}-1}}(\Omega))}\rightarrow |G(u)|_{L^{\frac{2\cdot2^{\ast}_{\alpha}}{2\cdot2^{\ast}_{\alpha}-1}}(\Omega,L^{\frac{2^{\ast}_{\alpha}}{2^{\ast}_{\alpha}-1}}(\Omega))}.
 \end{equation*}

Indeed, let $A=\left(\frac{1}{|x-y|^{\frac{2^{\ast}_{\alpha}}{2^{\ast}_{\alpha}-1}N-\frac{2\cdot2^{\ast}_{\alpha}}{2\cdot2^{\ast}_{\alpha}}\alpha}}\right)^{\frac{2^{\ast}_{\alpha}}{2^{\ast}_{\alpha}-1}}$. Then by using Vitali's convergence theorem we obtain
\begin{equation*}
\begin{aligned}
&\int_{\Omega}\left(\int_{\Omega}|G(u_{n})-G(u_{n}-u)|^{\frac{2^{\ast}_{\alpha}}{2^{\ast}_{\alpha}-1}}dy\right)^{\frac{2\cdot2^{\ast}_{\alpha}}{2\cdot2^{\ast}_{\alpha}-1}}dx\\
&=\int_{\Omega}\Bigg(\int_{\Omega}(A)
\left(|u_{n}(y)|^{2^{\ast}_{\alpha}}|u_{n}(x)|^{2^{\ast}_{\alpha}-1}-|u_{n}(y)-u(y)|^{2^{\ast}_{\alpha}}|u_{n}(x)-u (x)|^{2^{\ast}_{\alpha}-1}\right)^{\frac{2^{\ast}_{\alpha}}{2^{\ast}_{\alpha}-1}}dy\Bigg)^{\frac{2\cdot2^{\ast}_{\alpha}}{2\cdot2^{\ast}_{\alpha}-1}}dx\\
&=\int_{\Omega}\left(\int_{\Omega}(A)
\left(\int_{0}^{1}\frac{\mathrm{d}}{\mathrm{d}t}\left(|u_{n}(y)+(t-1)u(y)|^{2^{\ast}_{\alpha}}|u_{n}(x)-(t-1)u(x)|^{2^{\ast}_{\alpha}-1}\right)^{\frac{2^{\ast}_{\alpha}}{2^{\ast}_{\alpha}-1}}dt\right)dy\right)^{\frac{2\cdot2^{\ast}_{\alpha}}{2\cdot2^{\ast}_{\alpha}-1}}dx\\
%&=\int_{\Omega}\left(\int_{\Omega}(A)
%\left(\int_{0}^{1}\frac{\mathrm{d}}{\mathrm{d}t}\left(|u_{n}(y)+(t-1)u(y)|^{\frac{(2^{\ast}_{\alpha})^{2}}{2^{\ast}_{\alpha}-1}}|u_{n}(x)-(t-1)u(x)|^{2^{\ast}_{\alpha}}\right)dt\right)dy\right)^{\frac{2\cdot2^{\ast}_{\alpha}}{2\cdot2^{\ast}_{\alpha}-1}}dx\\
&=\int_{0}^{1}\Bigg[\int_{\Omega}\bigg(\int_{\Omega}(A)
\bigg(\frac{(2^{\ast}_{\alpha})^{2}}{2^{\ast}_{\alpha}-1}\left\langle|u_{n}(y)+(t-1)u(y)|^{\frac{(2^{\ast}_{\alpha})^{2}}{2^{\ast}_{\alpha}-1}-2}|u_{n}(x)+(t-1)u(x)|^{2^{\ast}_{\alpha}}(u_{n}(y)+(t-1)u(y)),u(y)
\right\rangle\\
&~~~~~~~~+2^{\ast}_{\alpha}\left\langle|u_{n}(y)+(t-1)u(y)|^{\frac{(2^{\ast}_{\alpha})^{2}}{2^{\ast}_{\alpha}-1}}|u_{n}(x)+(t-1)u(x)|^{2^{\ast}_{\alpha}-2}(u_{n}(x)+(t-1)u(x)),u(x)
\right\rangle\bigg)dy\bigg)^{\frac{2\cdot2^{\ast}_{\alpha}}{2\cdot2^{\ast}_{\alpha}-1}}dx\Bigg]dt\\
&\longrightarrow\int_{0}^{1}\Bigg[\int_{\Omega}\bigg(\int_{\Omega}(A)
\bigg(\frac{(2^{\ast}_{\alpha})^{2}}{2^{\ast}_{\alpha}-1}\left\langle|tu(y)|^{\frac{(2^{\ast}_{\alpha})^{2}}{2^{\ast}_{\alpha}-1}-2}|tu(x)|^{2^{\ast}_{\alpha}}(tu(y)),u(y)
\right\rangle\\
&~~~~~~~~~~~~~~~~~~~~~~~~~~~~~~~~~~~~~~~~~~~~~~~~~~~~~~~~~~~~~~~+2^{\ast}_{\alpha}\left\langle|tu(y)|^{\frac{(2^{\ast}_{\alpha})^{2}}{2^{\ast}_{\alpha}-1}}|tu(x)|^{2^{\ast}_{\alpha}-2}(tu(x)),u(x)
\right\rangle\bigg)dy\bigg)^{\frac{2\cdot2^{\ast}_{\alpha}}{2\cdot2^{\ast}_{\alpha}-1}}dx\Bigg]dt\\
&=\int_{\Omega}\left(\int_{\Omega}\left(A\right)
\left(|u(y)|^{2^{\ast}_{\alpha}}|u(x)|^{2^{\ast}_{\alpha}-1}\right)^{\frac{2^{\ast}_{\alpha}}{2^{\ast}_{\alpha}-1}}dy\right)^{\frac{2\cdot2^{\ast}_{\alpha}}{2\cdot2^{\ast}_{\alpha}-1}}dx=\int_{\Omega}\left(\int_{\Omega}|G(u)|^{\frac{2^{\ast}_{\alpha}}{2^{\ast}_{\alpha}-1}}dy\right)^{\frac{2\cdot2^{\ast}_{\alpha}}{2\cdot2^{\ast}_{\alpha}-1}}dx.
\end{aligned}
\end{equation*}
\end{proof}

\begin{lem}\label{PS}
Let $c_{\lambda}<c_{0}$ and $(u_{n})_{n\in\mathbb{N}}\subset\mathcal{N}_{\lambda}$ be the Palais-Smale sequence at $c_{\lambda}$. Then $u_{n}\longrightarrow u_{0}\neq0$ in $W^{\alpha,2^{\ast}_{\alpha}}_{0}(\rm{curl};\Omega)$ along a subsequence, where $u_{0}$ is the nontrivial weak limit in Lemma \ref{ground}.
\end{lem}

\begin{proof}
Let $(u_{n})$ be a $(PS)_{c_{\lambda}}-$sequence such that $(u_{n})\subset\mathcal{N}_{\lambda}$. By Lemma \ref{coercive}, $(u_{n})$ is bounded and we can assume that $u_{n}\rightharpoonup u_{0}$ in $W^{\alpha,2^{\ast}_{\alpha}}_{0}(\mathrm{curl};\Omega)$. Then as in the proof of Lemma \ref{weak4}, we have $J^{'}_{\lambda}(u_{0})=0$, this implies that $u_{0}$ is a solution for (\ref{nonlocal case}). Moreover, by the concentration-compactness lemma, we have $u_{n}\longrightarrow u_{0}$ in $L^{2^{\ast}_{\alpha}}_{loc}(\Omega)$, see the same analysis in Lemma \ref{weak4}. On the other hand, By the compactly perturbed analysis in Lemma \ref{ground}, the weak limits $u_{0}\neq0$. Then by the general principle for the refined nonlocal Brezis-Lieb identity in \cite[Proposition 4.3 (ii) (iii)]{Mercuri2016}, we have
\begin{equation*}
\mathop{\mathrm{lim}}\limits_{n\longrightarrow\infty}\left(\int_{\Omega}(I_{\alpha}\ast|u_{n}|^{2^{\ast}_{\alpha}})|u_{n}|^{2^{\ast}_{\alpha}}dx-
\int_{\Omega}(I_{\alpha}\ast|u_{n}-u_{0}|^{2^{\ast}_{\alpha}})|u_{n}-u_{0}|^{2^{\ast}_{\alpha}}dx\right)\longrightarrow \int_{\Omega}(I_{\alpha}\ast|u_{0}|^{2^{\ast}_{\alpha}})|u_{0}|^{2^{\ast}_{\alpha}}dx,
\end{equation*}
we hence have
\begin{equation*}
\mathop{\mathrm{lim}}\limits_{n\longrightarrow\infty}\left(J_{\lambda}(u_{n})-J_{\lambda}(u_{n}-u_{0})\right)=J_{\lambda}(u_{0})\geq0,
\end{equation*}
and by Lemma \ref{Brezis}
\begin{equation*}
\mathop{\mathrm{lim}}\limits_{n\longrightarrow\infty}\left(J'_{\lambda}(u_{n})-J'_{\lambda}(u_{n}-u_{0})\right)=J'_{\lambda}(u_{0})=0.
\end{equation*}
Since $J'(u_{n})\longrightarrow 0$ and $u_{n}\longrightarrow u_{0}$ in $L^{2}(\Omega,\mathbb{R}^{3})$, we have
\begin{equation}\label{000}
\mathop{\mathrm{lim}}\limits_{n\longrightarrow\infty}J'_{0}(u_{n}-u_{0})=0.
\end{equation}

Suppose $\mathop{\mathrm{lim}~\mathrm{inf}}\limits_{n\longrightarrow\infty}||u_{n}-u_{0}||>0$. Since
$\mathop{\mathrm{lim}}\limits_{n\longrightarrow\infty}J'_{0}(u_{n}-u_{0})(u_{n}-u_{0})=0$, we infer that
\begin{equation*}\label{00}
\mathop{\mathrm{lim}~\mathrm{inf}}\limits_{n\longrightarrow\infty}|\nabla\times(u_{n}-u_{0})|_{2}>0.
\end{equation*}
Let $u_{n}-u_{0}=v_{n}+\widetilde{w}_{\Omega}(v_{n})\in\mathcal{V}_{\Omega}\oplus\mathcal{W}_{\Omega}$ according to the Helmholtz decomposition in $W^{\alpha,2^{\ast}_{\alpha}}_{0}(\mathrm{curl};\Omega)$. If $v_{n}\longrightarrow0$ in $Q^{\alpha,2^{\ast}_{\alpha}}(\Omega,\mathbb{R}^{3})$, then by (\ref{000}) we have $J'_{0}(u_{n}-u_{0})v_{n}\longrightarrow0$, thus
\begin{equation*}
|\nabla\times(u_{n}-u_{0})|^{2}_{2}=|\nabla\times v_{n}|^{2}_{2}=J'_{0}(u_{n}-u_{0})v_{n}+\int_{\Omega}\left\langle\left(I_{\alpha}\ast|u_{n}-u_{0}|^{2^{\ast}_{\alpha}}\right)|u_{n}-u_{0}|^{2^{\ast}_{\alpha}-2}(u_{n}-u_{0}),v_{n}\right\rangle dx\longrightarrow0
\end{equation*}
as $n\longrightarrow\infty$, which is a contradiction. Therefore $|v_{n}|_{Q^{\alpha,2^{\ast}_{\alpha}}}$ is bounded away from 0. If $w_{n}:=\widetilde{w}_{\Omega}(u_{n}-u_{0})\in\mathcal{W}_{\Omega}$, then $(w_{n})$ is bounded and since $u_{n}-u_{0}+w_{n}=v_{n}+\widetilde{w}_{\Omega}(v_{n})\in\mathcal{V}_{\Omega}\oplus\mathcal{W}_{\Omega}$, $|u_{n}-u_{0}+w_{n}|_{Q^{\alpha,2^{\ast}_{\alpha}}}$ is bounded away from 0. Choose $t_{n}$ so that $t_{n}(u_{n}-u_{0}+w_{n})\in \mathcal{N}_{\Omega}$, see (\ref{N Omega}). As in (\ref{tn1}) we have
\begin{equation*}\label{tn2}
t_{n}^{2}=\frac{\left(\int_{\Omega}|\nabla\times   (u_{n}-u_{0})|^{2}dx\right)^{\frac{1}{2^{\ast}_{\alpha}-1}}}
{\left(\int_{\Omega}|I_{\alpha/2}\ast|u_{n}-u_{0}+w_{n}|^{2^{\ast}_{\alpha}}|^{2}dx\right)^{\frac{1}{2^{\ast}_{\alpha}-1}}},
\end{equation*}
and so $(t_{n})$ is bounded. Then using Lemma \ref{tn}, we have
\begin{equation*}
J_{0}(u_{n}-u_{0})\geq J_{0}(t_{n}(u_{n}-u_{0}+w_{n}))-J'_{0}(u_{n}-u_{0})\left[\frac{t^{2}_{n}-1}{2}(u_{n}-u_{0}+t^{2}_{n}w_{n})\right],
\end{equation*}
so by (\ref{000}) and since $u_{n}\longrightarrow u_{0}$ in $L^{2}(\Omega,\mathbb{R}^{3})$,
\begin{equation*}
c_{\lambda}=\mathop{\mathrm{lin}}\limits_{n\longrightarrow\infty} J_{\lambda}(u_{n}-u_{0})=\mathop{\mathrm{lim}}\limits_{n\longrightarrow\infty}J_{0}(u_{n}-u_{0})\geq\mathop{\mathrm{lim}}\limits_{n\longrightarrow\infty}J_{0}(t_{n}(u_{n}-u_{0}+w_{n}))\geq c_{0},
\end{equation*}
which is a contradiction. Therefore, passing to a subsequence, $u_{n}\longrightarrow u_{0}$, hence also in the $\mathcal{T}-$topology.
%, which contains condition (A3).
 %Moreover, since $u_{0}\in\mathcal{N}_{\lambda}$, $u_{0}\neq0$.
\end{proof}

Finally, we shall compare $c_{\lambda}$ and $c_{0}$ in some ranges of $\lambda$. Recall from the third identity in (\ref{three}), we note that $c_{0}=\frac{2^{\ast}_{\alpha}-1}{2\cdot2^{\ast}_{\alpha}} \bar{S}^{\frac{2^{\ast}_{\alpha}}{2^{\ast}_{\alpha}-1}}_{\mathrm{curl},HL}(\Omega)\geq\frac{2^{\ast}_{\alpha}-1}{2\cdot2^{\ast}_{\alpha}}S_{HL}^{\frac{2^{\ast}_{\alpha}}{2^{\ast}_{\alpha}-1}}$. \begin{lem}
Let $\lambda\in(-\lambda_{\nu},-\lambda_{\nu-1}]$ for some $\nu\geq1$. There holds
\begin{equation*}
\begin{aligned}
&c_{\lambda}=\mathop{\mathrm{inf}}\limits_{\mathcal{N}_{\lambda}}J_{\lambda}\leq\frac{2^{\ast}_{\alpha}-1}{2\cdot2^{\ast}_{\alpha}}(\lambda+\lambda_{\nu})^{\frac{2^{\ast}_{\alpha}}{2^{\ast}_{\alpha}-1}}|\mathrm{diam}\Omega|^{\frac{3\cdot2^{\ast}_{\alpha}-\alpha-3}{2^{\ast}_{\alpha}-1}},\\ &c_{\lambda}<c_{0}~~if~~\lambda<-\lambda_{\nu}+\bar{S}_{\mathrm{curl},HL}(\Omega)|\mathrm{diam}\Omega|^{-\frac{3\cdot2^{\ast}_{\alpha}-\alpha-3}{2^{\ast}_{\alpha}}}. \end{aligned}
\end{equation*}

\end{lem}
\begin{proof}
Let $e_{\nu}$ be an eigenvector corresponding to $\lambda_{\nu}$. Then $e_{\nu}\in\mathcal{V}_{\Omega}^{+}$. Choose $t>0$, $\widetilde{v}\in\widetilde{\mathcal{V}}_{\Omega}$ and $w\in\mathcal{W}_{\Omega}$ so that $u=v+w=t e_{\nu}+\widetilde{\nu}+w\in\mathcal{N}_{\lambda}$. Since $\lambda_{k}\leq\lambda_{\nu}$ for $k<\nu$,
\begin{equation*}
\begin{aligned}
c_{\lambda}&\leq J_{\lambda}(te_{\nu}+\widetilde{v}+w)\\
&=\frac{\lambda_{\nu}}{2}\int_{\Omega}|te_{v}|^{2}dx+ \frac{1}{2}\int_{\Omega}|\nabla\times \widetilde{v}|^{2}dx+\frac{\lambda}{2}\int_{\Omega}|u|^{2}dx-\frac{1}{2\cdot2^{\ast}_{\alpha}}\int_{\Omega}|I_{\alpha/2}\ast|u|^{2^{\ast}_{\alpha}}|^{2}dx\\
&\leq\frac{\lambda_{\nu}}{2}\int_{\Omega}|v|^{2}dx+\frac{\lambda}{2}\int_{\Omega}|u|^{2}dx-\frac{1}{2\cdot2^{\ast}_{\alpha}}\int_{\Omega}|I_{\alpha/2}\ast|u|^{2^{\ast}_{\alpha}}|^{2}dx\\
&\leq\frac{\lambda+\lambda_{\nu}}{2}\int_{\Omega}|u|^{2}dx-\frac{1}{2\cdot2^{\ast}_{\alpha}}\int_{\Omega}|I_{\alpha/2}\ast|u|^{2^{\ast}_{\alpha}}|^{2}dx,\\
&\leq\frac{\lambda+\lambda_{\nu}}{2}\int_{\Omega}|u|^{2}dx-\frac{1}{2\cdot2^{\ast}_{\alpha}}\frac{1}{|\mathrm{diam}\Omega|^{3-\alpha}}\left(\int_{\Omega}|u|^{2^{\ast}_{\alpha}}dx\right)^{2},
\end{aligned}
\end{equation*}
%By the semigroup property of the Riesz potential (see \cite{Riesz1949}), we obtain
%\begin{equation*}
%\begin{aligned}
%c_{\lambda}&\leq\frac{\lambda+\lambda_{\nu}}{2}\int_{\Omega}|u|^{2}dx-\frac{1}{2\cdot2^{\ast}_{\alpha}}\int_{\Omega}\left(\int_{\Omega}\frac{|u(y)|^{2^{\ast}_{\alpha}}}{|x-y|^{\frac{N+\mu}{2}}}dy\right)^{2}dx\\
%&\leq\frac{\lambda+\lambda_{\nu}}{2}\int_{\Omega}|u|^{2}dx-\frac{1}{2\cdot2^{\ast}_{\alpha}}\int_{\Omega}\int_{\Omega}\frac{|u(x)|^{2\cdot2^{\ast}_{\alpha}}}{|x-y|^{\mu}}dxdy\\
%&\leq\frac{\lambda+\lambda_{\nu}}{2}\int_{\Omega}|u|^{2}dx-\frac{1}{2\cdot2^{\ast}_{\alpha}}\frac{1}{|\mathrm{diam}\Omega|^{\mu}}\left(\int_{\Omega}|u|^{2^{\ast}_{\alpha}}dx\right)^{2},\\
%\end{aligned}
%\end{equation*}
where $|\mathrm{diam}\Omega|=\mathop{max}\limits_{x,y\in\Omega}|x-y|$. Then using the H\"{o}lder inequality, we get
\begin{equation*}
\begin{aligned}
c_{\lambda}&\leq\frac{\lambda+\lambda_{\nu}}{2}\left[\left(\int_{\Omega}|u|^{2^{\ast}_{\alpha}}dx\right)^{\frac{1}{2^{\ast}_{\alpha}}}\right]^{2}|\Omega|^{\frac{2^{\ast}_{\alpha}-2}{2^{\ast}_{\alpha}}}
-\frac{1}{2\cdot2^{\ast}_{\alpha}}\frac{1}{|\mathrm{diam}\Omega|^{3-\alpha}}\left(\int_{\Omega}|u|^{2^{\ast}_{\alpha}}dx\right)^{2}\\
&\leq\frac{\lambda+\lambda_{\nu}}{2}\left[\left(\int_{\Omega}|u|^{2^{\ast}_{\alpha}}dx\right)^{\frac{1}{2^{\ast}_{\alpha}}}\right]^{2}|\mathrm{diam}\Omega|^{3\cdot\frac{2^{\ast}_{\alpha}-2}{2^{\ast}_{\alpha}}}
-\frac{1}{2\cdot2^{\ast}_{\alpha}}\frac{1}{|\mathrm{diam}\Omega|^{3-\alpha}}\left(\int_{\Omega}|u|^{2^{\ast}_{\alpha}}dx\right)^{2}\\
%&\leq\frac{\lambda+\lambda_{\nu}}{2}\left[\left(\int_{\Omega}|u|^{2^{\ast}_{\alpha}}dx\right)^{\frac{1}{2^{\ast}_{\alpha}}}\right]^{2}\left[|\mathrm{diam}\Omega|^{-\frac{\mu}{2\cdot2^{\ast}_{\alpha}}}\right]^{2}|\mathrm{diam}\Omega|^{3\cdot\frac{2^{\ast}_{\alpha}-2}{2^{\ast}_{\alpha}}+\frac{2\mu}{2\cdot2^{\ast}_{\alpha}}}
%-\frac{1}{2\cdot2^{\ast}_{\alpha}}\left[|\mathrm{diam}\Omega|^{-\frac{\mu}{2\cdot2^{\ast}_{\alpha}}}\right]^{2\cdot2^{\ast}_{\alpha}}\left[\left(\int_{\Omega}|u|^{2^{\ast}_{\alpha}}dx\right)^{\frac{1}{2^{\ast}_{\alpha}}}\right]^{2\cdot2^{\ast}_{\alpha}}\\
&\leq\frac{2^{\ast}_{\alpha}-1}{2\cdot2^{\ast}_{\alpha}}(\lambda+\lambda_{\nu})^{\frac{2^{\ast}_{\alpha}}{2^{\ast}_{\alpha}-1}}|\mathrm{diam}\Omega|^{\frac{3\cdot2^{\ast}_{\alpha}-\alpha-3}{2^{\ast}_{\alpha}-1}},
\end{aligned}
\end{equation*}
where the last inequality follows from the inequality  $\frac{A}{2}t^{2}-\frac{1}{p}t^{p}\leq(\frac{1}{2}-\frac{1}{p})A^{\frac{p}{p-2}}~(A>0)$.
%$A=(\lambda+\lambda_{\nu})|\mathrm{diam}~\Omega|^{\frac{2\cdot2^{\ast}_{\alpha}-\alpha+5}{2^{\ast}_{\alpha}}}$ %and %$t=\left(\int_{\Omega}|u|^{2^{\ast}_{\alpha}}dx\right)^{\frac{1}{2^{\ast}_{\alpha}}}\cdot|\mathrm{diam}~\Omega|^{\frac{\alpha-3}{2\cdot2^{\ast}_{\alpha}}}$.

Since $c_{0}=\frac{2^{\ast}_{\alpha}-1}{2\cdot2^{\ast}_{\alpha}} \bar{S}^{\frac{2^{\ast}_{\alpha}}{2^{\ast}_{\alpha}-1}}_{\mathrm{curl},HL}(\Omega)$, the second inequality follows immediately.
\end{proof}

\begin{proof}[\bf  Complete of the Proof of Theorem \ref{main thm}] Note that if $\lambda<-\lambda_{\nu}+\bar{S}_{\mathrm{curl},HL}(\Omega)|\mathrm{diam}\Omega|^{-\frac{3\cdot2^{\ast}_{\alpha}-\alpha-3}{2^{\ast}_{\alpha}}}$, then $c_{\lambda}<c_{0}$, and by Lemma \ref{PS}, $J_{\lambda}$ satisfies the $(PS)_{c_{\lambda}}$ condition, hence satisfies the $(PS)_{c_{\lambda}}^{\mathcal{T}}$ condition. Then statement (a) follows from Lemma \ref{ground},
%Actually, since $u_{n}\longrightarrow u$ in $L^{2}(\Omega,\mathbb{R}^{3})$, we may use the compactly perturb method (see \cite{Mederski2018}) and obtained the same conclusion by setting $J_{cp}=J_{0}$.
and the remaining statements (b)-(d) are similar to \cite[Theorem 1.4]{Mederski2021} and can be proved by the same strategy.
\end{proof}

\subsection*{Acknowledgements}
Minbo Yang was supported by National Natural Science Foundation of China (No. 11971436, No. 12011530199) and Natural Science Foundation of Zhejiang Province (No. LZ22A010001).

\end{document}